\documentclass{article}
\usepackage[utf8]{inputenc}
\usepackage{comment}
\usepackage{amssymb,amsmath,amsthm,graphicx,mathtools,enumerate,mathrsfs}
\usepackage{fullpage,enumitem}
\usepackage{thmtools}
\usepackage{thm-restate}

\makeatletter
  \@ifundefined{newcounteralias}{}{%
    \renewcommand\thmt@autorefsetup{%
      \@xa\def\csname\thmt@envname autorefname\@xa\endcsname
        \@xa{\thmt@thmname}%
    }%
  }
  \makeatother

\usepackage{appendix}
\usepackage{comment}

\usepackage[a4paper,margin=2cm]{geometry}

\usepackage{setspace}
\usepackage[dvipsnames]{xcolor}

\usepackage[setpagesize=false, colorlinks=true, linkcolor=BrickRed, citecolor=OliveGreen, urlcolor=black, pdfencoding=auto, psdextra,]{hyperref}
\usepackage{cleveref}

\usepackage{tikz}
\usetikzlibrary{math}
\usetikzlibrary{snakes,calc,decorations.pathmorphing,through}
\tikzset{snake it/.style={decorate, decoration=snake}}

\theoremstyle{plain}
\newtheorem{theorem}{Theorem}[section]
\crefname{theorem}{Theorem}{Theorems}

\crefname{proposition}{Proposition}{Propositions}

\crefname{corollary}{Corollary}{Corollaries}

\newtheorem{lemma}[theorem]{Lemma}
\crefname{lemma}{Lemma}{Lemmas}

\newtheorem{conjecture}[theorem]{Conjecture}
\crefname{conjecture}{Conjecture}{Conjectures}

\crefname{problem}{Problem}{Problem}

\newtheorem{claim}{Claim}
\crefname{claim}{Claim}{Claims}

\crefname{observation}{Observation}{Observations}

\crefname{setup}{Setup}{Setups}

\crefname{fact}{Fact}{Facts}

\crefname{algorithm}{Algorithm}{Algorithms}

\crefname{remark}{Remark}{Remarks}

\crefname{example}{Example}{Examples}

\theoremstyle{definition}
\newtheorem{definition}[theorem]{Definition}
\crefname{definition}{Definition}{Definitions}

\crefname{construction}{Construction}{Constructions}

\crefname{question}{Question}{Questions}

\def\eps{\varepsilon}

\makeatletter
\newenvironment{claimproof}[1][Proof]{\par
	\pushQED{\qed}%
	
	\normalfont \topsep6\p@\@plus6\p@\relax
	\trivlist
	\item[\hskip\labelsep
	\textit{#1}\@addpunct{.}~]\ignorespaces
}{%
	\popQED\endtrivlist\@endpefalse
}
\makeatother

\usepackage{tikz}
\usetikzlibrary{math}
\usetikzlibrary{calc,decorations.pathmorphing,through}
\tikzset{snake it/.style={decorate, decoration=snake}}

\definecolor{DarkDesaturatedBlue}{HTML}{3A3556}
\definecolor{VividOrange}{HTML}{F15918}
\definecolor{PureOrange}{HTML}{FFBA00}
\definecolor{LightGrayishPink}{HTML}{EEC5D5}
\definecolor{VerySoftBlue}{HTML}{B5AFDB}

\newcounter{propcounter}
\stepcounter{propcounter}

\title{Hamilton cycles in pseudorandom graphs: resilience and approximate decompositions
}
\author{
Nemanja Dragani\'c\thanks{
Mathematical Institute, University of Oxford, UK. Research supported by SNSF project 217926.\\
\emph{E-mail}: \textbf{nemanja.draganic@maths.ox.ac.uk}} \and 
Jaehoon Kim\thanks{Department of Mathematical Sciences, KAIST, South Korea.
        \emph{E-mails:} \textbf{\{jaehoon.kim, hyunwoo.lee\}@kaist.ac.kr}. Supported by the National Research Foundation of Korea (NRF) grant funded by the Korea government(MSIT) No. RS-2023-00210430.}
\and Hyunwoo Lee\footnotemark[2]~\thanks{Extremal Combinatorics and Probability Group (ECOPRO), Institute for Basic Science (IBS). Supported by the Institute for Basic Science (IBS-R029-C4).}
\and David Munh\'a Correia \thanks{
Department of Mathematics, ETH, Z\"urich, Switzerland. Research supported in part by SNSF grant 200021-228014.
\emph{E-mails}: \textbf{\{david.munhacanascorreia,benjamin.sudakov\}@math.ethz.ch}.
}
\and Mat\'ias Pavez-Sign\'e\thanks{Departamento de Ingenier\'ia Matem\'atica, Universidad de Chile, and Centro de Modelamiento Matem\'atico (CNRS IRL2807), Chile. Supported by ANID Fondecyt Regular grant No. 1241398 and ANID Basal Grant CMM FB210005.  \emph{E-mail:} \textbf{\ mpavez@dim.uchile.cl}.}
\and Benny Sudakov\footnotemark[4]
}
\date{}

\begin{document}

\maketitle

\begin{abstract}
Dirac's classical theorem asserts that, for $n\ge 3$, any $n$-vertex graph with minimum degree at least $\frac{n}{2}$ is Hamiltonian.  Furthermore, if we additionally assume that such graphs are regular and sufficiently large, then, by the breakthrough work of Csaba, Kühn, Lo, Osthus, and Treglown, they admit a decomposition into Hamilton cycles and at most one perfect matching, solving the well-known Nash-Williams conjecture.
In the pseudorandom setting, it has long been conjectured that similar results hold in much sparser graphs.
We prove two overarching theorems for graphs that exclude excessively dense subgraphs, which yield asymptotically optimal resilience and Hamilton-decomposition results in sparse pseudorandom graphs.
First, our results imply that for every fixed $\gamma>0$, there exists a constant $C>0$ such that if $G$ is a spanning subgraph of an $(n,d,\lambda)$-graph satisfying
$\delta(G)\ge\bigl(\tfrac12+\gamma\bigr)d$ and $
d/\lambda\ge C,$
then $G$ must contain a Hamilton cycle.

Second, we show that for every $\varepsilon>0$, there is $C>0$ so that any $(n,d,\lambda)$-graph $G$ with $d/\lambda\ge C$ contains at least $\bigl(1-\varepsilon\bigr)\tfrac d2$ edge-disjoint Hamilton cycles, and, finally, we prove that the entire edge set of $G$ can be covered by no more than $\bigl(1+\varepsilon\bigr)\tfrac d2$ such cycles.
All bounds are asymptotically optimal and significantly improve earlier results on Hamiltonian resilience, packing, and covering in sparse pseudorandom graphs.
\end{abstract}

% \tableofcontents 

%--------------------------------------------------------

\section{Introduction}
A Hamilton cycle in a graph $G$ is a cycle passing through all vertices of $G$, and we say $G$ is Hamiltonian if it contains one. Since deciding whether a graph is Hamiltonian is NP-complete,  finding sufficient conditions that guarantee the existence of a Hamilton cycle in a graph is one of the fundamental problems in graph theory. In this direction, one of the most influential results is Dirac's theorem~\cite{Dirac} from 1952, which states that every graph $G$ on $n\ge 3$ vertices with minimum degree at least $n/2$ contains a Hamilton cycle (see~\cite{Gould-survey} and references therein for further sufficient conditions forcing Hamiltonicity). 
Most such conditions, however, apply only to very dense graphs, as structural properties are typically easier to uncover in dense settings. Therefore, developing analogous conditions for sparse graphs remains a major and challenging area of research.

Random graph models offer a natural and extensively explored framework for analysing sparse graphs. The appearance of Hamilton cycles in random graphs is a well-studied problem and dates back to 1976, when  P\'{o}sa~\cite{Posa} proved that the binomial random graph $G(n, p)$ is Hamiltonian with high probability once $p\geq C\log n/n$ for a sufficiently large constant $C$, introducing along the way the influential \textit{rotation-extension} technique for finding long cycles. This was later sharpened by Koml\'os and Szemer\'edi~\cite{KOMLOS198355} and by Korshunov~\cite{Kor76}, who showed that $p=(\log n+\log\log n+\omega(1))/n$ already suffices, which is precisely the threshold at which $G(n,p)$ attains minimum degree 2. Indeed, Bollob\'as~\cite{bollobas_randomreg} and Ajtai, Koml\'os and Szemer\'edi~\cite{AJTAI1985173} showed that the \textit{random graph process} becomes Hamiltonian as soon as it reaches minimum degree 2. For random regular graphs, Robinson and Wormald~\cite{RobinsonWormald} proved that the random $d$-regular graph is Hamiltonian with high probability for all fixed $d\ge 3$, improving upon earlier results of Bollob\'as~\cite{BOLLOBAS198397} and Fenner and Frieze~\cite{FENNER1984103}. The result of Robinson and Wormald was later extended to all $d=d(n)\geq 3$ in a series of works (see~\cite{ cooper2002random,krivelevich2001random}).

With Hamiltonicity in random graphs now well understood, our next goal is to tackle deterministic graph classes---first and foremost, pseudorandom graphs.
 Following Alon's definition, we say that $G$ is an $(n,d,\lambda)$-graph if $G$ is a $d$-regular graph on $n$ vertices such that all non-trivial eigenvalues of its adjacency matrix are bounded (in absolute value) by $\lambda$. This class has become a standard
model for pseudorandom graphs, as the \textit{Expander mixing lemma} implies that the edges of an $(n,d,\lambda)$-graph are evenly distributed among vertex subsets (see~\cite{krivelevich2006pseudo} for a comprehensive survey on pseudorandom graphs). Moreover, for each fixed $d\geq3$, random $d$-regular graphs are with high probability $(n,d,\lambda)$-graphs for $\lambda=O(\sqrt{d})$ (see e.g.~\cite{Spencer08}). Since random $d$-regular graphs are Hamiltonian for every fixed $d\ge 3$, it is natural to expect that $(n,d,\lambda)$-graphs are also Hamiltonian, even for $d$ constant. This was made precise by Krivelevich and Sudakov~\cite{Krivelevich-Sudakov}, who conjectured in 2003 that $(n, d, \lambda)$-graphs are Hamiltonian when the spectral ratio satisfies $d/\lambda\ge C$ for some absolute constant $C>0$. After significant efforts of several groups (see~\cite{Allen,Ferber,Glock,Hefetz,Krivelevich}), this conjecture was finally resolved by Dragani\'{c}, Montgomery, Munh\'{a} Correia, Pokrovskiy, and Sudakov~\cite{C-expander}.

\begin{theorem}[\cite{C-expander}]\label{thm:ndlambda-hamilton}
    There exists a constant $C > 0$ such that every $(n, d, \lambda)$-graph with $d/\lambda\ge C$ is Hamiltonian.
\end{theorem}
In fact, the authors of~\cite{C-expander} proved an even stronger result, showing that the key property of $(n,d,\lambda)$-graphs required for Hamiltonicity is that they are good \emph{expanders}, as defined below. An $ n $-vertex graph $ G $ with $n\geq3$ is called a \emph{$ C $-expander} if $|N_G(X)| \geq C|X|$
for every subset $ X \subset V(G) $ with $ |X| \leq n/(2C) $, and if there is an edge between any two disjoint sets $ X, Y \subset V(G) $ of size at least $ n/(2C) $.\footnote{For a graph $ G $ and a set $ X \subset V(G) $, the \emph{external neighbourhood} of $ X $, denoted $ N_G(X) $, is the set of vertices in $ V(G) \setminus X $ adjacent to vertices in $ X $.}
\begin{theorem}[\cite{C-expander}]\label{thm:C-expander}
   There is a constant $C>0$ such that every $C$-expander graph is Hamiltonian. Moreover, it is \emph{Hamilton-connected}, that is, any two of its vertices are the endpoints of a Hamilton path.
\end{theorem}
Since the Expander mixing lemma also shows that $(n,d,\lambda)$-graphs are $C$-expanders for $C=\Omega(d/\lambda)$, Theorem~\ref{thm:ndlambda-hamilton} follows immediately from Theorem~\ref{thm:C-expander}.

\subsection{Resilience}

In many graph-theoretic problems, it is not enough to show that a desirable property holds; one also wants to know how robust it is under adversarial changes. This viewpoint is captured by the notion of \emph{resilience}, which measures how much a graph can be altered---typically via edge or vertex deletions---while still retaining a given property. For instance, one may ask how many edges can be removed from a Hamiltonian graph before it ceases to contain any Hamilton cycle, or whether every spanning subgraph of sufficiently large minimum degree remains Hamiltonian. The systematic study of these questions was initiated by 
Sudakov and Vu~\cite{Local-resilience}, see also the survey \cite{Sudakov} for a thorough treatment.

Since establishing resilience typically requires significantly more intricate arguments than proving the existence of the property itself, results on the resilience of Hamiltonicity were proven only much later than the foundational results on Hamiltonicity.
For binomial random graphs, Lee and Sudakov~\cite{Choongbum-Lee} obtained an optimal version of Dirac's theorem. They showed that for every fixed $\varepsilon>0$ there is $C=C(\varepsilon)$ such that if $p\geq C\log n/n$, then with high probability every spanning subgraph of $G(n,p)$ with minimum degree at least $(1/2+\varepsilon)pn$ is Hamiltonian. For random regular graphs, Ben-Shimon, Krivelevich and Sudakov~\cite{Random-regular} showed that, given $\varepsilon>0$ and $d\ge d_0(\varepsilon)$ fixed, with high probability every spanning subgraph of a random $d$-regular graph with minimum degree at least $(1+\varepsilon)5d/6$ is Hamiltonian, and conjectured that the constant $5/6$ might be improved to $1/2$, which was later confirmed by Condon, Espuny D\'{i}az, Girão, Kühn, and Osthus~\cite{condon2021dirac}.

For pseudorandom graphs, Sudakov and Vu~\cite{Local-resilience} proved an optimal version of Dirac's theorem under the condition $d/\lambda \ge \log^{2} n$, namely, any spanning subgraph of an $(n,d,\lambda)$-graph with minimum degree at least $\left(1+o(1)\right)\tfrac{d}{2}$ contains a Hamilton cycle.
As in the conjecture of Krivelevich and Sudakov, the goal in this line of research is to show that Dirac's theorem holds for all suitable pseudorandom graphs, in particular for $(n,d,\lambda)$-graphs even when $d/\lambda $ is just at least a large constant, which implies Hamilton resilience for random $d$-regular graphs for large enough $d$.

Before stating our first result, let us note that strong expansion alone (as in Theorem~\ref{thm:C-expander}) is not sufficient for an analogue of Dirac's theorem. 
Let $G$ be a graph on a vertex set $V_1\cup V_2$ with $|V_1|=|V_2|=n/2$, where $G[V_1]$ and $G[V_2]$ are cliques, and $G[V_1, V_2]$ is a random $d$-regular bipartite graph for any $d<n/4$, say, with $d$ sufficiently large. Clearly, removing all edges in $G[V_1,V_2]$ disconnects the graph, thus the remaining graph is not Hamiltonian, although the initial graph is an excellent expander. The main issue with this construction is the presence of large dense parts. The next example shows that dense subgraphs on small vertex sets can also prevent Hamilton resilience.

Let $0<\gamma \ll 1$ and let $d$ be sufficiently large. Consider a random $2\gamma d$-regular graph $F$ on $m$ vertices with $m> 2\gamma d$ and a random $(1-2\gamma)d$-regular graph $H$ on $2n$ vertices, where $n$ is a sufficiently large multiple of $m$. Note that, with high probability, $H$ is an excellent expander and has no dense spots---essentially an optimal graph in that regard. Let $V(H)=V_1\cup V_2$ be a balanced partition of $H$ chosen so that $|N(v)\cap V_i|\geq (1/2-\gamma-O(d^{-1/2}\log d))d$ for each vertex $v$ and $i\in[2]$, which can be found by a standard application of the \textit{Lov\'asz local lemma}.
Now, take $n/m$ vertex-disjoint copies of $F$ and put them in each $V_i$ so that the resulting graph $G_0$ is $d$-regular (indeed, it is not hard to see that such a packing must exist by the Sauer--Spencer theorem~\cite{Sauer-Spencer}). Again, $G_0[V_1]\cup G_0[V_2]$ is a subgraph of $G_0$ with minimum degree at least $(1/2+\gamma - O(d^{-1/2}\log d))d$ but it is not even connected. This shows that even mildly denser spots in excellent expanders can destroy Hamilton resilience.

Similarly, using an appropriate $2m$-vertex random regular bipartite graph $F$ and appropriate $H$, one can also construct a strong expander $G_1$ on a vertex set $V_1\cup V_2$ with $|V_1|=|V_2|+1$ and $\delta(G_1[V_1,V_2])\geq (1/2+\gamma - O(d^{-1/2}\log d))d$. These examples suggest that we need a sparsity condition to rule out local dense structures between sets of size $m=\Omega(d)$ (like in the graph $F$ above). For this, the following definition is natural. 
\begin{definition}\label{def:sparse}
Say an $n$-vertex graph $G$ is \emph{$(\eta,\beta,p)$-sparse} if for all (not necessarily disjoint) sets $U,W\subset V(G)$ with $\eta pn \leq |U|\le |W|$, $$e(U,W) \leq (1+\beta) p\max\{\eta n,|U|\}|W|,$$ where $e(U,W)$ is the number of pairs $(u,w)\in U\times W$ with $uw\in E(G)$.
\end{definition}

This notion covers a wide class of graphs, including $(n,d,\lambda)$-graphs with mild conditions on $d/\lambda$ and the binomial random graph model. One can view $(\eta,\beta,p)$-sparsity as essentially the upper-bound condition that follows from the {Expander mixing lemma} (Theorem~2.11 in~\cite{krivelevich2006pseudo}), a condition that is satisfied by $(n,d,\lambda)$-graphs and the binomial random graph. Intuitively, the edge-density between subsets of vertices is not much larger than what we would expect in a random graph with edge probability $p$. This added flexibility allows our results to apply to various random graph models and a broader family of pseudorandom graphs rather than only to $(n,d,\lambda)$-graphs. This definition serves as a generalization of $(n,d,\lambda)$-graphs in our setting, mirroring the strong generalization of the Krivelevich--Sudakov conjecture to $C$-expanders in \Cref{thm:C-expander}.

We now state our first main result, the resilience of Hamiltonicity in $(\eta,\beta,p)$-sparse graphs.
\begin{theorem}\label{thm:sparse_mindegree_cycle}
For every $0<\gamma< 1/2$, there exist constants $\beta,\eta>0$ and $d_0\in\mathbb N$ such that the following holds. Let $n,d\in \mathbb{N}$ satisfy $n>d\ge d_0$ and suppose $G$ is an $n$-vertex $(\eta,\beta,d/n)$-sparse graph with $\delta(G)\geq (1/2+\gamma)d$. Then $G$ contains a Hamilton cycle.
\end{theorem}
As discussed, the Expander mixing lemma shows that all 
$(n,d,\lambda)$-graphs with $\lambda/d \leq \beta \eta$ are $(\eta,\beta,d/n)$-sparse (see Lemma~\ref{lemma:nd:sparse}). Therefore, \Cref{thm:sparse_mindegree_cycle}  implies an asymptotically optimal resilience result for $(n, d, \lambda)$-graphs.
\begin{theorem}\label{thm:main}
    For every $\gamma>0$, there is a constant $C> 0$ such that the following holds. Every spanning subgraph $G$ of an $(n,d,\lambda)$-graph with ${d}/{\lambda} \ge C$ and $\delta(G) \geq \left(1/2 + \gamma \right)d$ contains a Hamilton cycle.
\end{theorem}
We note that our proof requires $C$ to be at least a tower-type function of $1/\gamma$. On the other hand, one can use  Weyl's inequality (see Section 4.3 in~\cite{HJ}) to check that the example $G_0$  before Definition~\ref{def:sparse} is a $(2n,d,(2+o(1))\gamma d)$-graph. With $2\gamma$ in place of $\gamma$, this example shows that $C=\Omega(1/\gamma)$ is necessary. It is an interesting open question to determine the optimal dependency.
%-------------------------------------------------------------------------------

\subsection{Packing and covering}

The problem of decomposing regular structures into edge-disjoint copies of smaller substructures has attracted considerable attention across various areas of combinatorics. Classical examples include the folklore result that complete graphs of even order admit decompositions into perfect matchings, as well as the result that regular bipartite graphs admit such decompositions---so-called $1$-factorizations (see~\cite[Corollary~2.1.3]{Diestel}). Motivated by this, numerous extensions and generalizations of $1$-factorizations have been studied; see, for example~\cite{Chetwynd-Hilton,1-factorization,Perkovic-Reed}, including recent progress in the setting of pseudorandom graphs~\cite{ferber2024regular}.

Apart from $1$-factorization, in design theory, the well-known Kirkman's schoolgirl problem posed in 1850 asks for the existence of a $K_3$-factorization of $K_n$ for $n \equiv 3 \pmod{6}$, which was settled by Ray-Chaudhuri and Wilson~\cite{Chaudhuri-Wilson} in 1968. One generalization of Kirkman's schoolgirl problem is the Oberwolfach problem raised by Ringel~\cite{Guy}, which asks whether every odd order complete graph can be completely decomposed into edge-disjoint copies of any fixed $2$-factor. This problem was solved for all sufficiently large graphs by Glock, Joos, Kim, K\"{u}hn, and Osthus~\cite{Oberwolfach} and further generalized to directed graphs by Keevash and Staden~\cite{Keevash-Staden}. For more results on (hyper)graph decomposition problems, see~\cite{Decomposition-survey}.

Returning to Hamilton cycles, the problem of decomposing a graph into Hamilton cycles has a very long history. One of the earliest results is Walecki's 1890  construction of a Hamilton decomposition of complete graphs of odd order (see~\cite{Walecki}). This result was extended to many general families by various researchers and, specifically, Csaba, K\"{u}hn, Lo, Osthus, and Treglown~\cite{1-factorization} proved that, for all even $d$ with $d \geq n/2$, every sufficiently large $d$-regular graph has a Hamilton decomposition. Also, answering an old question of Kelly (see~\cite{Kelly-conjecture}), K\"{u}hn and Osthus~\cite{Kelly-proof} showed that every sufficiently large regular tournament has a Hamilton decomposition, where a Hamilton cycle in a directed graph is a Hamilton cycle with a consistent orientation of its edges. See also~\cite{regular-expander} and the survey~\cite{Hamdecomp-survey} for an overview of various related topics.

The Hamilton decompositions and coverings in random graphs have also been extensively studied. After Pósa's seminal result on the existence threshold for Hamilton cycles in $G(n,p)$, Bollob\'{a}s and Frieze~\cite{bollobas1985matchings} initiated the study of Hamilton decomposition of random graphs, which culminated in the optimal result that, for $p\geq C\log n/n$ with $C$ sufficiently large, with high probability $G(n,p)$ contains $\lfloor \delta(G(n,p))/2 \rfloor$ edge-disjoint Hamilton cycles (see e.g.~\cite{Knox, Knox2, krivelevich2012optimal}). More recently, attention has shifted to covering problems: Glebov, Krivelevich, and Szab\'{o}~\cite{Glebov2013} initiated the study of how few Hamilton cycles suffice to cover all the edges of $G(n, p)$, and asked for which $p$ the theoretically best possible number---$\lceil \Delta(G(n,p))/2 \rceil$---suffices. This was answered recently in~\cite{draganic2025optimal}, closing a long line of research and establishing that at the weak threshold for Hamiltonicity, a Hamilton cover of optimal size also exists. Furthermore, for every fixed even $d \geq 4$, Kim and Wormald~\cite{Kim-Wormald} proved that a random $d$-regular graph has a Hamilton decomposition with high probability.

This problem has also been studied in the context of pseudorandom graphs. Ferber, Kronenberg, and Long~\cite{ferber2017packing} proved that in a class of almost regular pseudorandom graphs with average degree $d=\omega(\log^{14} n)$, one can find $(1 - \varepsilon)d/2$ edge-disjoint Hamilton cycles. Their result relies on a rather strong edge distribution property, which does not necessarily hold in $(n,d,\lambda)$-graphs---even when $\lambda$ is as small as possible. While much weaker pseudorandom conditions were known to imply Hamiltonicity, this remained the state of the art for approximate Hamilton cycle decompositions in pseudorandom graphs. In this paper, we show that $(n, d, \lambda)$-graphs with a fairly mild assumption on the spectral ratio contain $(1-\varepsilon)d/2$ edge-disjoint Hamilton cycles. In fact, we prove that significantly weaker conditions---namely, $(\eta,\beta,p)$-sparsity and almost regularity\footnote{For a graph $G$ and a vertex $x\in V(G)$, we write $d_G(x)=(1\pm \gamma)d$ if $(1-\gamma)d\le d_G(x)\le (1+\gamma)d.$}---are already sufficient.

\begin{theorem}\label{thm: decomp}
    For every $\varepsilon>0$, there exist constants $\beta,\eta,\gamma>0$ and $d_0\in\mathbb N$ such that the following holds. Let $n,d\in\mathbb N$ satisfy $n>d\ge d_0$ and suppose $G$ is an $n$-vertex $(\eta,\beta, d/n)$-sparse graph with $d_G(v)=(1\pm\gamma) d$ for all $v\in V(G)$. Then, $G$ contains at least $(1-\eps)d/2$ edge-disjoint Hamilton cycles.
\end{theorem}

Finally, we also prove that in such graphs only $(1+\varepsilon)d/2$ Hamilton cycles are required to cover every edge, hence we also resolve the dual covering question asymptotically. 

\begin{theorem}\label{thm: covering}
    For every $\varepsilon>0$, there exist constants $\beta,\gamma, \eta>0$ and $d_0\in\mathbb N$ such that the following holds. Let $n,d\in\mathbb N$ satisfy $n>d\ge d_0$ and suppose $G$ is an $n$-vertex $(\eta,\beta, d/n)$-sparse graph with $d_G(v)=(1\pm\gamma) d$ for all $v\in V(G)$. Then the edges of $G$ can be covered by at most $(1+\varepsilon)d/2$ Hamilton cycles.
\end{theorem}

As discussed, this implies that $(n,d,\lambda)$-graphs also admit approximate packings and covers.

\begin{theorem}\label{thm:decompo:ndlambda}
    For every $\varepsilon>0$, there exists a constant $C>0$ such that the following holds. If $G$ is an $(n,d,\lambda)$-graph with $d/\lambda \geq C$, then $G$ contains at least $(1-\varepsilon)d/2$ edge-disjoint Hamilton cycles. Moreover, the edges of $G$ can be covered with at most $(1+\varepsilon)d/2$ Hamilton cycles. 
\end{theorem}
Since random regular graphs typically have a Hamilton decomposition, we believe that Theorem~\ref{thm:decompo:ndlambda} may be improved, that is, we conjecture that there is a universal constant $C > 0$ such that every $(n, d, \lambda)$-graph with even $d$ and $d/\lambda \ge C$ has a Hamilton decomposition.  
%--------------------------------------------------------
\subsection{Applications to Cayley graphs}
Given a finite group $\Gamma$ and a symmetric set $S\subset \Gamma \setminus \{1_{\Gamma}\}$, the \textit{Cayley graph} $\text{Cay}(\Gamma,S)$ is the graph with vertex set $\Gamma$ where two distinct elements $x,y\in \Gamma$ form an edge whenever $xy^{-1}\in S$. The Cayley graph $\text{Cay}(\Gamma,S)$ is $|S|$-regular and vertex-transitive, and is connected if $S$ generates $\Gamma$. A folklore conjecture, stated for Cayley graphs by Rapaport-Strasser~\cite{rapaportstrasser1959}, says that connectivity suffices to force Hamiltonicity in Cayley graphs. It is closely related to Lov\'asz's question about Hamilton paths in connected vertex-transitive graphs~\cite{lovasz1970problem}.
\begin{conjecture}\label{conj:lovasz}Every connected Cayley graph on at least 3 vertices is Hamiltonian.
\end{conjecture}
Conjecture~\ref{conj:lovasz} is still open, and only partial results are known (see e.g.~\cite{CurranGallian96,WitteGallian84}). However, by a result of Pak and Radoi\v{c}i\'c~\cite{PakRadoicic09}, we know that every finite group $\Gamma$ of order at least three admits a generating set $S\subset \Gamma$ of size $O(\log |\Gamma|)$ such that $\text{Cay}(\Gamma,S)$ is Hamiltonian. Moreover, they conjectured that a random set of generators of size $O(\log |\Gamma|)$ is sufficient to force Hamiltonicity. This was proved in~\cite{C-expander}, using a result of Alon and Roichman~\cite{alonroichman1994}. 

For the next two applications, set $n=|\Gamma|$ and let $$\mathcal{O}=\{ \{g,g^{-1}\}: g\in \Gamma\setminus \{1_{\Gamma}\} \}$$
be the set of non-trivial inverse-orbits. Choose $k=C\log n$ inverse-orbits uniformly at random without replacement, and let $S$ be their union.
Since each inverse-orbit has one or two elements, we have $k\leq|S|\leq2k$. For every fixed $\theta>0$, if $C$ is sufficiently large, then with high probability $\text{Cay}(\Gamma,S)$ is an $(n,|S|,\theta|S|)$-graph. To see this, consider the adjacency matrix $B_O$ of $\text{Cay}(\Gamma,O)$ for each $O\in\mathcal O$. For independent uniform choices $O_1,\dots,O_k$, we have $\mathbb E[B_{O_i}]=(J-I)/|\mathcal O|$ and $\mathbb E[B_{O_i}|_{\mathbf1^\perp}]=-I/|\mathcal O|$, where $J$ is the all-ones matrix and $\mathbf1$ is the all-ones vector.
Also, the centred matrices $B_{O_i}|_{\mathbf1^\perp}+I/|\mathcal O|$ have norm at most $3$, since $\|B_O\|\leq2$. Since $k/|\mathcal O|=o(k)$, the matrix Hoeffding inequality~\cite[Theorem~1.3]{Tropp2012} gives
\[\mathbb P\left(\left\|\textstyle\sum_{i\in[k]}B_{O_i}|_{\mathbf1^\perp}\right\|>\theta k\right)=o(1)\]
as $C$ is sufficiently large.
Moreover, the probability of a repeated inverse-orbit is at most $\binom{k}{2}/|\mathcal O|=o(1)$. Conditional on no repetition, $\{O_1,\dots,O_k\}$ is a uniformly random set of $k$ inverse-orbits and $\sum_{i\in[k]}B_{O_i}$ is the adjacency matrix of $\text{Cay}(\Gamma,O_1\cup\dots\cup O_k)$. Hence the same spectral estimate holds for our model, and $\theta k\leq\theta|S|$ gives the assertion.
Thus Theorem~\ref{thm:main} yields the following application.

\begin{theorem}\label{thm:resilience:Cayley}For every $\eps>0$, there exists a constant $C>0$  such that the following holds. If $\Gamma$ is an $n$-element group and $S\subset \Gamma$ is the union of $C\log n$ random inverse-orbits chosen as above, then with high probability every spanning subgraph $G\subset \text{Cay}(\Gamma,S)$ with $\delta(G)\ge (1/2+\eps)|S|$ contains a Hamilton cycle.     
\end{theorem}

Note that the random model based on inverse-orbits is natural, as $g$ and $g^{-1}$ give the same edges of the Cayley graph. However, the analogous resilience conclusion can fail for a uniform random set of group elements. To see this, let $\Gamma$ be the dihedral group of order $n$ and let $R$ be its rotation subgroup of order $n/2$. Choose a uniform random subset $S\subset\Gamma\setminus\{1_\Gamma\}$ of size $k=C\log n$. With high probability, the simple Cayley graph $G=\text{Cay}(\Gamma,S\cup S^{-1})$ has degree $d=(3/2+o(1))k$. Deleting the edges between $R$ and $\Gamma\setminus R$ leaves a disconnected spanning subgraph of degree $(1+o(1))k=(2/3+o(1))d$. Thus the conclusion of Theorem~\ref{thm:resilience:Cayley} fails in this model for every fixed $0<\varepsilon<1/6$.

Regarding decomposition into Hamilton cycles, a conjecture of Alspach~\cite{alspach1985problem} states that connected abelian Cayley graphs admit a Hamilton decomposition when the degree is even. As a consequence of Theorem~\ref{thm:decompo:ndlambda} and the above spectral computation, random Cayley graphs admit an approximate decomposition into Hamilton cycles.
\begin{theorem}For every $\varepsilon>0$, there exists a constant $C>0$ such that the following holds. If $\Gamma$ is an $n$-element group and $S\subset \Gamma$ is the union of $C\log n$ random inverse-orbits chosen as above, then with high probability $\text{Cay}(\Gamma,S)$ contains at least $(1-\varepsilon)|S|/2$ edge-disjoint Hamilton cycles.\end{theorem}
\subsection{Brief proof overview and organization}\label{sec:overview}
Here we give a high-level overview of the proof strategies. Each later section also starts with a short overview.
\paragraph{Resilience}
Suppose $G$ is a graph as in Theorem~\ref{thm:sparse_mindegree_cycle} and has even order. To prove our resilience result, we use two main ingredients.  
The first ingredient is a structural result about $G$, which states that $G$ can be partitioned into balanced bipartite expanders that can be arranged cyclically so that consecutive expanders are joined by an edge (see Theorem~\ref{thm:cyclic}). This step relies on the \textit{sparse regularity lemma}, which we use to find a cycle of sparse regular pairs covering almost every vertex of $G$. The second step is to clean each cluster in the cycle so that the pairs have large minimum degree, and the last and most difficult step is to reallocate the leftover vertices so that each vertex of $G$ belongs to some bipartite expander. Although the task is reminiscent of applications of the (dense) regularity lemma, the sparse setting is considerably more delicate: for instance, moving even a constant number of vertices between the sets of the partition can, without proper control, quickly leave some vertices isolated. Moreover, we need to carefully control how many vertices are moved to each part, as we need to ensure that each bipartite expander is balanced, that is, it has the same number of vertices in each part. The proof of our partition result is given in~\Cref{sec:decomposition}.

The second ingredient is \Cref{thm:bipartite-c-expander}, which extends~\Cref{thm:C-expander} by asserting that, in suitable balanced bipartite expanders, any two vertices in different parts are the endpoints of a Hamilton path. Its proof is given in Appendix~\ref{app:bipartite-connectivity} and follows the argument of~\cite{C-expander}. 
%If one is interested only in the bipartite-expander setting, these extra conditions can be discarded and the same argument (indeed, a simpler one) yields \Cref{thm:bipartite-c-expander} \emph{mutatis mutandis}.

The claim of \Cref{thm:sparse_mindegree_cycle} follows directly from these two ingredients. Indeed, we use Theorem~\ref{thm:cyclic} to partition $G$ into balanced bipartite expanders arranged cyclically. Then we use Theorem~\ref{thm:bipartite-c-expander} to find a Hamilton path in each of them, and we choose the endpoints so that these paths and the edges between consecutive expanders form a Hamilton cycle in $G$. If $G$ has odd order, one vertex is left outside the bipartite expanders and is inserted into the cycle using property~\ref{cyclic:4} of Theorem~\ref{thm:cyclic}. The formal deduction is given in~\Cref{sec:resilience}.

\paragraph{Packing and covering}

To streamline the outline, assume that $G$ is as in Theorem~\ref{thm: decomp} and our goal is to find $(1-\varepsilon)d/2$ edge-disjoint Hamilton cycles. Our overall strategy will be to iteratively remove Hamilton cycles, and to extract Hamilton cycles our approach will be similar to that of Theorem~\ref{thm:C-expander}, which can be summarized as follows:
\begin{enumerate}[label=\upshape{\textbf{S\arabic{enumi}}}]
    \item\label{step:1} Find an \textit{$(A,B)$-linking structure} with the property that, for any bijection $\varphi:A\to B$, one can find vertex-disjoint paths spanning the linking structure, with endpoints $a$ and $\varphi(a)$ for each $a\in A$.  
    \item\label{step:2} Find a spanning linear forest $\mathcal F$ on the vertices outside the linking structure together with $A\cup B$, whose endpoints are exactly $A\cup B$.
    \item\label{step:3}Connect the paths in $\mathcal F$ using the property of the linking structure so that it completes to a Hamilton cycle. 
\end{enumerate}
As in each round we remove a 2-regular graph, the remaining graph maintains minimum degree at least $(\varepsilon-\gamma)d$, which, together with the sparsity condition, guarantees local expansion throughout the process. During the iterative removal of cycles, however, the remainder might develop large bipartite holes, and, once this occurs, techniques from~\cite{C-expander} are no longer applicable to perform steps~\ref{step:1} and~\ref{step:2}. Our goal will then be to find Hamilton cycles that do not overuse any part of $G$ so that, at each step of the removal process, we preserve the property that any two sufficiently large vertex sets span at least one edge.

To make this approach work, we first perform a random edge-splitting: each edge of $G$ is placed into $G_1$ with probability $1-q$, and into $G_2$ otherwise. This ensures that $G_2$ contains many edges between any two sufficiently large vertex sets.
The plan is to extract Hamilton cycles through steps~\ref{step:1}--\ref{step:3} using (almost) only edges of $G_1$, resorting to $G_2$ sparingly. Specifically, we will use $G_2$-edges only when we need an edge across two large sets, so, if successful, $G_2$ will retain the `no large bipartite holes' property until the very end.

Our starting point for step~\ref{step:2} is to find a spanning linear forest with as few paths as possible, for which we use a bound on the fractional linear arboricity~\cite{christoph2025linear} (see Lemma~\ref{lem:pathpartition}). The graph $G_1$ is only nearly regular, so the bound on the number of paths given by Lemma~\ref{lem:pathpartition} depends on its degree variation. To reduce the number of paths, we efficiently connect the endpoints of those paths using \textit{rotations} until only a few paths remain, while also controlling their lengths, relying on $G_2$ each time we connect two paths. Therefore, the number of $G_2$-edges we are forced to use depends on how irregular $G_1$ becomes after previous deletions, which creates a positive feedback loop: the more we tap into $G_2$, the more irregular $G_1$ becomes, which in turn compels further use of $G_2$.

The central tool to control this loop is \Cref{lem:regularHamcycle}, which produces a Hamilton path with prescribed endpoints while tightly controlling how many (and where) $G_2$-edges are used. This key lemma is proved in \Cref{sec:Hamilton-conn-bipartite}, which also uses the ideas from~\cite{C-expander} to rotate the linear forest we found first so that all the endpoints are suitably placed on the linking structure we found in step~\ref{step:1} (see Lemmas~\ref{lem: partition}--\ref{lemma:linearforest decomp2}). 

In \Cref{sec:packing-covering} we iteratively apply Lemma~\ref{lem:regularHamcycle} to obtain our packing result (\Cref{thm: decomp}), and in the same section we also prove our covering result (\Cref{thm: covering}). After applying the packing theorem with $\varepsilon^5$ in place of $\varepsilon$, we are left with a remainder whose maximum degree is at most $2\varepsilon^5 d$. We then cover these leftover edges by at most $2\varepsilon^3d$ additional Hamilton cycles, again using edges of the original graph $G$.

More precisely, we first decompose the leftover into matchings; then we split each matching randomly into a large constant number of parts; finally, as every vertex has many neighbours outside the vertices of each part, each part can be completed to a Hamilton cycle via a careful application of the extendability method together with the tools from \cite{C-expander}. This process is captured in a slightly more general result (\Cref{lem: partition}). While in standard uses of the extendability method the connecting paths are of length logarithmic in $n$, here we require paths of constant average length, as the matchings are of linear size. A similar idea to overcome this issue was also employed in~\cite{draganic2025disjoint}.

\section{Tools and preliminaries}\label{sec:prelim}

For a positive integer $n\in \mathbb{N}$, we write $[n]$ to denote the set $\{1,\dots, n\}$. We use standard hierarchy notation throughout the paper, that is, we use `$a\ll b$' to mean that there is a non-decreasing function $f$ such that the following statement holds whenever $0<a\le f(b)$.  Hierarchies with more variables are defined analogously and are meant to be read from right to left. When $1/m$ appears in some hierarchy, we implicitly assume that $m$ is a positive integer. Also, we write $a=b\pm c$ if $b-c \leq a \leq b+c$ holds, and we often omit floors and ceilings and treat large numbers as integers if it does not affect the argument. All logarithms are natural. 

 We use standard graph theory notation. Given a graph $G$, let $V(G)$ and $E(G)$ denote the vertex set and edge set of $G$, respectively, and we write $|G|=|V(G)|$ and $e(G)=|E(G)|$ for the numbers of vertices and edges of $G$.  For vertex sets $U,W\subset V(G)$, we write $G[U]$ for the subgraph induced by $U$, and let $e_G(U,W)$ be the number of pairs $(u,w)\in U\times W$ with $uw\in E(G)$. In particular, when $U$ and $W$ are disjoint,  $e_G(U,W)$ counts the number of edges of the bipartite graph induced by $U$ and $W$, which we denote by $G[U,W]$. Also, note that $e_G(U,U)= 2e(G[U])$ counts all edges in the induced subgraph $G[U]$ twice.  The set $N_G(x)$ denotes the neighbours of $x\in V(G)$, and $d_G(x)=|N_G(x)|$ is the degree of $x$. Also, given a set $S\subset V(G)$, we let $d_G(x,S)$ denote the number of neighbours of $x$ in $S$. We omit the subscript $G$ in these notations whenever the graph $G$ is clear from the context. 

 For a set $E$ of edges, we write $V(E)$ to denote the set of vertices that the edges in $E$ are incident with, and we write $G+E$ to denote the graph whose vertex set is $V(G)\cup V(E)$ and edge set is $E(G)\cup E$, and $G\setminus E$ means the graph with vertex set $V(G)$ and edges $E(G)\setminus E$. Given a set $S\subset V(G)$, $G-S$ is the graph induced by $V(G)\setminus S$ in $G$. Given two graphs $G_1, G_2$, we use $G_1\cup G_2$ to denote the graph with vertex set $V(G_1)\cup V(G_2)$ and edge set $E(G_1)\cup E(G_2)$. Given distinct vertices $u,v\in V(G)$, say a path $P$ in $G$ is a $u,v$-path if $P$ has endpoints $u$ and $v$, and say the length of $P$ is its number of edges. A linear forest in $G$ is a subgraph of $G$ whose components are paths, where an isolated vertex is a path of length $0$, and it is spanning if its vertex set is $V(G)$. For a linear forest $\mathcal F$, let $End(\mathcal F)$ denote the set of endpoints of its paths. For a given graph $G$ with a bipartition $V_1\cup V_2$, we say that a set $U$ is balanced  if $|U\cap V_1|=|U\cap V_2|$.
\subsection{Probabilistic estimates}
We will need the following well-known probabilistic results.
\begin{lemma}[Chernoff's bound~\cite{alonspencer2016}]\label{lem:Chernoff}Let $X$ be a binomial random variable. Then, for any $0<\varepsilon\le 3/2$, 
\[\mathbb P\big(|X-\mathbb EX|\ge \varepsilon\mathbb EX\big)\le 2e^{-\varepsilon^2\mathbb EX/3}.\]
\end{lemma}

\begin{lemma}[Lov\'asz's local lemma~\cite{alonspencer2016}]\label{Lem:LLL}
Let $A_1,\dots, A_n$ be events such that each event $A_i$ is mutually independent of all the other events but at most $d$ of them. If $\mathbb{P}(A_i)\leq p$ for all $i\in [n]$ and $ep(d+1)<1$, then we have 
\[\textstyle\mathbb{P}\Big(\bigcap_{i\in [n]} \overline{A_i}\Big) \geq \left(1-\tfrac{1}{d+1}\right)^{n}.\]
\end{lemma}
\begin{lemma}[McDiarmid's inequality~\cite{mcdiarmid1989}]\label{lem:mcdiarmid}Let $X_1,\ldots, X_n$ be independent random variables taking values in a set $\Omega$. Let $c_1,\ldots, c_n\ge 0$ and let $f:\Omega^n\to\mathbb R$ be a function so that for every $i\in [n]$ and $x_1,\ldots, x_n,x_i'\in\Omega$, we have $\Big |f(x_1,\ldots, x_i,\ldots, x_n)-f(x_1,\ldots, x_i',\ldots, x_n)\Big|\le c_i$. Then, for all $t>0$, 
\[\mathbb P\left (\big|f(X_1,\ldots, X_n)-\mathbb E(f(X_1,\ldots, X_n))\big|\ge t\right)\le 2\exp\left(-\tfrac{t^2}{2\sum_{i\in [n]}c_i^2}\right).\]
    
\end{lemma}
In the next lemma, we split a set into $k$ parts so that each vertex has about the same number of neighbours in each part.
\begin{lemma}\label{lem:balanced-random-splitting}
Let $1/d\ll1/k$. Let $J$ be a graph with $\Delta(J)\leq2d$, and let $W\subset V(J)$ be partitioned into sets of size $k$. There is a partition $V_1,\ldots,V_k$ of $W$ such that each $V_i$ meets each of these sets in one vertex and
\[
 \bigl|d_J(v,V_i)-\tfrac1kd_J(v,W)\bigr|\leq d^{2/3}
\]
for every vertex $v\in V(J)$ and $i\in[k]$.
\end{lemma}
\begin{proof}
For each of the sets, choose independently and uniformly a bijection to $[k]$, and put each vertex into $V_i$ for its image $i$. Fix $v$ and $i$. Then $\mathbb E[d_J(v,V_i)]=d_J(v,W)/k$. Each set contributes at most one to this count, and thus at most $d_J(v)\leq2d$ choices affect it, each by at most one. By McDiarmid's inequality (Lemma~\ref{lem:mcdiarmid}), the bound fails with probability at most $2\exp(-d^{4/3}/(4d))=2\exp(-d^{1/3}/4)$. The failure event for $(v,i)$ is determined by the bijections chosen for the at most $d_J(v)\leq2d$ sets that meet $N_J(v)$, and hence it is mutually independent of all failure events for pairs $(u,j)$ such that no set meets both $N_J(v)$ and $N_J(u)$ (possibly in different vertices). Each set meets the neighbourhoods of at most $2kd$ vertices, and each vertex gives $k$ events. So each event is independent of all but at most $4k^2d^2$ others, and Lov\'asz's local lemma (Lemma~\ref{Lem:LLL}) applies as $e(4k^2d^2+1)\cdot2\exp(-d^{1/3}/4)<1$ by $1/d\ll1/k$.
\end{proof}
\subsection{Partitioning an almost-regular graph into few paths}
In Lemma~\ref{lem:pathpartition}, we partition an almost-regular graph into few paths. Christoph, Dragani\'c, Gir\~ao, Hurley, Michel and M\"uyesser~\cite{christoph2025linear} showed that the vertices of every $d$-regular $n$-vertex graph can be partitioned into at most $2n/(d+1)$ paths. Earlier, Montgomery, M\"uyesser, Pokrovskiy and Sudakov~\cite{montgomery2024approximate} showed that, for every $\eps>0$ and every $d$ large in terms of $\eps$, all but at most $\eps n$ vertices of such a graph can be partitioned into at most $n/(d+1)$ paths. We use the following bound of~\cite{christoph2025linear} on the fractional linear arboricity, stated in terms of random linear forests.

\begin{theorem}[{\cite[Theorem~1.3]{christoph2025linear}}]\label{thm:pathpartition}
Every graph $G$ with maximum degree $\Delta\geq2$ has a random linear forest $F\subset G$ that contains each edge of $G$ with probability at least $1/(\Delta/2+O(\log\Delta))$.
\end{theorem}
\begin{lemma}\label{lem:pathpartition}Let $1/d\ll\varepsilon$ and $n,d\in\mathbb N$. If $G$ is an $n$-vertex graph such that $d\le d(x)\le (1+\eps)d$ for all $x\in V(G)$, then $V(G)$ can be partitioned into at most $3\varepsilon n$ paths.
\end{lemma}
\begin{proof}
Let $\Delta=\Delta(G)$, so that $d\leq\Delta\leq(1+\eps)d$. By Theorem~\ref{thm:pathpartition}, there is a random linear forest $F\subset G$ that contains each edge of $G$ with probability at least $1/(\Delta/2+O(\log\Delta))\geq2/((1+2\eps)d)$, as $1/d\ll\eps$. As every vertex has degree at least $d$, we have $e(G)\geq nd/2$, and thus $\mathbb E[e(F)]\geq2e(G)/((1+2\eps)d)\geq n/(1+2\eps)$. We fix an outcome of $F$ with $e(F)\geq n/(1+2\eps)$. Its paths, together with the vertices outside $V(F)$ as paths of length $0$, partition $V(G)$ into $n-e(F)\leq2\eps n/(1+2\eps)\leq3\eps n$ paths.
\end{proof}

\subsection{The sparse regularity lemma}\label{sec:sparseregularity}
In this section, we collect some definitions and results regarding regular partitions of $(\eta,\beta,p)$-sparse graphs. For $p>0$, a graph $G$ and nonempty disjoint sets $A, B \subset V(G)$, the \emph{$p$-density} of the pair $(A, B)$ is defined as
    \[d_p(A,B)=\frac{e(A,B)}{p|A||B|}.\]
We say the pair $(A, B)$ is \emph{$(\eps, p)$-regular} if for all subsets $A'\subset A$ and $B'\subset B$, with $|A'|\ge \eps|A|$ and $|B'|\ge \eps|B|$, 
    \[e(A',B')=(d_p(A,B)\pm\eps)p|A'||B'|.\]
Kohayakawa~\cite{kohayakawa1997szemeredi} and R\"odl (unpublished) proved a sparse analogue of the \textit{regularity lemma} of Szemer\'edi~\cite{szemeredi1978regular}. We use the following minimum degree version of it.
\begin{lemma}[{\cite[Lemma~13]{ABET}}]\label{thm:regularity}
For all $\eps>0$, $c\in[0,1]$ and $k_0\in\mathbb N$ there are $K_0\ge k_0$ and $n_0$ such that the following holds for all $\alpha\in[0,1]$, $p>0$ and $n\geq n_0$. Let $G$ be an $n$-vertex graph with $\delta(G)\ge cpn$ such that $e(X,Y)\le(1+\eps^2/1000)p|X||Y|$ for all disjoint $X,Y\subset V(G)$ with $|X|,|Y|\ge\eps n/K_0$. Then there are a partition $V(G)=V_0\cup V_1\cup\dots\cup V_k$ with $k_0\le k\le K_0$, $|V_0|\le\eps n$ and $|V_1|=\dots=|V_k|$, and a graph $R$ on $[k]$ with $\delta(R)\ge(c-\alpha-\eps)k$, such that for each $ij\in E(R)$ there is $d'\ge\alpha$ with $e(A,B)=(d'\pm\eps)p|A||B|$ for all $A\subset V_i$, $B\subset V_j$ with $|A|\ge\eps|V_i|$ and $|B|\ge\eps|V_j|$.
\end{lemma}
As usual, such a partition is called a \emph{regular partition} of $G$, each $V_i$ with $i\in[k]$ is called a \emph{cluster}, and $R$ is called the \emph{reduced graph}.

\subsection{Expansion and \texorpdfstring{$(\eta,\beta,p)$}{(eta,beta,p)}-sparsity}
Recall that given a graph $G$ and a subset $U\subset V(G)$, the \textit{external neighbourhood} of $U$ in $G$, denoted $N_G(U)$, is the set of vertices in $ V(G)\setminus U$ that have neighbours in $U$. We write $\Gamma_G(U)=\bigcup_{u\in U}N_G(u)$ for the set of all neighbours of vertices in $U$. We will use different notions of expansion, which we define next.
\begin{definition}
    Let $G$ be a graph and let $U,V\subset V(G)$. 
    \begin{itemize}
        \item We say  $U$ is $(m,C)$-expanding  into $V$ in $G$ if for every $X\subset U$ with $|X|\leq m$, we have $|N(X)\cap V| \geq C|X|$. Moreover, we say $G$ is \emph{$(m,C)$-expanding} if $V(G)$ is $(m,C)$-expanding into $V(G)$.
        \item We say $G$ is \emph{$m$-joined} if there is an edge between $A$ and $B$ for any disjoint $A,B\subset V(G)$ with $|A|,|B|\geq m$.
        \item If $G$ is bipartite with parts $V_1$ and $V_2$, we say $G$ is \emph{$m$-bipartite-joined} if there is an edge between $A$ and $B$ for any $A\subset V_1$ and $B\subset V_2$ with $|A|, |B|\geq m$.
    \end{itemize}  
\end{definition}
We remark that an $n$-vertex graph $G$ with $n\geq3$ is a $C$-expander (as defined in the introduction) if $G$ is $(n/(2C),C)$-expanding and $n/(2C)$-joined.
\begin{lemma}\label{lem:few-heavy}
Let $G$ be an $N$-vertex $(\eta,\beta,p)$-sparse graph with $\beta\leq1$, let $T\subset V(G)$, and let $h>(1+\beta)p\max\{\eta N,|T|\}$. Then at most $2\eta pN|T|/h$ vertices of $G$ have at least $h$ neighbours in $T$.
\end{lemma}
\begin{proof}
Let $B$ be the set of these vertices. We may assume that $B\neq\emptyset$, so that $h\leq|T|$, and that $|B|>\eta pN$, as otherwise $|B|\leq\eta pN\leq2\eta pN|T|/h$. Then $|T|\geq h>\eta pN$ and $e_G(B,T)\geq h|B|$. If $|T|\leq|B|$, then $e_G(B,T)\leq(1+\beta)p\max\{\eta N,|T|\}|B|<h|B|$ by sparsity, a contradiction. So $|B|<|T|$, and sparsity gives $h|B|\leq(1+\beta)p\max\{\eta N,|B|\}|T|$. Hence $|B|<\eta N$ as $h>(1+\beta)p|T|$, and then $|B|\leq(1+\beta)\eta pN|T|/h\leq2\eta pN|T|/h$.
\end{proof}
The next lemma states that having some minimum degree condition in an $(\eta,\beta,p)$-sparse graph is enough to guarantee expansion.
\begin{lemma}\label{lem:sparse->expander}Let $D\ge 1$ and let $0<\alpha,\beta,\eta< 1$ satisfy $\eta < \alpha/(4D)$. Let $G$ be an $n$-vertex $(\eta,\beta,d/n)$-sparse graph with $d>0$ and let $X,Y\subset V(G)$ be two not necessarily disjoint sets. 
Suppose that every vertex $x\in X$ satisfies $|N_{G}(x)\cap Y|\geq \alpha d$. Then $X$ is $(\alpha n/(2D+2),D)$-expanding into $Y$ in $G$, and every $U\subset X$ with $|U|\geq\eta n$ satisfies $|\Gamma_G(U)\cap Y|\geq\alpha n/2$.
\end{lemma}
\begin{proof}[Proof of Lemma~\ref{lem:sparse->expander}]
Suppose that there is a nonempty set $U\subset X$ with $|U|\leq\alpha n/(2D+2)$ and $|N_G(U)\cap Y|<D|U|$, and let $W=U\cup(N_G(U)\cap Y)$. Then $|W|\leq(D+1)|U|\leq\alpha n/2$, and each vertex of $U$ has at least $\alpha d$ neighbours in $W$. Since $\beta<1$, $\eta<\alpha/4$ and $|W|\leq\alpha n/2$, we have $\alpha d>(1+\beta)(d/n)\max\{\eta n,|W|\}$, and Lemma~\ref{lem:few-heavy} with $T=W$ and $h=\alpha d$ gives $|U|\leq2\eta|W|/\alpha\leq2\eta(D+1)|U|/\alpha<|U|$, a contradiction.

For the second assertion, let $U\subset X$ with $|U|\geq\eta n$ and let $W=\Gamma_G(U)\cap Y$. If $|W|<\alpha n/2$, then $|U|\leq2\eta|W|/\alpha<\eta n$ by Lemma~\ref{lem:few-heavy} as before, a contradiction.
\end{proof}
We now show that the Expander mixing lemma (Theorem 2.11 in~\cite{krivelevich2006pseudo}) implies that $(n,d,\lambda)$-graphs are $(\eta, \beta, p)$-sparse for a suitable choice of $\eta,\beta$ and $p$.

\begin{lemma}\label{lemma:nd:sparse}If $G$ is an $(n,d,\lambda)$-graph with $\lambda /d\leq \beta \eta$, then $G$ is $(\eta,\beta,d/n)$-sparse.
    \end{lemma}
    
\begin{proof}[Proof of Lemma~\ref{lemma:nd:sparse}]
Take $U,W\subset V(G)$ with $\eta d\le |U|\le |W|$. The Expander mixing lemma implies $|e_G(U,W)-d|U||W|/n|\le\lambda\sqrt{|U||W|}$, so
\[
e_G(U,W)\le \frac{d|U||W|}{n}+\lambda\sqrt{|U||W|}\le \frac{d\max\{\eta n,|U|\}|W|}{n}+\lambda |W|\le (1+\beta)\frac{d}n\max\{\eta n,|U|\}|W|,
\]
where the last inequality follows from the fact that $\lambda /d\leq \beta \eta$.
\end{proof}

To deduce the results on $(n,d,\lambda)$-graphs in the introduction from those on sparse graphs, we use Lemma~\ref{lemma:nd:sparse}. For Theorem~\ref{thm:main}, note that a spanning subgraph of an $(\eta,\beta,d/n)$-sparse graph is $(\eta,\beta,d/n)$-sparse. The theorems on sparse graphs also require $d\geq d_0$, which follows from $d/\lambda\geq C\geq d_0$ since $\lambda\geq1$. Indeed, for $1\leq d<n$, the trace of the square of the adjacency matrix gives
\[
 (n-1)\lambda^2\geq nd-d^2=d(n-d)\geq n-1.
\]
Thus $\lambda\geq1$.

\begin{lemma}\label{lem:restriction-normalization}
Let $G$ be an $N$-vertex graph and let $W\subset V(G)$ with $|W|=M>0$. If $G$ is $(\eta,\beta,p)$-sparse, then $G[W]$ is $(\eta,\beta,pN/M)$-sparse.
\end{lemma}
\begin{proof}
Let $U,T\subset W$ with $\eta(pN/M)M=\eta pN\leq|U|\leq|T|$. Since $G$ is $(\eta,\beta,p)$-sparse, we have $e(U,T)\leq(1+\beta)p\max\{\eta N,|U|\}|T|$, and $p\max\{\eta N,|U|\}\leq(pN/M)\max\{\eta M,|U|\}$ as $M\leq N$.
\end{proof}

\section{Partition into bipartite expanders}\label{sec:decomposition}

In this section, we consider $(\eta,\beta, d/n)$-sparse graphs with minimum degree at least $(1/2+\gamma)d$, as in the setup of Theorem~\ref{thm:sparse_mindegree_cycle}, and prove a structural result for them. For such a graph $G$, we will find vertex-disjoint bipartite expanders $G_1,\ldots,G_K$ covering all but at most one vertex of $G$, with an edge between consecutive expanders in cyclic order. To state our result, we define bipartite expanders first; the notion we need is analogous to that of $C$-expanders as in \Cref{thm:C-expander}.

\begin{definition}
We say that a bipartite graph $H$, with parts $U_1$ and $U_2$ with $|U_1|=|U_2|=n'$, is a \textit{bipartite $C$-expander} if $H$ is $\frac{n'}{2C}$-bipartite-joined and $U_i$ is $(\frac{n'}{2C},C)$-expanding into $U_{3-i}$ for each $i\in [2]$.
\end{definition}

The main result of this section is the following.

\begin{theorem}\label{thm:cyclic}
Let $1/d_0\ll\eta\ll\beta\ll1/D,\gamma$ with $\gamma<1/2$, and let $n,d\in \mathbb{N}$ satisfy $n>d\ge d_0$. Suppose $G$ is an $(\eta,\beta,d/n)$-sparse graph with $n$ vertices and $\delta(G)\geq (1/2+ \gamma)d$. Then, there are $K\in\mathbb N$ and a partition $V(G)=U_0\cup U_1\cup \ldots \cup U_{2K}$ such that the following holds for each $i\in [K]$ (where $U_{2K+1}=U_1$). 
\begin{enumerate}[label= \upshape\textbf{\Alph{propcounter}\arabic{enumi}}]
        \item  $|U_{2i-1}|=|U_{2i}|$. \label{cyclic:1}
        \item $G[U_{2i-1},U_{2i}]$ is a bipartite $D$-expander. \label{cyclic:2}
        \item $e(U_{2i},U_{2i+1})>0$.\label{cyclic:3}
        \item  If $n$ is even, then $U_0$ is empty; otherwise, $U_0$ consists of a single vertex $u_0$ which has a neighbour in both $U_1$ and $U_{2K}$.\label{cyclic:4}
         \end{enumerate} \stepcounter{propcounter}
\end{theorem}
The proof of Theorem~\ref{thm:cyclic} relies on the \textit{regularity method} for sparse graphs, which we sketch now. Suppose that $G$ has a regular partition $V(G)=V_0\cup V_1\cup \ldots \cup V_{2K}$ with reduced graph $R$, as given by Lemma~\ref{thm:regularity}. (We may always assume that the number of clusters is even at the cost of moving one cluster to $V_0$.) Although $G$ is a sparse graph (we may think of $d$ as being a large constant), the reduced graph $R$ is dense, and its minimum degree is at least $|R|/2$, thus Dirac's theorem implies $R$ is Hamiltonian. Therefore, there is a Hamilton cycle in $R$ whose consecutive clusters correspond to regular pairs with large $p$-density. To prove Theorem~\ref{thm:cyclic}, we need to overcome two issues: a) regular pairs might not be bipartite expanders, and b) vertices in $V_0$ need to be incorporated into this cyclic structure somehow.
We address a) by showing that, after removing a few vertices, regular pairs contain almost-spanning subgraphs of large minimum degree. We move all the leftover vertices to $V_0$ and then show that each vertex in $V_0$ can be reallocated to some cluster in the cycle so that the pairs become balanced bipartite expanders, thus resolving b).  
While this plan of vertex relocation may seem similar to the approach used in the regularity method for dense graphs, several new problems arise in the sparse setting.  For example, moving even a constant number of vertices from a set $V_{2i}$ to $V_0$ may create isolated vertices in $G[V_{2i-1}, V_{2i}]$. A much more delicate analysis is needed to address such issues. Property~\ref{cyclic:4} comes into play when $n$ is odd, in which case it is not possible to partition $V(G)$ into an even number of parts satisfying~\ref{cyclic:1}. 
%-------------------------------------------
\begin{proof}[Proof of Theorem~\ref{thm:cyclic}] We introduce additional constants $\eps,\alpha$ with $\eps\ll\alpha\ll1/D,\gamma$.
Let $p=d/n$, and let $K_0$ and $n_0$ be given by Lemma~\ref{thm:regularity} for $\eps^{10}$, $1/2+\gamma$ and $\eps^{-10}$ in the roles of $\eps$, $c$ and $k_0$. As $K_0$ and $n_0$ depend only on $\eps$ and $\gamma$, we may assume that $1/d_0\ll\eta\ll\beta\ll1/K_0$ and $d_0\geq n_0$, so that $n>d\geq n_0$. Since $p<1$ and $\eta\ll\beta\ll1/K_0$, $G$ satisfies the condition of Lemma~\ref{thm:regularity} on $e(X,Y)$ by Definition~\ref{def:sparse}, and $\delta(G)\geq(1/2+\gamma)pn$. We may therefore apply Lemma~\ref{thm:regularity} with $2\alpha$ in the role of $\alpha$, and we obtain a partition $V(G)=V_0\cup V_1\cup\dots\cup V_k$ with $\eps^{-10}\leq k\leq K_0$, $|V_0|\leq\eps^{10}n$ and $|V_1|=\dots=|V_k|$, and a graph $R$ on $[k]$ with $\delta(R)\geq(1/2+\gamma-2\alpha-\eps^{10})k$. If $k$ is odd, we move $V_k$ to $V_0$ and delete $k$ from $R$, and we write $2K\geq\eps^{-10}-1$ for the number of remaining clusters. Then $|V_0|\leq\eps^{10}n+n/k\leq2\eps^{10}n$ and $\delta(R)\geq(1+\gamma)K$ as $\eps\ll\alpha\ll\gamma$, and for each $ij\in E(R)$ the pair $(V_i,V_j)$ is $(2\eps^{10},p)$-regular with $p$-density at least $2\alpha-\eps^{10}\geq\alpha$. Dirac's theorem yields a Hamilton cycle in $R$ and, relabelling if necessary, we may assume that $V_1 V_2 \ldots V_{2K}$ is that cycle. Let $m$ be the common size of $V_1,\ldots,V_{2K}$. Let $M$ be the perfect matching consisting of the edges $V_{2i-1}V_{2i}$, $i\in [K]$, and let $M(V_j)$ denote the neighbour of $V_j$ in $M$.

\paragraph{Cleaning the pairs.}
For $j\in[2K]$, write $V_{j'}=M(V_j)$, and let $V_j^\circ$ denote the cluster $V_j$ before the cleaning. Since $(V_j^\circ,V_{j'}^\circ)$ is $(2\varepsilon^{10},p)$-regular with $p$-density at least $\alpha$, less than $2\varepsilon^{10}m$ vertices of $V_j^\circ$ have less than $(\alpha-2\varepsilon^{10})pm$ neighbours in $V_{j'}^\circ$. As long as some vertex of some $V_j$ has less than $\alpha d/(10K)$ neighbours in the current set $V_{j'}$, we move such a vertex to $V_0$. Let $Z_j$ be the set of vertices moved from $V_j$ so far. We claim that $|Z_j|<4\varepsilon^{10}m$ for all $j\in[2K]$ throughout. Otherwise, consider the first move after which $|Z_j|\geq4\varepsilon^{10}m$ for some $j$. Then $|Z_{j'}|<4\varepsilon^{10}m$. Moreover, more than $2\varepsilon^{10}m$ vertices of $Z_j$ have at least $(\alpha-2\varepsilon^{10})pm$ neighbours in $V_{j'}^\circ$ but less than $\alpha d/(10K)$ in $V_{j'}$, since no vertex was added to $V_{j'}$ after they were moved. Since $m\geq(1-\varepsilon^{10})n/(2K+1)\geq n/(3K)$, we have $\alpha d/(10K)\leq3\alpha pm/10$, and hence each of these vertices has at least $h:=\alpha pm/2$ neighbours in $V_{j'}^\circ\setminus V_{j'}=Z_{j'}$. As $\eta\ll1/K_0$ and $\eps\ll\alpha$, we have $h>(1+\beta)p\max\{\eta n,|Z_{j'}|\}$, and Lemma~\ref{lem:few-heavy} with $T=Z_{j'}$ shows that at most $4\eta n|Z_{j'}|/(\alpha m)\leq|Z_{j'}|/2<2\varepsilon^{10}m$ vertices have at least $h$ neighbours in $Z_{j'}$, a contradiction. Hence, when no vertex can be moved, we have $|V_0|\leq2\varepsilon^{10}n+2K\cdot4\varepsilon^{10}m\leq6\varepsilon^{10}n$, for all $j\in[2K]$ we have $(1-4\varepsilon^{10})m\leq|V_j|\leq m\leq n/(2K)$, and $\delta(G[V_{2i-1},V_{2i}])\geq\alpha d/(10K)$ for all $i\in[K]$.

Our goal in the rest of the proof is to relocate the vertices in $V_0$ to $\bigcup_{i\in [2K]}V_i$ so that, for each $j\in [2K]$,  $G[V_j,M(V_j)]$ is a balanced bipartite graph with minimum degree at least $\alpha d/(100K)$. We update the sets $V_i$ as vertices are moved, keeping the same notation.

\begin{claim}\label{cl: degrees in half of the sets}
All but at most $2n/K^3$ vertices $v\in V(G)$ have at least $\alpha d/(10K)$ neighbours in at least $(1+\gamma/10)K$ sets $V_i$, $i\in [2K]$. We call those vertices good. 
\end{claim}
\begin{claimproof}[Proof of Claim~\ref{cl: degrees in half of the sets}]
Let $B$ be the set of vertices in $G$ that have at least $(1+\gamma/10)p|V_i|$ neighbours in some $V_i$ with $i\in [2K]$ or at least $\gamma d/100$ neighbours in $V_0$. Lemma~\ref{lem:few-heavy}, with $T=V_i$ and $h=(1+\gamma/10)p|V_i|$ for each $i\in[2K]$, and with $T=V_0$ and $h=\gamma d/100$, applies since $\eta n\le|V_i|$, $|V_0|\le6\varepsilon^{10}n$ and $\eta,\beta,\varepsilon\ll\gamma$, which implies $|B|\leq(2K+1)2\eta n\le 2n/K^3$ as $\eta\ll1/K_0$. It suffices to show that every $v\in V(G)\setminus B$ is good. Suppose that some $v\in V(G)\setminus B$ is not good. Then
    \begin{eqnarray*}
    d(v)\leq \gamma d/100 + (1+\gamma/10)K \cdot  (1+\gamma/10)d/(2K) + 2K \cdot \alpha d/(10K)
    &\leq& \gamma d/100+ (1/2+\gamma/9)d+\alpha d/5\\&<&(1/2+\gamma)d,
    \end{eqnarray*}
contradicting the bound on the minimum degree of $G$.
\end{claimproof}

\paragraph{Integration of $V_0$.}
Every nonempty $T\subset V_0$ contains a vertex with less than $d/4$ neighbours in $T$. Indeed, as $|T|\leq|V_0|\leq6\varepsilon^{10}n$, Lemma~\ref{lem:few-heavy} with $h=d/4$ shows that at most $8\eta|T|<|T|$ vertices have at least $d/4$ neighbours in $T$. We may thus order $V_0$ as $(v_1,\ldots,v_\ell)$ so that each $v_i$ has less than $d/4$ neighbours in $\{v_i,\ldots,v_\ell\}$, and hence more than $\delta(G)-d/4>d/4$ neighbours outside it. We move $v_1,\ldots,v_\ell$ in this order, so that when $v_i$ is moved, its neighbours outside $\{v_i,\ldots,v_\ell\}$ already lie in $\bigcup_{j\in[2K]}V_j$. If $v_i$ is good, we take $j\in[2K]$ such that $v_i$ has at least $\alpha d/(10K)$ neighbours in $V_j$ and $j$ was taken less than $\varepsilon^8n/K$ times before. Such a choice is possible, as no vertex was removed from the sets $V_j$ after Claim~\ref{cl: degrees in half of the sets}, and hence $v_i$ has at least $\alpha d/(10K)$ neighbours in at least $(1+\gamma/10)K$ of them, and at most $|V_0|K/(\varepsilon^8n)<K$ indices were taken $\varepsilon^8n/K$ times. Otherwise, we take $j$ such that $v_i$ has at least $d/(8K)\geq\alpha d/(10K)$ neighbours in $V_j$. In both cases we add $v_i$ to $M(V_j)$. This keeps the minimum degree of each pair at least $\alpha d/(10K)$. Also, each $V_j$ receives at most $\varepsilon^8n/K+2n/K^3\leq2\varepsilon^8n/K$ vertices, as at most $2n/K^3$ vertices are not good.

  \paragraph{Correcting the imbalance.}  
All vertices of $V_0$ now lie in $\bigcup_{j\in[2K]}V_j$, but the pairs may be unbalanced. For $j\in[2K]$, let $t_j:=\max\{|V_j|-|M(V_j)|,0\}$. After the cleaning we had $|V_j|-|M(V_j)|\leq4\varepsilon^{10}m$, and since then $V_j$ has received at most $2\varepsilon^8n/K$ vertices, so $0\leq t_j\leq4\varepsilon^{10}m+2\varepsilon^8n/K\leq 3\varepsilon^8n/K$. We remove from each $V_j$ a set $T_j$ of $t_j$ vertices and redistribute these vertices later.

\begin{claim}\label{cl:randomsubsets}
    There is a collection of sets  $(T_j)_{j\in [2K]}$ such that  for $U=\bigcup_{j\in [2K]} T_j$ the following holds.
    \begin{enumerate}[label=\upshape(\roman{enumi})]
        \item $T_j\subset V_j$ with $|T_j|=t_j$ for each $j\in [2K]$.
        \item For each $j\in[2K]$, every vertex with at least $\alpha d/(10K)$ neighbours in $V_j$ has at least $\alpha d/(100K)$ neighbours in $V_j\setminus T_j$.
        \item Every vertex in $U$ is good.
    \end{enumerate}
\end{claim}

\begin{claimproof}[Proof of Claim~\ref{cl:randomsubsets}]

Let $X$ be the set of vertices of $G$ with degree at least $10d$. Lemma~\ref{lem:few-heavy} with $T=V(G)$ and $h=10d$ implies $|X|\leq\eta n/5$.

For each $j\in [2K]$, let $L_j$ be the set of good vertices in $V_j\setminus X$. If $|L_j|$ is odd, we remove one vertex from $L_j$. Then $|L_j|\ge |V_j|-|X|-2n/K^3-1\ge n/(5K)$. Let $P$ be the set of $(v,j)\in V(G)\times[2K]$ such that $v$ has at least $\alpha d/(10K)$ neighbours in $V_j$ but less than $\alpha d/(100K)$ in $V_j\setminus L_j$. For each $(v,j)\in P$, we choose a set $S_{v,j}\subset N_G(v)\cap L_j$ of size $9\alpha d/(100K)$. Let $J$ be the bipartite graph between $P$ and $\bigcup_{j\in[2K]}L_j$ in which each $(v,j)\in P$ is joined to the vertices of $S_{v,j}$. Each $(v,j)\in P$ has degree $9\alpha d/(100K)\leq d$ in $J$, and each $w\in L_j$ has degree at most $d_G(w)<10d$ in $J$, as $w$ lies in $S_{v,j}$ only if $v\in N_G(w)$. We partition each $L_j$ into pairs and apply Lemma~\ref{lem:balanced-random-splitting} to $J$ and these pairs, with $5d$ playing the role of $d$. This yields, for each $j\in[2K]$, a set $R_j\subset L_j$ that meets each pair of $L_j$ in one vertex and satisfies $|S_{v,j}\cap R_j|\leq9\alpha d/(200K)+(5d)^{2/3}\leq8\alpha d/(100K)$ for all $(v,j)\in P$, as $1/d_0\ll\alpha,1/K_0$.

Then $|R_j|=|L_j|/2\geq n/(10K)\geq t_j$ as $t_j\le 3\varepsilon^8n/K$. Let $v\in V(G)$ have at least $\alpha d/(10K)$ neighbours in $V_j$. If $(v,j)\in P$, then $v$ has at least $|S_{v,j}\setminus R_j|\geq\alpha d/(100K)$ neighbours in $V_j\setminus R_j$. Otherwise, $v$ has at least $\alpha d/(100K)$ neighbours in $V_j\setminus L_j\subset V_j\setminus R_j$. Now, for each $j\in [2K]$, pick a set $T_j\subset R_j$ of size $|T_j|=t_j$, and let $U=\bigcup_{j\in [2K]}T_j$. This proves (i) and (ii), and (iii) holds as $R_j\subset L_j$.
\end{claimproof}

%Let $U=\bigcup T_i$ be the union of the vertices given by the previous claim. 
For each $j\in [2K]$, we remove $T_j$ from $V_j$, so that now $|V_j|=|M(V_j)|$. We still have $\delta(G[V_{2i-1},V_{2i}])\geq\alpha d/(100K)$ for all $i\in[K]$ by (ii), as $\delta(G[V_{2i-1},V_{2i}])\geq\alpha d/(10K)$ before the removal. Moreover, all vertices of $U$ are good, and no vertex was removed from the sets $V_j$ between Claims~\ref{cl: degrees in half of the sets} and~\ref{cl:randomsubsets}. Hence, by (ii), every vertex of $U$ has at least $\alpha d/(100K)$ neighbours in $V_j$ for at least $(1+\gamma/10)K$ indices $j\in[2K]$.

\paragraph{Reintegration of $U$ and completing the proof.}
Let $u,v\in U$. At least $\gamma K/10$ indices $j\in[K]$ are such that $u$ has at least $\alpha d/(100K)$ neighbours in $V_{2j-1}$ and $v$ in $V_{2j}$, or vice versa. Indeed, for $w\in\{u,v\}$ let $I(w)$ be the set of indices $i\in[2K]$ such that $w$ has at least $\alpha d/(100K)$ neighbours in $V_i$. Since $|I(u)|,|I(v)|\geq(1+\gamma/10)K$, at least $|I(u)|+|I(v)|-2K\geq\gamma K/5$ indices $i\in I(u)$ satisfy $i'\in I(v)$, where $V_{i'}=M(V_i)$. For each such $i$, the index $j=\lceil i/2\rceil$ is as claimed, and at most two such $i$ have the same $\lceil i/2\rceil$. We pair up the vertices of $U$, leaving at most one, and for each pair $\{u,v\}$ in turn we take such a $j$ that was taken less than $\varepsilon^6n/K$ times before. Such a choice is possible, as $|U|\le 2K\cdot 3\varepsilon^8n/K=6\varepsilon^8n$ and hence at most $(|U|/2)K/(\varepsilon^6n)\leq3\varepsilon^2K<\gamma K/10$ indices were taken $\varepsilon^6n/K$ times. We add $u$ to $V_{2j}$ and $v$ to $V_{2j-1}$, or vice versa, so that each has at least $\alpha d/(100K)$ neighbours on the other side. Then each $G[V_{2j-1},V_{2j}]$ stays balanced and has minimum degree at least $\alpha d/(100K)$. As each $G[V_{2j-1},V_{2j}]$ is balanced, $|U|$ and $n$ have the same parity. Hence one vertex $u_0\in U$ is left over if $n$ is odd, and no vertex is left over if $n$ is even.

If $n$ is odd, then $u_0$ has neighbours in at least $(1+\gamma/10)K>K$ of the sets $V_i$. Thus, by the pigeonhole principle, $u_0$ has a neighbour in both $V_{2j}$ and $V_{2j+1}$ for some $j\in[K]$, with $V_{2K+1}=V_1$. Relabelling cyclically, we may assume that $u_0$ has a neighbour in both $V_1$ and $V_{2K}$. Recall that $V_j^\circ$ is the cluster $V_j$ before the cleaning. Each $V_j$ lost at most $4\varepsilon^{10}m$ vertices in the cleaning and the $t_j\leq3\varepsilon^8n/K$ vertices of $T_j$. It received at most $2\varepsilon^8n/K$ vertices in the integration of $V_0$ and at most $\varepsilon^6n/K$ in the reintegration of $U$. Hence, for each $j\in[2K]$, we have
\[
 |V_j^\circ\setminus V_j|,\ |V_j\setminus V_j^\circ|
 \leq4\varepsilon^{10}m+3\varepsilon^6n/K
 \leq\varepsilon^4m.
\]
The last inequality holds as $m\geq n/(3K)$. Since $V_1^\circ V_2^\circ\ldots V_{2K}^\circ$ is a Hamilton cycle in $R$, each pair $(V_j^\circ,V_{j+1}^\circ)$ is $(2\varepsilon^{10},p)$-regular with $p$-density at least $\alpha$. Let $j\in[2K]$, and let $A\subset V_j$ and $A'\subset V_{j+1}$ have at least $\varepsilon^2m$ vertices each. By the display above, $|A\cap V_j^\circ|,|A'\cap V_{j+1}^\circ|\geq\varepsilon^2m/2\geq2\varepsilon^{10}m$, so regularity gives, as $\alpha>2\varepsilon^{10}$,
\begin{equation}\label{eq:final-pairs}
 e(A,A')\geq e(A\cap V_j^\circ,A'\cap V_{j+1}^\circ)>0.
\end{equation}
In particular, $e(V_j,V_{j+1})>0$ for all $j\in[2K]$, which proves~\ref{cyclic:3}. Moreover, for odd $j$ we have  $\delta(G[V_j, V_{j+1}]) \geq \alpha d / (100K)$ by (ii) of Claim~\ref{cl:randomsubsets} and the reintegration of vertices of $U$.

It remains to show that, for each $j\in[K]$, $H:=G[V_{2j-1},V_{2j}]$ is a bipartite $D$-expander. By~\eqref{eq:final-pairs}, $H$ is $|V_{2j}|/(2D)$-bipartite-joined, as $|V_{2j}|/(2D)\geq\varepsilon^2m$. We show that $V_{2j}$ is $(|V_{2j}|/(2D),D)$-expanding into $V_{2j-1}$, and the same argument applies with $V_{2j-1}$ and $V_{2j}$ exchanged. Fix a set $X\subset V_{2j}$ of size at most $|V_{2j}|/(2D)$. If $|X|\geq\varepsilon^2m$, then by~\eqref{eq:final-pairs} less than $\varepsilon^2m$ vertices of $V_{2j-1}$ have no neighbour in $X$, and thus $|N(X)\cap V_{2j-1}|\geq|V_{2j-1}|-\varepsilon^2m\geq|V_{2j-1}|/2\geq D|X|$. Otherwise, $|X|<\varepsilon^2m\leq\frac{\alpha}{100K(2D+2)}n$. Since every vertex of $V_{2j}$ has at least $\alpha d/(100K)$ neighbours in $V_{2j-1}$, Lemma~\ref{lem:sparse->expander}, with $\alpha/(100K)$ in the role of $\alpha$, implies that the set $V_{2j}$ is $(\frac{\alpha}{100K(2D+2)}n,D)$-expanding into $V_{2j-1}$ in $G$. Hence $|N(X)\cap V_{2j-1}|\geq D|X|$. Therefore, \ref{cyclic:1}--\ref{cyclic:4} hold for $U_j=V_j$ with $j\in [2K]$, where $U_0=\{u_0\}$ if $n$ is odd and $U_0=\emptyset$ if $n$ is even. This completes the proof. \end{proof}

\section{Resilience}\label{sec:resilience}

In this section, we prove our resilience result, Theorem~\ref{thm:sparse_mindegree_cycle}. We will employ Theorem~\ref{thm:cyclic} to find a partition $V(G)=U_0\cup U_1\cup\ldots\cup U_{2K}$ so that, for each $i\in [K]$, $G_i:=G[U_{2i-1}, U_{2i}]$ is a bipartite $D$-expander and there is an edge $e_i$ connecting $U_{2i}$ with $U_{2i+1}$ (working modulo $2K$). By Theorem~\ref{thm:bipartite-c-expander}, for large $D$, a balanced bipartite $D$-expander has a Hamilton path between any two vertices in different parts. So for each $i\in [K]$, we can find a Hamilton path $P_i\subset G_i$ that connects the endpoint of $e_{i-1}$ in $U_{2i-1}$ and the endpoint of $e_{i}$ in $U_{2i}$. These paths $P_1, P_2,\ldots, P_K$ are then put together to form a Hamilton cycle in $G$ (if $n$ is odd, the vertex of $U_0$ is used in place of $e_K$, by~\ref{cyclic:4}).

\begin{theorem}\label{thm:bipartite-c-expander}
    Let $C>0$ be sufficiently large, and let $G$ be a bipartite $C$-expander. Then for all $x,y\in V(G)$ in different parts, there is a Hamilton path in $G$ with $x$ and $y$ as endpoints.
\end{theorem}
  Theorem~\ref{thm:bipartite-c-expander} is a bipartite version of the second statement of Theorem~\ref{thm:C-expander}. We prove it in Appendix~\ref{app:bipartite-connectivity}, following the arguments of~\cite{C-expander}. We write $C_{\ref{thm:bipartite-c-expander}}$ for a constant such that Theorem~\ref{thm:bipartite-c-expander} holds for all $C\geq C_{\ref{thm:bipartite-c-expander}}$.

\begin{proof}[Proof of Theorem~\ref{thm:sparse_mindegree_cycle}]
    Let $D\geq C_{\ref{thm:bipartite-c-expander}}$ and choose $\beta,\eta,d_0$ with
    \[1/d_0 \ll \eta \ll \beta \ll 1/D,\gamma,\]
    so that Theorem~\ref{thm:cyclic} applies. Let $n,d\in\mathbb N$ with $n>d\ge d_0$, and let $G$ be an $(\eta,\beta, d/n)$-sparse graph with $n$ vertices and $\delta(G)\ge (1/2+\gamma)d$.
    Theorem~\ref{thm:cyclic} gives some $K\in\mathbb N$ and a partition $V(G) = U_0 \cup U_1\cup \ldots \cup U_{2K}$ so that~\ref{cyclic:1}--\ref{cyclic:4} hold. For each $i\in[K-1]$ there are vertices $u_{2i}\in U_{2i}$ and $u_{2i+1}\in U_{2i+1}$ with $u_{2i}u_{2i+1}\in E(G)$ by~\ref{cyclic:3}. If $n$ is even, \ref{cyclic:3} also gives $u_{2K}\in U_{2K}$ and $u_1\in U_1$ with $u_{2K}u_1\in E(G)$, and we let $Q=u_{2K}u_1$. If $n$ is odd, by~\ref{cyclic:4} the vertex $u_0$ of $U_0$ has neighbours $u_1\in U_1$ and $u_{2K}\in U_{2K}$, and we let $Q=u_{2K}u_0u_1$. For $i\in [K]$, as $G[U_{2i-1},U_{2i}]$ is a bipartite $D$-expander by \ref{cyclic:2} and $D\ge C_{\ref{thm:bipartite-c-expander}}$, we can use \Cref{thm:bipartite-c-expander} to find a Hamilton path $P_i$ in $G[U_{2i-1},U_{2i}]$ with endpoints $u_{2i-1}$ and $u_{2i}$. Then $u_1P_1u_2u_3P_2u_4\ldots u_{2K-1}P_Ku_{2K}$ followed by $Q$ is a Hamilton cycle in $G$.
\end{proof}

 %------------------------------------------------------

\section{Hamilton connectivity via sparse augmentation}\label{sec:Hamilton-conn-bipartite}
As outlined in the introduction, our strategy for proving the packing result (Theorem~\ref{thm: decomp}) is to extract Hamilton cycles one by one while maintaining the expansion properties in the remaining graph. To do so, we split the edges of the initial graph into two parts: a graph $G_1$, which will have local expansion, and a graph $G_2$, which will have many edges between every pair of disjoint sets of size at least $\beta n$. While extracting Hamilton cycles, we aim to use as few $G_2$-edges as possible, and, more importantly, to control the maximum degree of the subgraph of $G_2$ formed by the edges we do use. By doing this, we can ensure that, throughout the process, $G_2$ will have edges between any two disjoint sets of size $\beta n$, and $G_1$ will have good local expansion properties, as at each step the degree in $G_1$ drops by at most 2.
Given a graph $H$ (which we think of as the set of $G_2$-edges already used), the main result of this section (Lemma~\ref{lem:regularHamcycle}) guarantees the existence of a Hamilton path using the remaining edges, whose $G_2$-edges ideally avoid the high-degree vertices of $H$. If avoiding them altogether is not possible, then the goal becomes to ensure that the maximum degree of $H$ remains low. To formalize this control, we introduce the notion of the $m$-max degree of $H$ in the definition below.

\begin{definition}\label{def:degree}
    For a graph $H$, let $t(H)$ be the number of vertices of $H$ with degree exactly $\Delta(H)$ (so $t(H)=|V(H)|$ if $H$ has no edges). For an integer $m\ge 1$, the \emph{$m$-max degree} of $H$ is defined as $\Delta_m(H)=m\cdot \Delta(H) + \min\{t(H),m\}$, and we say a collection $E$ of pairs of vertices of $H$ is \emph{$(i,m,H)$-bounded} if $\Delta_m(H+E) \leq \Delta_m(H) + i$.
\end{definition}

Note that  
$\Delta(H)\leq \Delta_{m}(H)/m \leq \Delta(H) + 1$ always holds. Since the $m$-max degree of $H$ is monotone under adding edges, if a collection $E$ of vertex pairs is $(i,m,H)$-bounded, then every subset of $E$ is also $(i,m,H)$-bounded.
Furthermore, there always exists a set $Z\subset V(H)$ of size at most $m$ such that for every edge $e$ not incident to any vertex in $Z$, we have $\Delta_m(H + e) \leq \Delta_m(H) + 2$. Indeed, if the set of maximum-degree vertices has size at most $m$, it can serve as such a set $Z$; otherwise, $Z=\emptyset$ suffices. We will also use the fact that if $E_1$ is $(i_1,m,H)$-bounded and $E_2$ is $(i_2,m,H+E_1)$-bounded, then $E_1\cup E_2$ is $(i_1+i_2,m,H)$-bounded.
These facts are proved in Lemma~\ref{lem:max-degree}.
When $m$ is written as $\beta n$ or $5\beta n$, it is computed with the original value of $n$, and it does not change when we restrict to subgraphs.
\par Now we state the main result of the section, which we will prove after we collect several supporting results.
\begin{restatable}{lemma}{RegularHamCycleLemma}\label{lem:regularHamcycle}
Let $1/d_0\ll\eta\ll\beta\ll\alpha\ll1$ and $1/d_0\ll\gamma\ll\alpha$, and let $n\in\mathbb N$ and $d\geq d_0$. Let $G_1$ and $G_2$ be bipartite graphs with parts $V_1,V_2$ of size $n$, where $G_1$ is $(\eta,\beta,d/n)$-sparse with $\delta(G_1)\geq\alpha d$ and $\Delta(G_1)\le(1+\gamma)\delta(G_1)$, and $G_2$ is $\beta n$-bipartite-joined. Let $H$ be any graph with $V(H)\supset V_1\cup V_2$.
Then, for all $x_1 \in V_1$ and $x_2 \in V_2$, the graph $G_1 \cup G_2$ contains a Hamilton path $P$ with endpoints $x_1$ and $x_2$ such that $E(P)\setminus E(G_1)$ is $(10^6 \gamma \log(1/\gamma) n, \beta n, H)$-bounded.
\end{restatable}
 We follow the approach of~\cite{C-expander}, which uses a linking structure and rotations of linear forests. Here $G_1$ expands only small sets and need not be $\beta n$-bipartite-joined, and the edges of $G_2$ that we use must form a bounded set as in Lemma~\ref{lem:regularHamcycle}. So we rotate inside sets that expand into each other in $G_1$, and we use edges of $G_2$ only when needed. Many of the results in this section are suitable adaptations of auxiliary results from~\cite{C-expander}, and we prove those in the appendix. Following steps~\ref{step:1}--\ref{step:3} of Section~\ref{sec:overview}, we find the linking structure in Section~\ref{sec:linking} (Lemma~\ref{lem: partition}) and the linear forest in Section~\ref{sec:forests} (Lemmas~\ref{lemma:linearforest decomp1} and~\ref{lemma:linearforest decomp2}), and we join them in Section~\ref{sec:G1G2} using Lemma~\ref{lem:cycle-assembly}.

\subsection{Linking structures in bipartite expanders}\label{sec:linking}
The first supporting result we need is the existence of linking structures in bipartite expanders. The notion of linking structures in sparse expanders was first used by Hyde, Morrison, M\"{u}yesser, and Pavez-Sign\'{e}~\cite{hyde2023spanning} to find spanning trees in $(n,d,\lambda)$-graphs, and it was then used by Dragani\'c, Montgomery, Munh\'{a} Correia, Pokrovskiy, and Sudakov~\cite{C-expander} to solve the Hamiltonicity problem in expander graphs. Here, we use the terminology introduced in~\cite{C-expander}, but we do not require the paths $P_a$ to have equal length.
\begin{definition}\label{def: linking structure}
Let $H$ be a graph and let $A,B\subset V(H)$ be disjoint sets with $|A|=|B|$. Say $H$ is an \textit{$(A,B)$-linking structure} if for every bijection $\varphi:A\to B$ there are vertex-disjoint paths $P_a$, $a\in A$, which partition $V(H)$ and such that $P_a$ has endpoints $a$ and $\varphi(a)$. We call the vertices of $A\cup B$ the \textit{terminals} of the linking structure.
\end{definition}
We use a linking structure to join the paths of a linear forest into one path, as in step~\ref{step:3} of Section~\ref{sec:overview}.
\begin{lemma}\label{lem:cycle-assembly}
Let $L$ be an $(A\cup\{x_1\},B\cup\{x_2\})$-linking structure with $x_1,x_2\notin A\cup B$, and let $\mathcal F$ be a linear forest with no isolated vertices such that $V(\mathcal F)\cap V(L)=End(\mathcal F)=A\cup B$. Then $L\cup\mathcal F$ has a spanning $x_1,x_2$-path containing every edge of $\mathcal F$.
\end{lemma}
\begin{proof}
If $A=\emptyset$, then $\mathcal F$ is empty and this is the definition of a linking structure. Otherwise let $A'=A\cup\{x_1\}$ and $B'=B\cup\{x_2\}$, and let $M_{\mathcal F}$ be the perfect matching on $A'\cup B'$ formed by the pair $x_1x_2$ and the pairs of endpoints of the paths of $\mathcal F$. Among all $A'$--$B'$ perfect matchings $M$, we choose one minimizing the number of components of the multigraph $M\cup M_{\mathcal F}$. The components of $M\cup M_{\mathcal F}$ are alternating cycles, where a pair lying in both $M$ and $M_{\mathcal F}$ gives a cycle of length two. If there are at least two components, choose $ab\in M$ in one and $a'b'\in M$ in another, where $a,a'\in A'$. Replacing $ab,a'b'$ by $ab',a'b$ merges these two cycles, which contradicts the choice of $M$.

Since $L$ is an $(A',B')$-linking structure, it has vertex-disjoint paths that partition $V(L)$ and join the pairs of $M$. The cycle $M\cup M_{\mathcal F}$ has at least four vertices since $|A'|\geq2$, so $M$ and $M_{\mathcal F}$ share no pair. The paths of $L$ and of $\mathcal F$ have disjoint interiors, as $V(\mathcal F)\cap V(L)=A\cup B$. Replacing each edge of the cycle other than $x_1x_2$ by its path therefore gives a spanning $x_1,x_2$-path containing all edges of $\mathcal F$.
\end{proof}
In our setting, we will need to find a linking structure in the union of two bipartite graphs $G_1$ and $G_2$ so that it uses only a few edges from $G_2$.
Furthermore, the terminals of the linking structure will be paired up and connected by disjoint paths; these paths will also contain a prescribed linear forest, which might be of linear size. 
We now state the lemma which captures this; the statement is slightly technical, so we discuss it shortly after.

\begin{lemma}\label{lem: partition}
Let $D\geq8$ be an integer, and let $\rho>0$, $0<\beta\leq\rho/10^8$, and $0<\delta<1/10$ satisfy $\delta\log(e/\delta)\leq\rho/10^4$. Let $n$ be sufficiently large in terms of these constants, and set $b=\beta n$, $r=1000b/D$, and $t=n^{0.9}$. Let $G_1,G_2$ be bipartite graphs with parts $V_1,V_2$ of size at most $n$. Let $U_1,U_2$ be balanced sets partitioning $V_1\cup V_2$. Let $H$ be any graph with $V(H)\supset V_1\cup V_2$. Fix $x_1\in U_2\cap V_1$ and $x_2\in U_2\cap V_2$. Suppose that the following holds.
\begin{enumerate}[label=\upshape(\roman{enumi})]
\item\label{eq: partition condition 1} For $i\in[2]$ and nonempty $U\subset V_i$ with $|U|\leq2r$, we have $|\Gamma_{G_1}(U)\cap U_2|\geq2D|U|$, and $|\Gamma_{G_1}(U)\cap U_2|\geq\rho n$ if $|U|\geq r$.

\item\label{eq: partition condition 2} $G_1[U_1]$ has a spanning linear forest $\mathcal F$ with $\ell$ paths, where $t\leq\ell\leq\delta n$, each of length at most $n^{0.15}$.

\item\label{eq: partition condition 3} $G_2$ is $b$-bipartite-joined.
\end{enumerate}
Then there are a graph $F\subset G_1\cup G_2$, balanced sets $X,Y,Z$ partitioning $U_2$, and sets $A\subset X\cap V_1$ and $B\subset X\cap V_2$, such that the following holds (see Figure~\ref{fig:linking-output}).
\begin{enumerate}[label=\upshape\textbf{\Alph{propcounter}\arabic{enumi}}]
\item\label{eq: partition 1} $F$ spans $U_1\cup X\cup Y$, contains every edge of $\mathcal F$, and $|X\cup Y|\leq200\delta\log(e/\delta)n$.

\item\label{eq: partition 2} $x_1,x_2\in X\setminus(A\cup B)$, $|A|=|B|=t$, and $F[X]$ is an $(A\cup\{x_1\},B\cup\{x_2\})$-linking structure.

\item\label{eq: partition 3} $F[U_1\cup Y\cup A\cup B]$ is a forest of $t$ paths, each of length at most $3n^{0.15}$, with one endpoint in $A$ and one in $B$.

\item\label{eq: partition 4} $V_1\cup V_2$ is $(2r,D/2)$-expanding into $Z$ in $G_1$.

\item\label{eq: partition 5} $E(F)\setminus E(G_1)$ is $(10\delta n,b,H)$-bounded.
\end{enumerate}
\stepcounter{propcounter}
\end{lemma}

\begin{figure}[htbp]
  \centering
  \resizebox{0.7\linewidth}{!}{%
  \begin{tikzpicture}[
      x=1cm, y=.75cm,
      font=\small,
      region/.style={
          rounded corners=4pt,
          draw=black!40,
          line width=.45pt
      },
      vertex/.style={
          circle,
          fill=black!85,
          inner sep=0pt,
          minimum size=3.5pt
      },
      retained/.style={
          draw=black!85,
          line width=1.7pt,
          line cap=round
      },
      paths/.style={
          draw=blue!55!black,
          line width=.9pt,
          rounded corners=9pt,
          line cap=round
      },
      lab/.style={inner sep=2pt}
  ]

  \draw[region, dashed, dash pattern=on 3pt off 2pt]
      (-.4,-.15) rectangle (13.65,6.65);

  \path[region, fill=yellow!8]
      (0,3.55) rectangle (9.7,6.25);
  \path[region, fill=blue!3]
      (0,.25) rectangle (9.7,3.2);
  \path[region, fill=black!3]
      (10.25,.25) rectangle (13.2,6.25);
  \path[region, fill=black!2]
      (0,-3.15) rectangle (9.7,-.65);

  \draw[draw=black!35, line width=.4pt, fill=white]
      (2.3,3.97) ellipse [x radius=1.6cm, y radius=.23cm];
  \draw[draw=black!35, line width=.4pt, fill=white]
      (7.4,3.97) ellipse [x radius=1.6cm, y radius=.23cm];

  \draw[paths]
      (1.1,3.97) -- (1.1,-2.65) -- (8.6,-2.65) -- (8.6,3.97);
  \draw[paths]
      (2.3,3.97) -- (2.3,-1.9) -- (7.4,-1.9) -- (7.4,3.97);
  \draw[paths]
      (3.5,3.97) -- (3.5,-1.15) -- (6.2,-1.15) -- (6.2,3.97);

  \foreach \xx in {1.1,2.3,3.5,6.2,7.4,8.6}
      \node[vertex, fill=blue!55!black, minimum size=2.6pt]
          at (\xx,1.1) {};

  \draw[retained] (1.95,-2.65) -- (4.4,-2.65);
  \draw[retained] (5.0,-2.65) -- (7.75,-2.65);
  \draw[retained] (3.0,-1.9) -- (6.7,-1.9);
  \draw[retained] (4.15,-1.15) -- (5.55,-1.15);

  \foreach \p in {
      (1.95,-2.65),(4.4,-2.65),(5.0,-2.65),(7.75,-2.65),
      (3.0,-1.9),(6.7,-1.9),(4.15,-1.15),(5.55,-1.15)
  }
      \node[vertex, minimum size=2.6pt] at \p {};

  \node[vertex] at (4.45,3.97) {};
  \node[vertex] at (5.25,3.97) {};
  \node[lab, anchor=south] at (4.45,4.04) {$x_1$};
  \node[lab, anchor=south] at (5.25,4.04) {$x_2$};

  \foreach \p in {
      (1.1,3.97),(2.3,3.97),(3.5,3.97),
      (6.2,3.97),(7.4,3.97),(8.6,3.97)
  }
      \node[vertex] at \p {};

  \node[lab, fill=white] at (12.7,6.65) {$U_2$};

  \node[lab, anchor=north west] at (.2,6.08) {$X$};
  \node[lab, align=center] at (5.6,5.45)
      {$F[X]$\\[-1pt]linking structure};

  \node[lab] at (2.3,4.65) {$A$};
  \node[lab] at (7.4,4.65) {$B$};

  \node[lab, anchor=north west] at (.2,3.04) {$Y$};
  \node[lab, align=center, text=blue!55!black] at (4.85,1.72)
      {paths joining\\$A$ to $B$};

  \node[lab, align=center] at (11.72,3.45)
      {$Z$\\[4pt]unused by $F$};

  \node[lab, anchor=north west] at (.2,-.82) {$U_1$};

  \draw[retained] (2.8,-3.8) -- (3.65,-3.8);
  \node[lab, anchor=west] at (3.8,-3.8)
      {original forest $\mathcal F$};

  \end{tikzpicture}%
  }
  \caption{The sets and the graph $F$ given by Lemma~\ref{lem: partition}.
  The linking structure $F[X]$ and the paths of~\ref{eq: partition 3} intersect in $A\cup B$,
  and the paths contain the forest $\mathcal F$ (bold segments).
  The edges of $F[X]$ are not drawn; it is an $(A\cup\{x_1\},B\cup\{x_2\})$-linking structure.}
  \label{fig:linking-output}
  \end{figure}
We prove Lemma~\ref{lem: partition} in Appendix~\ref{app:linking}. The second bound in~\ref{eq: partition condition 1} ensures that each set $U$ in one part with $r\le|U|\le2r$ keeps many neighbours in $U_2$ outside the paths that we build.
We now explain how Lemma~\ref{lem: partition} is used in the proof of Lemma~\ref{lem:regularHamcycle}. We start with conditions~\ref{eq: partition condition 1} and~\ref{eq: partition condition 2}, noting that  condition~\ref{eq: partition condition 3} is exactly the hypothesis on $G_2$ in Lemma~\ref{lem:regularHamcycle}. First, in the proof of Lemma~\ref{lem:regularHamcycle} we split $V(G_1)$ into three balanced sets $U_1,U_2,U_3$ so that each vertex has about a third of its neighbours in each of them, and we apply Lemma~\ref{lem: partition} inside $U_1\cup U_2$. Then condition~\ref{eq: partition condition 1} follows from the hypotheses on $G_1$ and Lemma~\ref{lem:sparse->expander}. Second, as explained in Section~\ref{sec:overview}, condition~\ref{eq: partition condition 2} can be achieved by using Lemma~\ref{lem:pathpartition} and then cutting the resulting paths into shorter ones, giving us the linear forest $\mathcal F$.

By~\ref{eq: partition 2}, the subgraph $F[X]$ of $F$ is an $(A\cup\{x_1\},B\cup\{x_2\})$-linking structure. The paths from $\mathcal F$ are joined together efficiently through the set $Y$ so that $F[U_1\cup Y\cup A\cup B]$ has a spanning linear forest with $|A|=|B|=n^{0.9}$ paths whose endpoints are precisely $A\cup B$. This is property~\ref{eq: partition 3}. Moreover, no path in $\mathcal F$ is broken, and only a few new vertices are used along the connecting procedure; this is property~\ref{eq: partition 1}. Property~\ref{eq: partition 5} says that the edges from $G_2$ that we use while constructing $F$ form a $(10\delta n,b,H)$-bounded edge set. This is one of the new ingredients we need, for which we carefully use the \textit{extendability method}, an embedding technique introduced by Friedman and Pippenger~\cite{FP1987} to embed trees in expanders. Lastly, $F$ does not use the vertices of $Z$, so that later we can find a spanning linear forest on $Z$ and the vertices outside $U_1\cup U_2$, which we combine with the paths in~\ref{eq: partition 3} and then rotate to obtain a spanning linear forest on $V(G_1)\setminus(X\setminus(A\cup B))$ with endpoint set $A\cup B$, while keeping the set of $G_2$-edges used bounded (this last step is similar to that in~\cite{C-expander}, but with the additional consideration of how to rely on $G_2$ so that no bipartite hole is produced while rotating the forest). This final step will then be enough, together with the properties of the linking structure, in order to produce a Hamilton path in $G_1\cup G_2$ with endpoints $x_1$ and $x_2$ while controlling how the edges of $G_2$ are used. Property~\ref{eq: partition 4} is not needed for Lemma~\ref{lem:regularHamcycle}; we use it in the proof of Theorem~\ref{thm: covering} and in Appendix~\ref{app:bipartite-connectivity}.

%----------------------------------------------------

\subsection{Removal of endpoints in a linear forest}\label{sec:forests}

Using Lemma~\ref{lemma:linearforest decomp1}, we reduce the number of paths of a spanning linear forest to at most $n^{0.8}$, and then, using Lemma~\ref{lemma:linearforest decomp2}, we move all endpoints into $A\cup B$. We prove both lemmas in Appendix~\ref{app:forests}.
In Lemmas~\ref{lemma:linearforest decomp1} and~\ref{lemma:linearforest decomp2}, let $1/d_0\ll\eta,\beta\ll\alpha\ll1$ and let $n\in\mathbb N$ and $d\geq d_0$. Let $G_1$ be an $(\eta,\beta,d/n)$-sparse bipartite graph with parts $V_1,V_2$ of size $n$ and $\delta(G_1)\geq\alpha d$, let $G_2$ be a $\beta n$-bipartite-joined graph with parts $V_1,V_2$, and let $H$ be any graph with $V(H)\supset V_1\cup V_2$.

\begin{lemma}\label{lemma:linearforest decomp1}
Given a spanning linear forest $\mathcal{F}\subset G_1$ with $k$ paths, there is a $(400k, \beta n, H)$-bounded set $E\subset E(G_2)$ and a spanning linear forest $\mathcal{F}'\subset G_1+E$ with at most $n^{0.8}$ paths and no isolated vertices.
\end{lemma}
\begin{lemma}\label{lemma:linearforest decomp2}
Suppose $\mathcal F=\mathcal P\cup\mathcal Q$ is a spanning linear forest on $V_1\cup V_2$ with no isolated vertices, all of whose edges join $V_1$ and $V_2$, where
\begin{enumerate}[label=\upshape(\roman{enumi})]
\item $\mathcal P$ has $t\leq n^{0.95}$ paths of lengths at most $n^{0.19}$, and $End(\mathcal P)=A\cup B$ with $A\subset V_1$, $B\subset V_2$ and $|A|=|B|=t$,
\item $R=V(\mathcal P)$ satisfies $|R\cap V_1|\geq n/10$ and $\delta(G_1[R])\geq\alpha d$, and
\item $\mathcal Q$ has $q\leq n^{0.8}$ paths.
\end{enumerate}
Then there is a $(10q,\beta n,H)$-bounded set $E\subset E(G_2)$ and a spanning linear forest $\mathcal F'\subset G_1+\mathcal F+E$ with no isolated vertices and $End(\mathcal F')=A\cup B$.
\end{lemma}

%------------------------------------------------

\subsection{Proof of the main lemma}\label{sec:G1G2}
\begin{proof}[Proof of Lemma~\ref{lem:regularHamcycle}]
We split $V(G_1)$ into sets $U_1,U_2,U_3$ and use Lemma~\ref{lem: partition} to find a linking structure $F[X]$ with $X\subset U_2$ and paths from $A$ to $B$ covering $U_1$. We cover the remaining vertices by paths ending in $A\cup B$ using Lemmas~\ref{lemma:linearforest decomp1} and~\ref{lemma:linearforest decomp2}, and then we join these paths and the linking structure into a Hamilton $x_1,x_2$-path by Lemma~\ref{lem:cycle-assembly}.
Since $\delta(G_1)\geq\alpha d$, we have $n\geq\alpha d\geq\alpha d_0$, and thus $n$ is sufficiently large. By sparsity, we have 
\[\delta(G_1)|V_1|\le e(V_1,V_2)\le (1+\beta)d|V_1|\]
and thus $\Delta(G_1)\le (1+\gamma) \delta(G_1)\le (1+\beta)(1+\gamma)d\le 1.1d$ as $\beta,\gamma\ll 1$.

\begin{claim}\label{cl: three-U}
There is a partition $V(G_1)=U_1\cup U_2\cup U_3$ such that the following holds.    
\begin{enumerate}[label= \upshape\textbf{\Alph{propcounter}\arabic{enumi}}]
    \item Each $U_i$ is a balanced set of size $\tfrac{2n}{3}\pm 4$ and $x_1,x_2\in U_2$.\label{eq: partition into six 1}
    \item For every vertex $v\in V(G_1)$ and $i\in [3]$, 
    $|N_{G_1}(v)\cap U_i| = \tfrac{1}{3} d_{G_1}(v) \pm 2d^{2/3} \geq \alpha d/4$. \label{eq: partition into six 2}
\end{enumerate}\end{claim}
\begin{claimproof}[Proof of Claim~\ref{cl: three-U}]
Let $L$ consist of $x_1,x_2$ and the same number, at most two, of further vertices in each part, chosen so that each part has a multiple of three vertices outside $L$. Partition $W=V(G_1)\setminus L$ into triples, each lying in $V_1$ or in $V_2$. As $\Delta(G_1)\le1.1d$, by Lemma~\ref{lem:balanced-random-splitting} applied to $G_1$ with $k=3$ there is a partition $W=U_1'\cup U_2'\cup U_3'$ such that each $U_i'$ meets each triple in one vertex and $|N_{G_1}(v)\cap U_i'|=\tfrac13|N_{G_1}(v)\cap W|\pm d^{2/3}$ for all $v\in V(G_1)$ and $i\in[3]$. Let $U_1=U_1'$, $U_2=U_2'\cup L$ and $U_3=U_3'$. For~\ref{eq: partition into six 1}, each $U_i$ is balanced, as $L$ is, and $|U_i|=\frac{2n}3\pm4$ as $|U_i'|=\frac23(n-1-j)$, where $j\le2$ is the number of vertices of $L\setminus\{x_1,x_2\}$ in each part. Moreover, as $|L|\le6$, we have $|N_{G_1}(v)\cap U_i|=\frac13d_{G_1}(v)\pm(d^{2/3}+6)$, and~\ref{eq: partition into six 2} follows as $d\ge d_0$ is large.
\end{claimproof}

By~\ref{eq: partition into six 2}, each vertex has at least $\bar d:=\frac{1}{3}\delta(G_1)-2d^{2/3}$ and at most $\frac{1}{3}(1+\gamma)\delta(G_1)+2d^{2/3}\le (1+2\gamma)\bar d$ neighbours in each $U_i$. Thus all degrees in $G_1[U_1]$ lie between $\bar d$ and $(1+2\gamma)\bar d$, and those in $G_1[U_2\cup U_3]$ between $2\bar d$ and $2(1+2\gamma)\bar d$. We may apply~\Cref{lem:pathpartition} to both graphs, with $2\gamma$ playing the role of $\eps$, to obtain spanning linear forests $\mathcal F^0_1$ in $G_1[U_1]$ and $\mathcal F^0_2$ in $G_1[U_2\cup U_3]$. They have at most $6\gamma|U_1|\le5\gamma n$ and $6\gamma|U_2\cup U_3|\le9\gamma n$ paths, respectively. Moreover, by cutting the paths in $\mathcal F^0_1$ into paths of length at most $n^{0.1}/2$, we may assume that $\mathcal F_1^0$ consists of $\ell$ paths each of length at most $n^{0.1}/2$, where $n^{0.9}\le \ell\le 5\gamma n+2n^{0.9}\le 6\gamma n$. Set $D=8$, $b=\beta n$, $r=1000b/D$ and $\rho=\alpha/100$. Lemma~\ref{lem:sparse->expander}, with $2n$, $2d$, $\alpha/8$ and $2D$ playing the roles of $n$, $d$, $\alpha$ and $D$, together with~\ref{eq: partition into six 2} for $U_2$, implies condition~\ref{eq: partition condition 1} of Lemma~\ref{lem: partition}, as $2\eta n\leq r$ and $2r\leq\alpha n/200$ by $\eta\ll\beta\ll\alpha$. Finally, $6\gamma\log(e/6\gamma)\leq\rho/10^4$ and $\beta\leq\rho/10^8$ as $\beta,\gamma\ll\alpha$. Thus we can apply \Cref{lem: partition} in $G_1[U_1\cup U_2]$ and $G_2[U_1\cup U_2]$, with $6\gamma$ playing the role of $\delta$, to obtain a partition $U_2=X\cup Y\cup Z$, a subgraph $F\subset (G_1\cup G_2)[U_1\cup U_2]$, and sets $A,B\subset X$ so that \ref{eq: partition 1}--\ref{eq: partition 5} hold. By~\ref{eq: partition 3}, $\mathcal F_1^1:=F[U_1\cup Y\cup A\cup B]$ is a linear forest of $n^{0.9}$ paths with $End(\mathcal F_1^1)=A\cup B$.
Moreover, the set $E_1=E(F)\setminus E(G_1)$ is $(60\gamma n, \beta n, H)$-bounded.

We will apply Lemmas~\ref{lemma:linearforest decomp1} and~\ref{lemma:linearforest decomp2} inside balanced sets $S\supset U_3$. For such $S$, set $m_S=|S|/2$ and $\beta_S=b/m_S$. We have $\beta\leq\beta_S\leq4\beta$, since $n/3-2\leq m_S\leq n$ and $n$ is large. Every vertex has at least $\alpha d/4$ neighbours in $U_3\subset S$ by~\ref{eq: partition into six 2}, so $\delta(G_1[S])\geq\alpha d/4$. By Lemma~\ref{lem:restriction-normalization} and as $\beta\leq\beta_S$, $G_1[S]$ is $(\eta,\beta_S,d/m_S)$-sparse, and $G_2[S]$ is $\beta_Sm_S$-bipartite-joined as $\beta_Sm_S=b$. The hypotheses stated before the two lemmas hold for $G_1[S]$ and $G_2[S]$, with $m_S$, $\beta_S$ and $\alpha/4$ playing the roles of $n$, $\beta$ and $\alpha$, as $\beta_S\leq4\beta$ and $\beta\ll\alpha$.

Let $\mathcal{F}^1_2=\mathcal F^0_2-(X\cup Y)$, a spanning linear forest in $G_1[Z\cup U_3]$. As each removed vertex adds at most one new path and $|X\cup Y|\le 1200\gamma\log(e/\gamma)n$ by~\ref{eq: partition 1}, $\mathcal{F}_2^1$ has at most $9\gamma n+|X\cup Y|\leq1212\gamma\log(e/\gamma)n$ paths. The set $Z\cup U_3$ is balanced; set $m_0=|Z\cup U_3|/2$. We apply Lemma~\ref{lemma:linearforest decomp1} to this forest with $S=Z\cup U_3$ and $H+E_1$ in place of $H$. Since $400\cdot1212\log(e/\gamma)\leq600000\log(1/\gamma)$ as $\gamma\ll\alpha$, we obtain a $(600000\gamma\log(1/\gamma)n, \beta n,H+E_1)$-bounded set $E_2\subset E(G_2)$ and a spanning linear forest $\mathcal{F}'\subset G_1[Z\cup U_3]+E_2$ with at most $m_0^{0.8}\leq n^{0.8}$ paths and no isolated vertices.

Let $\mathcal F=\mathcal F'\cup \mathcal F_1^1$, let $W=V(\mathcal F)$ and note that $A\cup B\subset W$ as $A\cup B= End(\mathcal F_1^1)$. Our goal now is to use Lemma~\ref{lemma:linearforest decomp2} to reduce the number of endpoints of $\mathcal F$. The set $W$ is balanced, as its complement is $X\setminus(A\cup B)$, and it contains $U_3$. Let $m=|W|/2$. Note that $\mathcal F$ is a spanning linear forest in $(G_1\cup G_2)[W]$ with parts $V_1\cap W$ and $V_2\cap W$, and has no isolated vertices, as $\mathcal F'$ has no isolated vertices and each path in $\mathcal F_1^1$ has one endpoint in $A$ and one in $B$. To check the hypotheses of Lemma~\ref{lemma:linearforest decomp2}, let $\mathcal P=\mathcal F_1^1$ and $\mathcal Q=\mathcal F'$, so that $R=V(\mathcal P)=U_1\cup Y\cup A\cup B$. By~\ref{eq: partition 3}, $\mathcal P$ has $t=n^{0.9}$ paths of lengths at most $3n^{0.15}$. As $U_1\cup U_3\subset W$, we have $m\geq n/2$, and thus $t\leq m^{0.95}$ and $3n^{0.15}\leq m^{0.19}$. The set $R$ contains $U_1$ and has at least $n/3-2\geq m/10$ vertices in each part. Every vertex has at least $\alpha d/4$ neighbours in $U_1$, so $\delta(G_1[R])\geq\alpha d/4$. Also $\mathcal Q$ has $q\leq m_0^{0.8}\leq m^{0.8}$ paths. Therefore we may apply Lemma~\ref{lemma:linearforest decomp2} with $S=W$ and $H+E_1+E_2$ in place of $H$, which yields a $(10q,\beta n,H+E_1+E_2)$-bounded set $E_3\subset E(G_2[W])$ and a spanning linear forest $\mathcal F^*\subset G_1+\mathcal F+E_3$ with no isolated vertices and $End(\mathcal F^*)=A\cup B$.

To conclude, note that $X\cap W=A\cup B$ and $X\cup W=V(G_1)$. Since $\mathcal F^*$ has no isolated vertices and endpoint set $A\cup B$, by~\ref{eq: partition 2} and Lemma~\ref{lem:cycle-assembly} there is a spanning $x_1,x_2$-path $P$ of $F[X]\cup\mathcal F^*$, and $P$ is a Hamilton path.
Every edge of $P$ outside $G_1$ belongs to $E_1\cup E_2\cup E_3$. The bounds above give
\[
 60\gamma n+600000\gamma\log(1/\gamma)n+10n^{0.8}
 \leq10^6\gamma\log(1/\gamma)n
\]
for sufficiently large $n$. By the facts stated after Definition~\ref{def:degree}, $E(P)\setminus E(G_1)$ is therefore $(10^6\gamma\log(1/\gamma)n,\beta n,H)$-bounded.
\end{proof}
%------------------------------------------------

\section{Packing and covering}\label{sec:packing-covering}
Suppose $G$ is as in Theorem~\ref{thm: decomp} and we want to obtain an approximate decomposition of $G$ into Hamilton cycles, that is, we want to find $(1-\varepsilon)d/2$ edge-disjoint Hamilton cycles in $G$ for a given $\eps\in(0,1)$. A natural approach is to find and then delete Hamilton cycles from $G$ one by one until only $O(\varepsilon d n)$ edges remain. Using Theorem~\ref{thm:sparse_mindegree_cycle}, we can safely remove $(1/4 - \varepsilon)d$ edge-disjoint Hamilton cycles, as at each step the degree of each vertex drops by exactly 2. After this point, however, there is no guarantee that we can continue, as the minimum degree may drop below $d/2$. 
As mentioned in Section~\ref{sec:overview}, the key obstacle is that after linearly many cycles have been removed from $G$, the remaining subgraph may contain two large disjoint vertex sets with no edge between them. To address this issue, we first  partition $G$ into two random subgraphs $G=G_1\cup G_2$ so that a)~$G_1$ is almost regular with degrees close to $d$, and b) $G_2$ is `robustly $5\beta n$-joined' for some $0<\beta\ll \varepsilon$, that is, $G_2$ has the property that between any two disjoint subsets of size at least $5\beta n$ there are at least $\beta^3 dn$ edges. Our goal will then be to find Hamilton cycles $C_1, C_2, \dots, C_{(1 - \eps)d/2}$ in $G_1 \cup G_2$ so that, for each $i\in [(1-\varepsilon)d/2]$ in turn, $C_i$ uses only a few edges from $G_2$ so that the unused edges of $G_2$ still form a $5\beta n$-joined graph throughout the process. This plan is viable because of Lemma~\ref{lem:regularHamcycle} from the previous section, which produces a Hamilton path with prescribed endpoints, and hence a Hamilton cycle $C_i$, while controlling the distribution of $E(C_i)\cap E(G_2)$. This allows us to keep finding cycles until we have found $(1-\varepsilon)d/2$ of them, giving the desired packing.

\subsection{The packing result}
\begin{proof}[Proof of Theorem~\ref{thm: decomp}]
It suffices to prove the statement for sufficiently small $\eps>0$. We choose constants with $1/d_0\ll\eta,\gamma\ll\beta\ll\eps\ll1$, and set $\alpha=\eps/10$.

First, we partition $G$ randomly into two edge-disjoint subgraphs $G_1$ and $G_2$ such that
\begin{itemize}
       \item  $d_{G_1}(v)= (1-\alpha\pm 2\gamma) d$ for all $v\in V(G)$, and \label{eq: G121}
        \item for all disjoint sets $X,Y\subset V(G)$ with $|X|,|Y|\geq 5\beta n$, we have $e_{G_2}(X,Y)\geq \beta^3 d n$. \label{eq: G122}
\end{itemize}
Indeed, form $G_1$ and $G_2$ so that each edge $e\in E(G)$ is added to $G_1$ with probability $1-\alpha$, and to $G_2$ with probability $\alpha$, and all choices are independent. For each vertex $v\in V(G)$, let $X_v$ be the event that $d_{G_1}(v)\neq  (1-\alpha\pm 2\gamma)d$. For a fixed vertex $v\in V(G)$,  note that Chernoff's bound (Lemma~\ref{lem:Chernoff}) implies that 
\[\mathbb P(X_v)\le \mathbb P\left(|d_{G_1}(v)-\mathbb E[d_{G_1}(v)]|>\gamma d\right)\le e^{ -\gamma^{2} d/4}.\]
Moreover, $X_v$ is independent of all other $X_u$ except for at most $2d$ choices of $u\in N_G(v)$. Then, since $(2d+1) e^{1-\gamma^2 d/4} <1$ holds for large enough $d$,  Lemma~\ref{Lem:LLL} gives that $\mathbb P(\cap_{v\in V(G)}\overline{X_v})\ge (1- \frac{1}{2d+1})^n \geq e^{-n/d}$. To check the second bullet point, note that for each pair of disjoint sets $X,Y\subset V(G)$ with $|X|=|Y|= 5\beta n$, as $G$ is $(\eta,\beta,d/n)$-sparse we have
\begin{eqnarray*}e_{G}(X,Y)\ge (1-\gamma)d|X| - e_G(X,V(G)\setminus Y) &\geq&  (1-\gamma)d|X| - (1+\beta)\frac{d|X|(n-5\beta n)}{n} \\
&\geq& (4\beta+5\beta^2-\gamma )d|X|\geq 15\beta^2 d n,\end{eqnarray*}
as $\gamma\ll\beta$. Then $\mathbb E[e_{G_2}(X,Y)]\ge 15\alpha \beta^2dn\ge 2\beta^3dn $ as $\beta\ll\eps$, and hence, by Chernoff's bound, $e_{G_2}(X,Y)<\beta^3dn$ with probability at most $e^{-\beta^3dn/12}$. By a union bound over all choices of $X,Y$, the second bullet point fails with probability at most $2^{2n}e^{-\beta^3dn/12}<e^{-n}<e^{-n/d}$, as $1/d\ll\beta$. Therefore, with positive probability, both bullet points hold at the same time. Fix such a choice of $G_1$ and $G_2$.

We now find the Hamilton cycles one at a time. Suppose that we have found a collection $\mathcal{C}$ of $i\ge 0$ edge-disjoint Hamilton cycles in $G$, and let $H$ be the spanning subgraph of $G$ with edge set $E(G_2)\cap E(\mathcal{C})$. We keep the $5\beta n$-max degree $\Delta_{5\beta n}(H)$ (Definition~\ref{def:degree}) small, so that the unused edges of $G_2$ still join any two disjoint sets of $5\beta n$ vertices.
\begin{claim}\label{cl: iteration}Suppose that $i\le (1-\eps) d/2$ and $\Delta_{5\beta n}(H) < \beta^4 d n$. Then, there is a Hamilton cycle $C$ in $G$ that is edge-disjoint from  $\mathcal{C}$, such that the spanning subgraph $H'$ of $G$ with edge set $E(G_2)\cap(E(\mathcal{C})\cup E(C))$ satisfies
\[\Delta_{5\beta n}(H') \leq  \Delta_{5\beta n}(H) + 5\cdot10^5 \gamma^*\log(1/\gamma^*) n ,\]
where $\gamma^*= \varepsilon^{-1}(24 \gamma + \Delta_{5\beta n}(H)/(\beta d n))$.
\end{claim}
\begin{claimproof}[Proof of Claim~\ref{cl: iteration}]
Let us denote $M = \Delta_{5\beta n}(H)$, $G_1'=G_1\setminus E(\mathcal C)$ and $G_2'=G_2\setminus E(\mathcal C)$. The following facts hold.
\begin{itemize}
    \item $\Delta(H) \le M/(5\beta n)\le \beta^3d/5$.
    \item $d_{G_1'}(v)=d_1\pm d_1'$ for all $v\in V(G_1')$, where $d_1=(1-\alpha)d-2i+\frac{M}{10\beta n}$ and $d_1'=2\gamma d+\frac{M}{10\beta n}$.
    \item $G_2'$ has an edge between any two disjoint sets of at least $5\beta n$ vertices.
\end{itemize}
The first bullet point follows directly from Definition \ref{def:degree}. For the second, note that $E(\mathcal{C})$ is $2i$-regular and, since $G_1$ and $G_2$ partition $G$, for any $v\in V(G_1')$ we have 
\[d_{G_1'}(v)=d_{G_1}(v)-2i+d_H(v)=(1-\alpha)d-2i+d_H(v)\pm 2\gamma d=d_1\pm d_1',\]
where we used that $d_H(v)\le \Delta(H) \leq M/(5\beta n)$ by the first bullet point. For the third bullet point, note that if $X,Y \subset V(G_2')$ are disjoint sets of size $5\beta n$, then 
\[e_{G_2'}(X,Y) \geq e_{G_2}(X,Y) - e_{H}(X,Y) \geq \beta^3 dn - |X| \Delta(H)>0,\]
as $|X|\Delta(H)\le 2\beta^4 dn$ by the first bullet point for sufficiently large $n$.

 We now use Lemma~\ref{lem:regularHamcycle} to find $C$. If $n$ is odd, choose $z\in V(G)$ and distinct $x_1,x_2\in N_{G_1'}(z)$, and let $R=V(G)\setminus\{z\}$. If $n$ is even, choose an edge $x_1x_2\in E(G_1')$ and let $R=V(G)$. These choices exist, as $2i\leq(1-\varepsilon)d$, $\alpha=\eps/10$ and $\gamma\ll\varepsilon$ give $d_1-d_1'\geq\varepsilon d/2$. In both cases, we set $m=|R|/2=\lfloor n/2\rfloor$ and partition $R$ into pairs so that $\{x_1,x_2\}$ is one of the pairs.

Lemma~\ref{lem:balanced-random-splitting}, applied to $G_1'[R]$ with $k=2$ and $W=R$, gives a partition $R=V_1\cup V_2$ with $|V_1|=|V_2|=m$ that splits each pair and satisfies $d_{G_1'}(v,V_\ell)=d_{G_1'[R]}(v)/2\pm d^{2/3}$ for all $v\in R$ and $\ell\in[2]$. We may assume that $x_1\in V_1$ and $x_2\in V_2$. As $d_{G_1'[R]}(v)=d_{G_1'}(v)\pm1$ and $d^{2/3}+1\leq d_1'$ by $1/d_0\ll\gamma$, all degrees in $G_1'[V_1,V_2]$ are $d_1/2\pm2d_1'$, and
\[
 \frac{\Delta(G_1'[V_1,V_2])}{\delta(G_1'[V_1,V_2])}
 \leq1+\frac{8d_1'}{d_1-4d_1'}
 \leq1+\frac{10}{\varepsilon}\left(2\gamma+\frac{M}{10\beta nd}\right)
 \leq1+\gamma^*.
\]
Here $d_1-4d_1'\geq4\varepsilon d/5$ follows from $2i\leq(1-\varepsilon)d$ and $M<\beta^4dn$, as $\gamma\ll\beta\ll\varepsilon$.

Set $\beta'=5\beta n/m$ and $d'=d/2$. By Lemma~\ref{lem:restriction-normalization}, $G[R]$ is $(\eta,\beta,d'/m)$-sparse, and thus its subgraph $G_1'[V_1,V_2]$ is $(\eta,\beta',d'/m)$-sparse as $\beta\leq\beta'$. Its minimum degree is at least $(\varepsilon/4)d'$, and $G_2'[V_1,V_2]$ is $5\beta n$-bipartite-joined by the third bullet point. Moreover, $\gamma\leq\gamma^*\leq(24\gamma+\beta^3)/\varepsilon$ as $M<\beta^4dn$, and $\beta\leq\beta'\leq11\beta$ as $n$ is large. Thus $\varepsilon/4$, $\beta'$ and $\gamma^*$ satisfy the hierarchy of Lemma~\ref{lem:regularHamcycle} in the roles of $\alpha$, $\beta$ and $\gamma$ for all $i$. We apply Lemma~\ref{lem:regularHamcycle} to $G_1'[V_1,V_2]$, $G_2'[V_1,V_2]$ and $H$, with $\varepsilon/4$, $m$, $d'$, $\beta'$ and $\gamma^*$ playing the roles of $\alpha$, $n$, $d$, $\beta$ and $\gamma$, and obtain a Hamilton $x_1,x_2$-path $P$ on $R$ such that $E(P)\setminus E(G_1')$ is
\[
 (10^6\gamma^*\log(1/\gamma^*)m,5\beta n,H)\text{-bounded}.
\]
Let $C=P+x_1z+zx_2$ if $n$ is odd, and $C=P+x_1x_2$ if $n$ is even. Since the added edges belong to $G_1'$, the graph $C$ is a Hamilton cycle in $G$ which is edge-disjoint from $\mathcal C$. Moreover, as $E(P)\cap E(G_2')\subset E(P)\setminus E(G_1')$ and $m\leq n/2$, we have
\[
 H'=H+(E(P)\cap E(G_2')),
 \qquad
 \Delta_{5\beta n}(H')\leq\Delta_{5\beta n}(H)+5\cdot10^5\gamma^*\log(1/\gamma^*)n.
\]
\end{claimproof}
We now apply Claim~\ref{cl: iteration} repeatedly. Suppose that we have found edge-disjoint Hamilton cycles $C_1,\ldots,C_i$. Let $H_i$ be the spanning subgraph of $G$ with edge set $E(G_2)\cap E(C_1\cup\dots\cup C_i)$, and set $\gamma_i=24\gamma+\Delta_{5\beta n}(H_i)/(\beta dn)$, which is $\eps\gamma^*$ in Claim~\ref{cl: iteration} for $H=H_i$. If $i<(1-\eps)d/2$ and $\gamma_i<\beta^3$, then $\Delta_{5\beta n}(H_i)<\beta^4dn$, and Claim~\ref{cl: iteration} gives a Hamilton cycle $C_{i+1}$ with
\[\gamma_{i+1}\leq\gamma_i+\frac{5\cdot10^5}{\beta d}\cdot\frac{\gamma_i}{\eps}\log\Bigl(\frac{\eps}{\gamma_i}\Bigr)\leq\gamma_i\Bigl(1+a\log\frac1{\gamma_i}\Bigr)\leq\gamma_i^{1-a},\]
where $a=5\cdot10^5/(\eps\beta d)$ and we used $1+x\leq e^x$ in the last step. By the hierarchy, we may assume that $a\le1/2$, $d\geq6/\gamma$ and
\[(25\gamma)^{\exp(-10^6/(\eps\beta))}<\beta^3.\]
Note that $H_0$ is the empty graph on $V(G)$, so $\Delta_{5\beta n}(H_0)=\lceil5\beta n\rceil$ and $\gamma_0\leq24\gamma+6/d\leq25\gamma<\beta^3$. As $\log(1-a)\geq-2a$ and $i+1\leq d$, induction on $i$ shows that, for every integer $0\le i<(1-\eps)d/2$, the cycle $C_{i+1}$ exists and
\[\gamma_{i+1}\leq\gamma_0^{(1-a)^{i+1}}\leq(25\gamma)^{\exp(-10^6/(\eps\beta))}<\beta^3.\]
In particular, we find $(1-\eps)d/2$ edge-disjoint Hamilton cycles in $G$.
\end{proof}

\subsection{The covering result}
\begin{proof}[Proof of Theorem~\ref{thm: covering}]
It suffices to prove the statement for sufficiently small $\varepsilon>0$. Let $\eta,d_0,\beta,\gamma>0$ be given by Theorem~\ref{thm: decomp} with $\eps^5$ playing the role of $\eps$. If necessary, we decrease $\eta,\beta,\gamma$ and increase $d_0$ so that $1/d_0\ll\eta,\gamma\ll\beta\ll\eps\ll1$. Let $n,d\in\mathbb N$ with $d_0\le d<n$ and let $G$ be an $n$-vertex $(\eta,\beta,d/n)$-sparse graph with $d_G(v)=(1\pm\gamma)d$ for every $v\in V(G)$. 
By Theorem~\ref{thm: decomp}, there are $s=\lfloor(1-\eps^5)d/2\rfloor$ edge-disjoint Hamilton cycles $C_1,\dots,C_s$ in $G$. Let $G'$ be the graph consisting of all edges not covered by $\bigcup_{i\in [s]}C_i$ and note that  
    $$(\eps^5 - \gamma) d\leq \delta(G')\leq \Delta(G')\leq (\eps^5+ \gamma) d+2.$$
By Vizing's theorem, the edge set of $G'$ can be partitioned into $t\le \Delta(G')+1\le 2\eps^5d$ matchings, say $M_1,\dots, M_t$. Let $k=1/\eps^2$. For $i\in[t]$, let $J_i$ be the bipartite graph with parts $V(G)$ and $E(M_i)$ in which $v$ is adjacent to $e$ if $e$ meets $N_G(v)$. As $M_i$ is a matching, $\Delta(J_i)\leq2\Delta(G)\leq4d$. We partition all but at most $k-1$ edges of $M_i$ into sets of size $k$, apply Lemma~\ref{lem:balanced-random-splitting} to $J_i$ and these sets with $2d$ playing the role of $d$ (possible as $1/d_0\ll\eps$), and add the remaining edges to distinct classes. The resulting matchings $M_{i,1},\dots,M_{i,k}$ partition $M_i$ and satisfy $d_{J_i}(v,E(M_{i,j}))\leq d_G(v)/k+(2d)^{2/3}+1$ for all $v\in V(G)$ and $j\in[k]$, and $e(M_{i,j})\leq n/(2k)+1\leq\eps^2n$ as $e(M_i)\leq n/2$. Thus $d_G(v,V(M_{i,j}))\leq2d_{J_i}(v,E(M_{i,j}))\leq2(1+\gamma)\eps^2d+2(2d)^{2/3}+2\leq4\eps^2d$ for all $v\in V(G)$, $i\in[t]$ and $j\in[k]$.

For fixed $i\in [t]$ and $j\in [k]$, we now find a Hamilton cycle $C_{i,j}$ that contains all edges of $M_{i,j}$. 

If $n$ is even, choose an edge $x_1x_2$ of $G-V(M_{i,j})$. If $n$ is odd, choose $z\in V(G)\setminus V(M_{i,j})$ and distinct $x_1,x_2\in N_G(z)\setminus V(M_{i,j})$. Such a choice is possible, as every vertex has at least $(1-\gamma-4\eps^2)d\geq2$ neighbours outside $V(M_{i,j})$. It suffices to find a Hamilton $x_1,x_2$-path of $G$, or of $G-z$ if $n$ is odd, that contains $M_{i,j}$, and then to add the edge $x_1x_2$, or the path $x_1zx_2$. In the odd case, by Lemma~\ref{lem:restriction-normalization}, $G-z$ is $(\eta,\beta,d/(n-1))$-sparse with degrees $(1\pm2\gamma)d$. For this pair $i,j$, replacing $G$, $n$ and $\gamma$ by $G-z$, $n-1$ and $2\gamma$ in the odd case, we may assume that $n$ is even.
\begin{claim}\label{cl:matching-split}There are a set $S\subset V(G)\setminus(V(M_{i,j})\cup\{x_1,x_2\})$ of size $2\lceil n^{0.9}\rceil$ and a partition $V(G)=V_1\cup V_2$ with $|V_1|=|V_2|=n/2$, $x_1\in V_1$ and $x_2\in V_2$ such that, for $U_1=V(M_{i,j})\cup S$ and $U_2=V(G)\setminus U_1$, the following hold.
\begin{enumerate}[label=\upshape(\roman{enumi})]
    \item\label{random_partition_matching:1} $d_G(v,V_\ell\cap U_2)\ge d/8$ and $d_G(v,V_\ell)\ge (1-2\gamma)d/2$ for all $v\in V(G)$ and $\ell\in [2]$.
    \item\label{random_partition_matching:2} $U_1$ is balanced, and every edge of $M_{i,j}$ has one vertex in $V_1$ and one in $V_2$.
\end{enumerate}
\end{claim}
\begin{claimproof}[Proof of Claim~\ref{cl:matching-split}]
The set $V(G)\setminus(V(M_{i,j})\cup\{x_1,x_2\})$ has at least $0.9n$ vertices. We partition all but at most $9$ of them into sets of size $10$ and, as $\Delta(G)\leq2d$, apply Lemma~\ref{lem:balanced-random-splitting} to $G$ and these sets. Each of the ten resulting sets has at least $0.9n/10-1\geq2\lceil n^{0.9}\rceil$ vertices, and each $v\in V(G)$ has at most $(1+\gamma)d/10+d^{2/3}\leq d/8$ neighbours in it. Let $S$ consist of $2\lceil n^{0.9}\rceil$ vertices of one of them.

We partition $V(G)$ into pairs so that $\{x_1,x_2\}$ and the edges of $M_{i,j}$ are pairs and each remaining pair lies in $S$ or in $V(G)\setminus(U_1\cup\{x_1,x_2\})$, which is possible as $|S|$ and $n$ are even. We apply Lemma~\ref{lem:balanced-random-splitting} to $G$ and these pairs to obtain a partition $V(G)=V_1\cup V_2$ that splits each pair and satisfies $d_G(v,V_\ell)=d_G(v)/2\pm d^{2/3}=(1\pm2\gamma)d/2$ for all $v\in V(G)$ and $\ell\in[2]$. As each pair is split, \ref{random_partition_matching:2} holds, and we may assume that $x_1\in V_1$. Since every $v\in V(G)$ has at most $4\eps^2d$ neighbours in $V(M_{i,j})$ and at most $d/8$ in $S$, we have $d_G(v,V_\ell\cap U_2)\geq(1-2\gamma)d/2-4\eps^2d-d/8\geq d/8$, which proves~\ref{random_partition_matching:1}.
\end{claimproof}

    Let $D=\lceil1000C_{\ref{thm:bipartite-c-expander}}\rceil$. As $\gamma\ll\beta\ll1$, for sets $W_1\subset V_1$ and $W_2\subset V_2$ with $|W_1|=|W_2|=\beta^{1/2} n/D$, we have 
\begin{align*}
 e_G(W_1,W_2)&\geq(1-2\gamma)d|W_1|/2-e_G(W_1,V_2\setminus W_2)\\
 &\geq(1-2\gamma)d|W_1|/2
       -(1+\beta)(1-2\beta^{1/2}/D)d|W_1|/2
 \geq\beta dn/(2D^2).
\end{align*}
Let $G_1=G_2=G[V_1,V_2]$. Then $M_{i,j}\subset G_1$ by~\ref{random_partition_matching:2}, and $G_2$ is $\beta^{1/2}n/D$-bipartite-joined by the display above. To apply Lemma~\ref{lem: partition}, set $b=\beta^{1/2}n$, $r=1000b/D$ and $\rho=1/1000$. By~\ref{random_partition_matching:1}, each vertex has at least $d/8$ neighbours in $U_2$ in $G_1$. We have $\eta n\leq r$ and $2r\leq n/(32D+16)$, and thus Lemma~\ref{lem:sparse->expander}, with $1/8$ and $2D$ playing the roles of $\alpha$ and $D$, yields condition~\ref{eq: partition condition 1} of Lemma~\ref{lem: partition}.
Also, $M_{i,j}$ and the vertices of $S$, as paths of length $0$, form a spanning linear forest $\mathcal F$ of $G_1[U_1]$ with $e(M_{i,j})+|S|\leq10\varepsilon^2n$ paths. As $1/d_0\ll\beta\ll\varepsilon\ll1$, we have $10\varepsilon^2\log(e/(10\varepsilon^2))\leq\rho/10^4$ and $\beta^{1/2}\leq\rho/10^8$, and $n\geq d_0$ is large in terms of these constants. Thus we can apply Lemma~\ref{lem: partition}, with $10\varepsilon^2$, $\beta^{1/2}$ and the empty graph on $V(G)$ playing the roles of $\delta$, $\beta$ and $H$, to obtain a partition $U_2=X\cup Y\cup Z$, sets $A,B\subset X$ and a graph $F$ satisfying~\ref{eq: partition 1}--\ref{eq: partition 5}.
By~\ref{eq: partition 1}--\ref{eq: partition 3} and Lemma~\ref{lem:cycle-assembly}, applied as at the end of the proof of Lemma~\ref{lem:regularHamcycle}, $F$ has a spanning $x_1,x_2$-path $P$ containing all edges of $M_{i,j}$, and $V(P)=V(G)\setminus Z$ by~\ref{eq: partition 1}.
Note that one of the first two edges of $P$ is an edge $uv\notin E(M_{i,j})$ with $u\in V_1$ and $v\in V_2$, since $M_{i,j}$ is a matching.

It remains to cover $Z$. Let $C'=D/1000\ge C_{\ref{thm:bipartite-c-expander}}$. The set $Z$ is balanced, and by~\ref{eq: partition 1},
\[|Z|\ge(1-20\eps^2)n-2000\eps^2\log(e/(10\eps^2))n \ge(1-4\eps)n,\] as $\eps\ll1$. As $\beta^{1/2}n/D\leq|Z\cap V_1|/(2C')$, the graph $G_2[Z]$ is $|Z\cap V_1|/(2C')$-bipartite-joined.
We claim that $(G_1\cup G_2)[Z\cap V_1,Z\cap V_2]$ is a bipartite $C'$-expander. Indeed, by~\ref{eq: partition 4}, the only thing we have to check is that for each $\ell\in[2]$, subsets of $Z\cap V_\ell$ of size between $2r$ and $|Z\cap V_\ell|/(2C')$ expand by a factor of $C'$. So let $U\subset Z\cap V_\ell$ have size $2r\le |U|\le |Z\cap V_\ell|/(2C')$. As $G_2$ is $\beta^{1/2}n/D$-bipartite-joined and $|Z\cap V_{1}|=|Z\cap V_2|$, we have
\[|N_{G_2}(U)\cap (Z\cap V_{3-\ell}) |\geq |Z\cap V_{3-\ell}|-\beta^{1/2} n/D \geq C'|U|,\]
as $\beta^{1/2} n/D +C'|U|\le (C'+1)|U|\le (C'+1)|Z\cap V_{\ell}|/(2C')\le |Z\cap V_{3-\ell}|$. 
By~\ref{eq: partition 4}, we may choose
$z_1\in N_{G_1}(v)\cap Z\cap V_1$ and
$z_2\in N_{G_1}(u)\cap Z\cap V_2$.
Theorem~\ref{thm:bipartite-c-expander} gives a Hamilton path $Q$ in $(G_1\cup G_2)[Z]$ with endpoints $z_1$ and $z_2$. We replace the edge $uv$ of $P$ by the path $uz_2Qz_1v$, and this gives a Hamilton $x_1,x_2$-path of $G$ containing $M_{i,j}$, as desired.

Finally, the cycles $C_1,\ldots,C_s$ and $C_{i,j}$ for $i\in[t]$ and $j\in[k]$ cover all edges of $G$, and their number is at most
\[
 s+tk\leq(1-\varepsilon^5)d/2+2\varepsilon^3d
 \leq(1+\varepsilon)d/2.
\]
\end{proof}
%-------------------------------------------------------

\section*{Use of AI tools}
The mathematical content of this paper was developed without the use of AI tools. During the revision, AI tools (ChatGPT and Claude) were
used only for proofreading and correcting minor issues, with the exception of the appendix. When writing the appendix, we used AI tools to expand our detailed proof sketches from the original version—adapting results from~\cite{C-expander} to our setting—into complete proofs and to verify that the adapted arguments were correct.

\bibliographystyle{plain}

\clearpage
\pagenumbering{Roman}
\appendix
\section*{Appendix}
The linking and rotation arguments below follow~\cite{C-expander}. We give them in full, since we need bipartite versions of them. In these versions the endpoints are prescribed and the edges taken from $G_2$ form a bounded set as in Lemma~\ref{lem:regularHamcycle}.

\begin{itemize}
\item In Appendix~\ref{app:extendability} we present the extendability method, and in Appendix~\ref{app:template} we build a bipartite linking template. We embed the template and prove Lemma~\ref{lem: partition} in Appendix~\ref{app:linking}.
\item Appendix~\ref{app:rotations} contains the rotation lemmas used in Appendix~\ref{app:forests}. In Appendix~\ref{app:forests}, we use rotations and a few edges of $G_2$ to reduce a linear forest to at most $n^{0.8}$ paths (Lemma~\ref{lemma:linearforest decomp1}) and then to prescribe its endpoints (Lemma~\ref{lemma:linearforest decomp2}).
\item Appendix~\ref{app:bipartite-connectivity} proves Theorem~\ref{thm:bipartite-c-expander} by applying these results with $G_1=G_2=G$.
\end{itemize}

\section{The \texorpdfstring{$m$}{m}-max degree and bipartite extendability}\label{app:extendability}

In this section, we prove the basic properties of the $m$-max degree and present the extendability method as in~\cite{araujo2025ramsey, C-expander}. In Appendix~\ref{app:linking}, we build the graph $F$ of Lemma~\ref{lem: partition} by adding leaves (Lemma~\ref{lem:add-leaf}) and connecting paths (Lemma~\ref{lem:connect}). At the end we use Lemma~\ref{lem:unused-expansion} to prove property~\ref{eq: partition 4} of Lemma~\ref{lem: partition}. In Lemma~5.3 of~\cite{C-expander} the connecting paths lie inside one $C$-expander, but our $G_1$ is not necessarily $\beta n$-bipartite-joined, and so each connecting path of Lemma~\ref{lem:connect} uses one edge of $G_2$.

\subsection{The $m$-max degree}

\begin{lemma}\label{lem:max-degree}
Let $H$ be a graph, let $b\geq1$ be an integer, and let $\Delta_b(H)$ and $t(H)$ be as in Definition~\ref{def:degree} with $m=b$.
\begin{enumerate}[label=\upshape(\roman{enumi})]
\item If $H$ is a spanning subgraph of $H'$, then $\Delta_b(H)\leq\Delta_b(H')$.
\item $\Delta(H)\leq\Delta_b(H)/b\leq\Delta(H)+1$.
\item There is a set $Z\subset V(H)$ with $|Z|\leq b$ such that $\Delta_b(H+e)\leq\Delta_b(H)+2$ for every pair $e$ of distinct vertices in $V(H)\setminus Z$.
\item If $E_1$ is $(i_1,b,H)$-bounded and $E_2$ is $(i_2,b,H+E_1)$-bounded, then every $E\subset E_1\cup E_2$ is $(i_1+i_2,b,H)$-bounded.
\end{enumerate}
\end{lemma}
\begin{proof}
Recall from Definition~\ref{def:degree} that $t(H)$ is the number of vertices of degree $\Delta(H)$, that
\[
\Delta_b(H)=b\Delta(H)+\min\{t(H),b\},
\]
and that a set $E$ of pairs is $(i,b,H)$-bounded if $\Delta_b(H+E)\leq\Delta_b(H)+i$.

For (i), if $\Delta(H')>\Delta(H)$, then $\Delta_b(H')\geq b\Delta(H)+b+1>\Delta_b(H)$, since $1\leq\min\{t,b\}\leq b$ for both graphs. Otherwise $\Delta(H')=\Delta(H)$, and $t(H')\geq t(H)$ since every vertex of maximum degree in $H$ has maximum degree in $H'$.

Part (ii) holds as $0\leq\min\{t(H),b\}\leq b$.

For (iii), let $Z$ be the set of vertices of maximum degree if $t(H)\leq b$, and $Z=\emptyset$ otherwise. Let $e\notin E(H)$ avoid $Z$. If $t(H)\leq b$, the ends of $e$ have degree below $\Delta(H)$. Then $\Delta(H+e)=\Delta(H)$ and $t(H+e)\leq t(H)+2$, which implies that $\Delta_b$ grows by at most $2$. If $t(H)>b$ and $\Delta$ does not grow, then $t(H+e)\geq t(H)>b$ and $\Delta_b$ does not change. If $t(H)>b$ and $\Delta$ grows, it grows by one. The new maximum degree is then attained only at the ends of $e$. Thus $t(H+e)\leq2$, and $\Delta_b$ grows by at most $b+2-b=2$.

For (iv), we add the two inequalities and apply (i) to $H+E\subset H+E_1+E_2$.
\end{proof}

\subsection{Bipartite extendability}
\begin{definition}
Let $G$ be bipartite with parts $V_1,V_2$, and let $S$ be a graph with $V(S)\subset V(G)$ whose edges join $V_1$ to $V_2$ but do not have to belong to $G$. For integers $D\geq3$ and $r\geq1$, we say that $S$ is \emph{$(D,r)$-extendable with respect to $G$} if $\Delta(S)\leq D$ and
\begin{equation}\label{eq:extendability}
 |\Gamma_G(U)\setminus V(S)|\geq
 (D-1)|U|-\sum_{u\in U\cap V(S)}(d_S(u)-1)
\end{equation}
for every $U$ of size at most $2r$ inside $V_1$ or $V_2$.
\end{definition}

The right side of~\eqref{eq:extendability} is $\sum_{u\in U\cap V(S)}(D-d_S(u))+(D-1)|U\setminus V(S)|$. Thus~\eqref{eq:extendability} asks for one neighbour outside $V(S)$ for each edge that a vertex of $U\cap V(S)$ may still get, and for $D-1$ such neighbours for each vertex of $U\setminus V(S)$.

Lemma~A.2 of~\cite{araujo2025ramsey} requires $S\subset G$. However, in Appendix~\ref{app:linking} the graph $S$ contains edges of $G_2$, so we prove the version we need.

\begin{lemma}\label{lem:add-leaf}
Let $G$ be bipartite with parts $V_1,V_2$, and suppose $S$ is $(D,r)$-extendable with respect to $G$ and
\begin{equation}\label{eq:leaf-condition}
 |\Gamma_G(U)\setminus V(S)|\geq 2Dr+1
\end{equation}
whenever $U$ lies in one part and $r\leq|U|\leq2r$. For every $x\in V(S)$ with $d_S(x)\leq D-1$, there is a vertex $y\in\Gamma_G(x)\setminus V(S)$ such that $S+xy$ is $(D,r)$-extendable with respect to $G$.
\end{lemma}
\begin{proof}
Adding a leaf $xy$ can only break~\eqref{eq:extendability} for a set $U$ with equality in~\eqref{eq:extendability} and $y\in\Gamma_G(U)$. We take the union $T$ of all sets with equality that are contained in the part of $x$ and avoid $x$. We show that $|T|<r$, and we choose $y$ outside $\Gamma_G(T)$.

We may assume by symmetry that $x\in V_1$. For $U\subset V_1$, let $f(U)$ be the left side minus the right side of~\eqref{eq:extendability}, that is, 
\[f(U)=|\Gamma_G(U)\setminus V(S)|-
 (D-1)|U|+\sum_{u\in U\cap V(S)}(d_S(u)-1).\]
 The function $f$ is submodular, since $U\mapsto|\Gamma_G(U)\setminus V(S)|$ is submodular and the other two terms are sums over $U$. If $|U|\leq2r$, then $f(U)\geq0$ by extendability.

Since the right side of~\eqref{eq:extendability} is at most $D|U|$, \eqref{eq:leaf-condition} implies $f(U)>0$ if $r\leq|U|\leq2r$. Thus every $U\subset V_1$ with $|U|\leq2r$ and $f(U)=0$ has size less than $r$. For two such sets $U,U'$ we have $|U\cup U'|<2r$, and hence $f(U\cup U')$ and $f(U\cap U')$ are nonnegative. Submodularity implies that they sum to at most $f(U)+f(U')=0$, and so $f(U\cup U')=0$. Hence the union $T$ of all such sets in $V_1\setminus\{x\}$ satisfies $f(T)=0$ and $|T|<r$.

Let $y\in\Gamma_G(x)\setminus V(S)$. Then $\Delta(S+xy)\leq D$ as $d_S(x)\leq D-1$. For $U\subset V_2$, adding the leaf $xy$ to $S$ changes neither side of~\eqref{eq:extendability}, as the only new vertex $y$ is in $V_2$ and has degree one. For $U\subset V_1$, it changes $f(U)$ by
\[
\mathbf1_{\{x\in U\}}-\mathbf1_{\{y\in\Gamma_G(U)\}}.
\]
As this change is at least $-1$ and $f$ is integral, $S+xy$ is extendable unless $y\in\Gamma_G(U)$ for some $U\subset V_1\setminus\{x\}$ with $|U|\leq2r$ and $f(U)=0$. In particular, it is extendable if $y\notin\Gamma_G(T)$.

It remains to find such a $y$. If $\Gamma_G(x)\setminus V(S)\subset\Gamma_G(T)$, then $f(T\cup\{x\})=f(T)-(D-d_S(x))<0$. This contradicts extendability, since $|T\cup\{x\}|\leq r$.
\end{proof}

\begin{lemma}\label{lem:connect}
Let $D\geq3$ and $r,b,s\geq1$ be integers with $s\leq2b$, and set $k=\lceil\log_2(4b/s)\rceil$.
Let $G_1,G_2$ be bipartite graphs on the same parts such that $G_2$ is $b$-bipartite-joined, let $S$ be $(D,r)$-extendable with respect to $G_1$, and let $H$ be any graph with $V(H)\supset V(G_1)$.

Let $a_1,\ldots,a_s,a_1',\ldots,a_s'$ be distinct vertices of $S$ of degree at most $D-1$ in $S$, with all $a_i$ in one part and all $a_i'$ in one part, possibly the same one. Let $\ell\geq2k+3$ have the parity required for a path from the part of the $a_i$ to the part of the $a_i'$, and suppose that
\begin{equation}\label{eq:connect-condition}
 |\Gamma_{G_1}(U)\setminus V(S)|\geq
 2Dr+1+50b+s(\ell-2k-1)
\end{equation}
whenever $U$ lies in one part and $r\leq|U|\leq2r$.

Then, for some $i,j\in[s]$, there are an $a_i,a_j'$-path $P$ of length $\ell$ and an edge $e\in E(P)\cap E(G_2)$ such that the following hold.
\begin{enumerate}[label=\upshape(\roman{enumi})]
\item All internal vertices of $P$ lie outside $V(S)$, and every edge of $P$ other than $e$ belongs to $G_1$.
\item The set $\{e\}$ is $(2,b,H)$-bounded.
\item $S\cup P$ is $(D,r)$-extendable with respect to $G_1$.
\end{enumerate}
\end{lemma}
\begin{figure}[htbp]
  \centering
  \begin{tikzpicture}[x=1cm,y=1cm,
      v/.style={circle,fill=black,inner sep=1pt},
      lab/.style={font=\scriptsize,inner sep=1.5pt},
      wavy/.style={decorate,decoration={snake,amplitude=.3mm,segment length=1.6mm}},
      pa/.style={DarkDesaturatedBlue,line width=1pt},
      off/.style={black!35,densely dashed,line width=.4pt}]
    \node[v] (a) at (0,0) {}; \node[lab,left] at (a) {$a_i$};
    \node[v] (r) at (3,0) {};
    \node[v] (b) at (9,0) {}; \node[lab,right] at (b) {$a_j'$};
    \node[v] (s) at (8,0) {};
    \foreach \y/\n in {.6/1,-.6/2}{\node[v] (r\n) at (4,\y) {}; \node[v] (s\n) at (7,\y) {};}
    \foreach \y/\n in {.9/1,.3/2,-.3/3,-.9/4}{\node[v] (rl\n) at (5,\y) {}; \node[v] (sl\n) at (6,\y) {};}
    \draw[pa,wavy] (a)--(r); \draw[pa] (r)--(r1)--(rl2);
    \draw[off] (r)--(r2) (r1)--(rl1) (r2)--(rl3) (r2)--(rl4);
    \draw[pa] (b)--(s)--(s1)--(sl2);
    \draw[off] (s)--(s2) (s1)--(sl1) (s2)--(sl3) (s2)--(sl4);
    \draw[VividOrange,line width=1.2pt] (rl2)--node[lab,above]{$e$}(sl2);
  \end{tikzpicture}
  \caption{The proof of Lemma~\ref{lem:connect} for $k=2$. The trees at $a_i$ and $a_j'$ are joined by an edge $e$ of $G_2$, and the dashed parts are deleted. The solid path, with $e$, is $P$.}
  \label{fig:connect}
\end{figure}
\begin{proof}
We grow a binary tree at the end of a path from each $a_i$ and each $a_i'$ by Lemma~\ref{lem:add-leaf}. Then we join two leaves by an edge of $G_2$ and delete the vertices we do not use (see Figure~\ref{fig:connect}). More precisely, we grow a path of length $\ell-2k-2$ from each $a_i$ and a path of length one from each $a_i'$. At the end of each of these $2s$ paths, we attach a complete binary tree of depth $k$ rooted at that end. As $4b/s\leq2^k<8b/s$, the paths and trees add
\[
 s(\ell-2k-1)+2s(2^{k+1}-2)<s(\ell-2k-1)+32b
\]
vertices. Hence the graph built so far always has less than $s(\ell-2k-1)+32b$ vertices outside $V(S)$, and~\eqref{eq:connect-condition} implies condition~\eqref{eq:leaf-condition} of Lemma~\ref{lem:add-leaf} with $G_1$ and this graph playing the roles of $G$ and $S$. We can therefore build these paths and trees by applying Lemma~\ref{lem:add-leaf} iteratively, as its degree condition holds too: each $a_i$ and $a_i'$ gets one new edge, and new vertices have one parent and at most two children.

A leaf of the tree at $a_i$ has distance $\ell-k-2$ from $a_i$, and a leaf of the tree at $a_j'$ has distance $k+1$ from $a_j'$. These distances sum to $\ell-1$. The parity of $\ell$ implies that the leaves of the trees at the $a_i$ are in one part, and the leaves of the trees at the $a_j'$ in the other part.

Each leaf set has $s2^k\geq4b$ vertices, and the set $Z$ of Lemma~\ref{lem:max-degree}(iii) for $H$ has at most $b$ vertices. Since $G_2$ is $b$-bipartite-joined, we can choose an edge $e\in E(G_2)$ that avoids $Z$ and joins the two leaf sets, say a leaf of the tree at $a_i$ and a leaf of the tree at $a_j'$. Then $\{e\}$ is $(2,b,H)$-bounded, because $e$ avoids $Z$. Moreover, $e$ and the paths from $a_i$ and $a_j'$ to its ends form an $a_i,a_j'$-path $P$ of length $\ell$.

Deleting the other new vertices leaf by leaf preserves extendability, because their edges belong to $G_1$. Indeed, let $y$ be a leaf with neighbour $x\in V_1$, say. Deleting $y$ raises the left side of~\eqref{eq:extendability} for $U\subset V_1$ by $\mathbf1_{\{y\in\Gamma_{G_1}(U)\}}\geq\mathbf1_{\{x\in U\}}$ and the right side by $\mathbf1_{\{x\in U\}}$. For $U\subset V_2$, deleting $y$ changes neither side, as $y$ contributed zero to the sum.

After these deletions, the ends of $e$ have degree one. Adding $e$ does not change the left side and does not increase the right side, and the maximum degree is still at most $D$. This gives (iii). Part (i) holds as the internal vertices of $P$ are new and every edge of $P$ other than $e$ was added by Lemma~\ref{lem:add-leaf} with $G_1$.
\end{proof}

\begin{lemma}\label{lem:unused-expansion}
If $S$ is $(D,r)$-extendable with respect to a bipartite graph $G$ and $\Delta(S)\leq3$, then every set $U\subset V(G)$ of size at most $2r$ satisfies
\[
 |N_G(U)\setminus V(S)|\geq(D-4)|U|.
\]
\end{lemma}
\begin{proof}
Each summand $d_S(u)-1$ in~\eqref{eq:extendability} is at most $2$, and thus $|\Gamma_G(U')\setminus V(S)|\geq(D-3)|U'|$ for $U'=U\cap V_1$ and $U'=U\cap V_2$. These two neighbourhoods are subsets of $V_2$ and $V_1$ and are therefore disjoint. If we remove $U$ from their union, we lose at most $|U|$ vertices, and what remains is $N_G(U)\setminus V(S)$.
\end{proof}

\section{A bipartite linking template}\label{app:template}

We adapt the construction of Lemma~5.5 in~\cite{C-expander}, which replaces the copies of $K_{2,2}$ in a sorting network by the gadgets of~\cite{hyde2023spanning}. The linking structure $H^*$ of Lemma~\ref{lem:template}, which we call the \emph{template}, must be bipartite, with $A^*$ and $B^*$ in opposite parts. For this we choose the lengths of the paths of the gadget so that the gadget is bipartite (Lemma~\ref{lem:gadget}). We build the template from paths, since Appendix~\ref{app:linking} embeds it one path at a time by Lemmas~\ref{lem:add-leaf} and~\ref{lem:connect}.

\begin{definition}
Let $H$ be a graph and $R\subset V(H)$. We say that $H$ is \emph{$R$-path-constructible} if it is obtained from the edgeless graph on $R$ by adding edge-disjoint paths, each with an endpoint already present and all its internal vertices new.
\end{definition}

\begin{lemma}\label{lem:gadget}
For every sufficiently large integer $h\equiv2\pmod4$, there is a bipartite $(\{a_1,a_2\},\{b_1,b_2\})$-linking structure (Definition~\ref{def: linking structure}) with $a_1,a_2$ in one part and $b_1,b_2$ in the other. This linking structure has maximum degree three and $O(h^2)$ vertices, and $a_1,a_2,b_1,b_2$ have degree one. Moreover, it is $\{a_1,a_2\}$-path-constructible and $\{a_1,a_2,b_1,b_2\}$-path-constructible by paths of length $h$ or $h+1$.
\end{lemma}
\begin{figure}[htbp]
  \centering
  \begin{tikzpicture}[x=.82cm,y=1.35cm,
      v/.style={circle,fill=black,inner sep=1pt},
      lab/.style={font=\scriptsize,inner sep=1.5pt},
      wavy/.style={decorate,decoration={snake,amplitude=.3mm,segment length=1.6mm}},
      pa/.style={DarkDesaturatedBlue,line width=1pt},
      pb/.style={VividOrange,line width=1pt},
      off/.style={black!35,densely dashed,line width=.4pt}]
    \foreach \panel/\dx in {0/0,1/9.6}{
      \begin{scope}[xshift=\dx cm]
        \foreach \i in {1,...,6}{
          \node[v] (t\i) at (\i,1) {};
          \node[v] (u\i) at (\i,0) {};
        }
        \node[v] (a1) at (-.4,1) {}; \node[lab,left] at (a1) {$a_1$};
        \node[v] (a2) at (-.4,0) {}; \node[lab,left] at (a2) {$a_2$};
        \node[v] (b1) at (7.4,1) {}; \node[lab,right] at (b1) {$b_1$};
        \node[v] (b2) at (7.4,0) {}; \node[lab,right] at (b2) {$b_2$};
        \foreach \i/\p in {1/above,2/above left,3/above right,4/above left,5/above right,6/above}{
          \node[lab,\p] at (t\i) {$t_{\i}$};}
        \foreach \i/\p in {1/below,2/below left,3/below right,4/below left,5/below right,6/below}{
          \node[lab,\p] at (u\i) {$u_{\i}$};}
        \ifnum\panel=0
          \foreach \x/\y in {t1/u2,t6/u5,t2/u4,t4/u6,t3/u1,t5/u3} \draw[off] (\x)--(\y);
          \draw[pa,wavy] (a1)--(t1); \draw[pa] (t1)--(t2);
          \draw[pa,wavy] (t2) to[bend left=80,looseness=1.5] (t3);
          \draw[pa] (t3)--(t4);
          \draw[pa,wavy] (t4) to[bend left=80,looseness=1.5] (t5);
          \draw[pa] (t5)--(t6); \draw[pa,wavy] (t6)--(b1);
          \draw[pb,wavy] (a2)--(u1); \draw[pb] (u1)--(u2);
          \draw[pb,wavy] (u2) to[bend right=80,looseness=1.5] (u3);
          \draw[pb] (u3)--(u4);
          \draw[pb,wavy] (u4) to[bend right=80,looseness=1.5] (u5);
          \draw[pb] (u5)--(u6); \draw[pb,wavy] (u6)--(b2);
          \node[font=\small] at (3.5,-.85) {$Q\cup M_0$};
        \else
          \foreach \x/\y in {t1/t2,t3/t4,t5/t6,u1/u2,u3/u4,u5/u6} \draw[off] (\x)--(\y);
          \draw[pa,wavy] (a1)--(t1); \draw[pa] (t1)--(u2);
          \draw[pa,wavy] (u2) to[bend right=80,looseness=1.5] (u3);
          \draw[pa] (u3)--(t5);
          \draw[pa,wavy] (t4) to[bend left=80,looseness=1.5] (t5);
          \draw[pa] (t4)--(u6); \draw[pa,wavy] (u6)--(b2);
          \draw[pb,wavy] (a2)--(u1); \draw[pb] (u1)--(t3);
          \draw[pb,wavy] (t2) to[bend left=80,looseness=1.5] (t3);
          \draw[pb] (t2)--(u4);
          \draw[pb,wavy] (u4) to[bend right=80,looseness=1.5] (u5);
          \draw[pb] (u5)--(t6); \draw[pb,wavy] (t6)--(b1);
          \node[font=\small] at (3.5,-.85) {$Q\cup M_1$};
        \fi
      \end{scope}
    }
  \end{tikzpicture}
  \caption{The gadget of Lemma~\ref{lem:gadget} for $h=6$: $Q\cup M_0$ (left) and $Q\cup M_1$ (right). Wavy lines are paths of $Q$, and dashed edges belong to the other matching.}
  \label{fig:gadget}
\end{figure}
\begin{proof}
Write $h=2q$, where $q$ is odd. On a top row $t_1,\ldots,t_h$ and a bottom row $u_1,\ldots,u_h$ of vertices, take the matchings
\[
 M_0=\{t_{2i-1}t_{2i},u_{2i-1}u_{2i}:1\leq i\leq q\},\quad
 M_1=\{t_1u_2,t_hu_{h-1}\}\cup\{t_{2i}u_{2i+2},t_{2i+1}u_{2i-1}:1\leq i<q\}.
\]
Let $Q$ be a set of paths whose internal vertices are new and pairwise disjoint, as follows. For $1\leq i<q$, it has a path of length $h$ from $t_{2i}$ to $t_{2i+1}$ and one from $u_{2i}$ to $u_{2i+1}$. Moreover, $Q$ has four paths of length $h+1$, from $a_1$ to $t_1$, from $a_2$ to $u_1$, from $t_h$ to $b_1$ and from $u_h$ to $b_2$. The graph is $Q\cup M_0\cup M_1$, and it has $O(h^2)$ vertices. Figure~\ref{fig:gadget} shows it for $h=6$.

For the linking property, we list only the vertices $a_1,a_2,b_1,b_2$ and $t_i,u_i$ of each path. Then $Q\cup M_0$ consists of the two paths $a_1,t_1,\ldots,t_h,b_1$ and $a_2,u_1,\ldots,u_h,b_2$. The graph $Q\cup M_1$ consists of the paths
\[
\begin{split}
 &a_1,t_1,u_2,u_3,t_5,t_4,u_6,u_7,t_9,t_8,\ldots,u_h,b_2,\\
 &a_2,u_1,t_3,t_2,u_4,u_5,t_7,t_6,u_8,u_9,\ldots,t_h,b_1.
\end{split}
\]
The pattern has period four and ends as shown since $h\equiv2\pmod4$. Modulo four, the first path uses the top indices $0,1$ and the bottom indices $2,3$, and the second path uses the others. The two paths are therefore disjoint and cover all $t_i,u_i$, and hence all of $Q$. Thus the graph is an $(\{a_1,a_2\},\{b_1,b_2\})$-linking structure.

For bipartiteness, if we give $t_{2i-1},u_{2i-1}$ colour $i\bmod2$ and $t_{2i},u_{2i}$ the other colour, then edges of $M_0\cup M_1$ join different colours and the even paths of $Q$ join equal colours. The vertices $t_1,u_1$ have colour $1$ and $t_h,u_h$ have colour $0$, as $q$ is odd. Since the paths of length $h+1$ at $a_1,a_2,b_1,b_2$ are odd, $a_1,a_2$ are in one part and $b_1,b_2$ in the other. Moreover, the vertices $a_1,a_2,b_1,b_2$ have degree one, and all other vertices have degree at most three.

For the path construction, we want to add $M_0\cup M_1$ as two paths of length $h$, so we need $M_0\cup M_1$ to be a $2h$-cycle. For $h=6$, it consists of the straight edges, solid and dashed, in either panel of Figure~\ref{fig:gadget}. In general this is easy to check after contracting the edges $t_{2i-1}t_{2i}$ and $u_{2i-1}u_{2i}$ to vertices $t_i'$ and $u_i'$. Then $M_1$ becomes the set of edges $t_1'u_1'$, $t_q'u_q'$, $t_i'u_{i+1}'$ and $t_{i+1}'u_i'$ for $1\leq i<q$. Since $q$ is odd, these form the cycle $t_1',u_2',t_3',\ldots,t_q',u_q',t_{q-1}',u_{q-2}',\ldots,t_2',u_1',t_1'$.

We add the path from $a_1$ to $t_1$, then the two halves of the cycle between $t_1$ and its opposite vertex, and then the other three paths of length $h+1$ and the rest of $Q$. At the start either $\{a_1,a_2\}$ or $\{a_1,a_2,b_1,b_2\}$ is present. In both cases each added path has an endpoint already present, and its internal vertices are new. These $h+4$ paths have lengths $h$ and $h+1$.
\end{proof}

\begin{lemma}\label{lem:template}
For every sufficiently large integer $N$, there is a balanced bipartite graph $H^*$ with at most $N^{1.1}$ vertices, and sets $A^*,B^*$ in opposite parts, such that $H^*$ is an $(A^*,B^*)$-linking structure and the following hold.
\begin{enumerate}[label=\upshape(\roman{enumi})]
\item $|A^*|=|B^*|=N$, and the set $A^*\cup B^*$ of terminals is independent, with every terminal of degree one.
\item $\Delta(H^*)\leq3$.
\item $H^*$ is $(A^*\cup B^*)$-path-constructible with paths of length between $10\log N$ and $40\log N$.
\end{enumerate}
\end{lemma}
\begin{proof}
Let $G^*$ be the graph of the sorting network of~\cite{AKS} with $m=O(\log N)$ rounds of comparisons. The vertex set of $G^*$ is partitioned into sets $L_0,\ldots,L_m$ of size $N$, called layers, where $L_j$ corresponds to the $N$ positions after $j$ rounds. A pair of positions compared in round $j$ gives a copy of $K_{2,2}$ between $L_{j-1}$ and $L_j$, and each position not compared in round $j$ gives an edge. Thus every edge of $G^*$ joins two consecutive layers, and $G^*$ is an $(N,m)$-sorting network in the sense of Definition~5.6 of~\cite{C-expander}. In particular, $G^*$ is an $(L_0,L_m)$-linking structure. We may assume that $m$ is odd, as otherwise we append a layer joined to the last one by a perfect matching.

Let $h$ be the least integer with $h\equiv2\pmod4$ and $h\geq10\log N$. Note that $h,h+1\in[10\log N,40\log N]$ for large $N$. Let $A^*=L_0$ and $B^*=L_m$, so $|A^*|=|B^*|=N$. We obtain $H^*$ from $G^*$ by replacing each $K_{2,2}$ between consecutive layers $L_{j-1},L_j$ by a new copy of the gadget of Lemma~\ref{lem:gadget}, with $a_1,a_2\in L_{j-1}$ and $b_1,b_2\in L_j$. We replace each remaining edge of $G^*$ by a new path of length $h+1$. It is easy to see that $H^*$ is an $(A^*,B^*)$-linking structure, since $G^*$ and each gadget are linking structures.

We colour the layers alternately. This colouring extends to each gadget and each new path, since the gadget has $a_1,a_2$ and $b_1,b_2$ in opposite parts and the new paths have odd length. As $m$ is odd, $A^*$ and $B^*$ are in opposite parts. Hence each of the $N$ paths given by the linking property is balanced. Since these paths partition $V(H^*)$, the graph $H^*$ is balanced.

For (i) and (ii), vertices of $A^*\cup B^*$ have degree one, vertices of other layers have degree at most two, and all other vertices have degree at most three. No gadget has an edge between two of $a_1,a_2,b_1,b_2$, and no new path is a single edge. This implies that $A^*\cup B^*$ is independent. The graph $H^*$ has $O(N\log^3N)\leq N^{1.1}$ vertices, since it consists of $O(N\log N)$ gadgets and new paths with $O(h^2)$ vertices each.

It remains to prove (iii). We start with $A^*\cup B^*$ and, for $j\in [m]$ in turn, add the new paths and the gadgets between $L_{j-1}$ and $L_j$. Since $L_0=A^*$ and the paths and gadgets added in step $j$ cover $L_j$, all of $L_{j-1}$ is present in step $j$. Thus each new path has an endpoint already present. Each gadget has $a_1,a_2$ present, and its vertices $b_1,b_2$ are new if $j<m$ and belong to $B^*$ if $j=m$. Hence, by Lemma~\ref{lem:gadget}, the gadget can be constructed from $\{a_1,a_2\}$ if $j<m$ and from $\{a_1,a_2,b_1,b_2\}$ if $j=m$.
\end{proof}

\section{Embedding the linking structure}\label{app:linking}

We prove Lemma~\ref{lem: partition} by the extendability method, as in the proof of Lemma~5.1 in~\cite{C-expander}. Here the paths and the linking structure are bipartite. Recall what Lemma~\ref{lem: partition} gives (see Figure~\ref{fig:linking-output}). The short paths of $\mathcal F$, which cover $U_1$, are joined into $t$ paths through a small set $Y\subset U_2$. A linking structure $F[X]$ on their ends and on $x_1,x_2$ spans a small set $X\subset U_2$. Small sets still expand into the rest $Z$ of $U_2$.

\begin{proof}[Proof of Lemma~\ref{lem: partition}]
Recall that $b=\beta n$, $r=1000b/D$ and $t=n^{0.9}$. We extend the paths of $\mathcal F$, join them into $t$ paths, and link their ends and $x_1,x_2$ by a copy of the template of Lemma~\ref{lem:template}. We start with an empty set $E$, and in each application of Lemma~\ref{lem:connect} we use $H+E$ in place of $H$ and add to $E$ the edge of $G_2$ given by that lemma.

Write $S$ for the graph constructed so far. Initially $S$ consists of $\mathcal F$ and the isolated vertices $x_1,x_2$. This graph is $(D,r)$-extendable with respect to $G_1$. Indeed, let $U$ be a nonempty set in one part with $|U|\leq2r$. As $V(S)\cap U_2=\{x_1,x_2\}$, condition~\ref{eq: partition condition 1} implies $|\Gamma_{G_1}(U)\setminus V(S)|\geq2D|U|-2\geq D|U|$. The right side of~\eqref{eq:extendability} is at most $D|U|$. In every later step we apply Lemma~\ref{lem:add-leaf} or~\ref{lem:connect}, and thus $S$ stays extendable.

Lemmas~\ref{lem:add-leaf} and~\ref{lem:connect} also need sets of size between $r$ and $2r$ to have many neighbours outside $V(S)$. We get this by keeping $S$ small. Every new vertex belongs to $U_2$, since $U_1\subset V(S)$. We bound the number of vertices that each stage below adds, and at the end we check that, before each application of these lemmas, $S$ contains less than $200\delta\log(e/\delta)n$ vertices of $U_2$, including $x_1$ and $x_2$. This is less than $\rho n/2$, as $\delta\log(e/\delta)\leq\rho/10^4$. Therefore, condition~\ref{eq: partition condition 1} implies that throughout
\begin{equation}\label{eq:invariant}
 |\Gamma_{G_1}(U)\setminus V(S)|>\rho n/2>2Dr+1+54b+40\log(t+1)
\end{equation}
for every $U$ in one part with $r\leq|U|\leq2r$. The second inequality holds as $Dr\leq1000b+D$, $\beta\leq\rho/10^8$ and $n$ is large. In particular, \eqref{eq:leaf-condition} always holds. In each application of Lemma~\ref{lem:connect} below, the right side of~\eqref{eq:connect-condition} is $2Dr+1+50b$ plus a last term. We show that this last term is at most $4b$ when we join paths and at most $40\log(t+1)$ when we embed the template.

\paragraph{Extending and joining the paths.}
The paths of $\mathcal F$ lie in $U_1$, but the ends that we link must be in $U_2$, in opposite parts. We therefore first add one leaf to each isolated vertex of $\mathcal F$, and then extend both ends of each path of $\mathcal F$ by Lemma~\ref{lem:add-leaf}. We can choose the part of the new endpoint by adding one or two edges at an end. We choose the two endpoints in opposite parts, using at most three new vertices for each path. We now have $\ell$ disjoint paths with endpoints in opposite parts of $U_2$, each of length at least $3$ with at most $n^{0.15}+4$ vertices.

We next reduce the number of paths from $\ell$ to $t$. We apply Lemma~\ref{lem:connect} to many short paths at once, since then $k$ is small and the connecting path is short. Suppose that $j\ge t+1$ paths remain. Since the paths are disjoint and $|V_1\cup V_2|\leq2n$, at least $j/2$ of them have at most $4n/j$ vertices. We choose two disjoint families of $s$ such paths and apply Lemma~\ref{lem:connect} with odd length $2k+3$, where
\[
 s=\min\{j/4,2b\}\qquad\text{and}\qquad k=\lceil\log_2(4b/s)\rceil,
\]
and with the endpoints in $V_1$ of the first family as the $a_i$ and those in $V_2$ of the second as the $a_i'$. Since the length is $2k+3$, the last term of~\eqref{eq:connect-condition} is $2s\leq4b$, and~\eqref{eq:connect-condition} holds by~\eqref{eq:invariant}.

The connecting path joins two distinct paths into one path with endpoints in opposite parts, and it has $2k+2$ internal vertices. We now bound $k$. If $s=2b$, then $k=1$. Otherwise $s=j/4$, so $k\leq1+\log_2(4b/s)\leq1+\log_2(16n/j)$. In both cases, as $j\leq\delta n<n/10$,
\[
 2k+3\leq5+3\log(16n/j)\leq10\log(en/j).
\]
Hence the connecting path has at most $10\log(en/j)$ internal vertices.

Summing over $j$, this stage adds at most
\[
 10\sum_{j=t+1}^{\ell}\log(en/j)
 \leq10\ell\log(e^2n/\ell)\leq10\delta n\log(e^2/\delta)
 \leq20\delta\log(e/\delta)n
\]
vertices, as $\log(\ell!)\geq\ell\log\ell-\ell$, $\ell\leq\delta n$ and $x\log(e^2/x)$ increases on $(0,1]$.

We now bound the orders of the resulting paths. A path that is never joined has at most $n^{0.15}+4$ vertices. A join combines two paths with at most $4n/j<4n/t\leq4n^{0.1}$ vertices each and a connecting path with $O(\log n)$ vertices. The new path therefore has at most $9n^{0.1}\leq3n^{0.15}$ vertices.

Let $\mathcal P$ be the resulting $S$ without $x_1$ and $x_2$. Then $\mathcal P$ consists of $t$ paths, each with endpoints in opposite parts of $U_2$ and length between $3$ and $3n^{0.15}$, and $\mathcal P$ contains every edge of $\mathcal F$. Let $A=End(\mathcal P)\cap V_1$ and $B=End(\mathcal P)\cap V_2$.

\paragraph{Embedding the template.}
We take the template of Lemma~\ref{lem:template} with $N=t+1$, identify $A^*$ with $A\cup\{x_1\}$ and $B^*$ with $B\cup\{x_2\}$, and embed the paths of Lemma~\ref{lem:template}(iii) in order. If only one endpoint of a path is present, we embed the path by adding leaves (Lemma~\ref{lem:add-leaf}). If both endpoints are present, we apply Lemma~\ref{lem:connect} to them with $s=1$ and the prescribed length.

We check that Lemma~\ref{lem:connect} applies. The prescribed length is at most $40\log(t+1)$ by Lemma~\ref{lem:template}(iii), so the last term of~\eqref{eq:connect-condition} is at most $40\log(t+1)$ and~\eqref{eq:connect-condition} holds by~\eqref{eq:invariant}. This length has the required parity, as the template is bipartite. Lemma~\ref{lem:connect} also needs this length to be at least $2k+3$, where $k=\lceil\log_2(4b)\rceil$ is the value of $k$ for $s=1$. This holds for large $n$, as the length is at least $10\log(t+1)$.

\paragraph{Checking~\eqref{eq:invariant}.}
At most $3\ell\leq3\delta n$ vertices are added to $S$ when we extend the paths and at most $20\delta\log(e/\delta)n$ when we join them, and the template has at most $(t+1)^{1.1}=o(n)$ vertices. As $\log(e/\delta)\geq1$ and $n$ is large, less than $200\delta\log(e/\delta)n$ vertices of $U_2$, including $x_1$ and $x_2$, are in $S$ at the end, and so before each application of Lemmas~\ref{lem:add-leaf} and~\ref{lem:connect}, as $S$ only grows.

We have $\Delta(S)\leq3$ throughout, as the template has maximum degree three and meets $\mathcal P$ only in its terminals, which have degree one. As $D\geq8$, this gives the degree conditions of Lemmas~\ref{lem:add-leaf} and~\ref{lem:connect}.

\paragraph{Properties~\ref{eq: partition 1}--\ref{eq: partition 5} of Lemma~\ref{lem: partition}.}
Let $\mathcal L$ be the embedded template, and let
\[
 F=\mathcal P\cup\mathcal L,\qquad X=V(\mathcal L),\qquad
 Y=V(\mathcal P)\setminus(U_1\cup A\cup B),\qquad
 Z=U_2\setminus(X\cup Y).
\]
The graphs $\mathcal L$ and $\mathcal P$ meet exactly in $A\cup B$. The terminals are independent in $\mathcal L$, and no edge of $\mathcal P$ joins two terminals since its paths have length at least $3$. Hence $F[X]=\mathcal L$ and $F[U_1\cup Y\cup A\cup B]=\mathcal P$.

The template and each path of $\mathcal P$ are balanced. Since $U_1$ and $A\cup B$ are balanced, $Y$ is balanced, and then $Z$ is balanced too. We have $|X\cup Y|\leq200\delta\log(e/\delta)n$, as $X\cup Y\subset U_2\cap V(S)$. This proves~\ref{eq: partition 1}--\ref{eq: partition 3}.

The graph $F=S$ is $(D,r)$-extendable with respect to $G_1$, $\Delta(F)\leq3$ and $V(F)=(V_1\cup V_2)\setminus Z$, so Lemma~\ref{lem:unused-expansion} gives
\[
 |N_{G_1}(U)\cap Z|\geq(D-4)|U|\geq(D/2)|U|
\]
for every $U\subset V_1\cup V_2$ of size at most $2r$, as $D\geq8$. This proves~\ref{eq: partition 4}.

Finally, we applied Lemma~\ref{lem:connect} $\ell-t$ times to join paths and $o(n)$ times for the template, since the template has $o(n)$ vertices and each of its paths has a new internal vertex. We therefore have $|E|\leq\ell-t+o(n)\leq2\delta n$. Each addition to $E$ increased $\Delta_b(H+E)$ by at most two, by Lemma~\ref{lem:connect}(ii). Thus by Lemma~\ref{lem:max-degree}(iv) applied $|E|-1$ times, every subset of $E$ is $(2|E|,b,H)$-bounded. As $E(F)\setminus E(G_1)\subset E$ and $2|E|\leq4\delta n\leq10\delta n$, this proves~\ref{eq: partition 5}.
\end{proof}

\section{Rotations of bipartite linear forests}\label{app:rotations}
We prove bipartite versions of Lemmas~3.3--3.6 of~\cite{C-expander}. A rotation adds an edge at an endpoint and deletes an edge of the forest, so the endpoint moves. When a set expands into its interior, rotations in it with edges of $G_1$ from one endpoint reach at least $20b$ new endpoints (Lemma~\ref{lem:endpoint-growth}). Then one edge of $G_2$ joins two endpoints (Lemma~\ref{lem:independent-rotations}), and at most one edge of $G_2$ moves an endpoint into a given set (Lemma~\ref{lem:moving}).

Throughout this appendix and Appendix~\ref{app:forests}, $G_1$ and $G_2$ are bipartite graphs with parts $V_1,V_2$, where $|V_1|=|V_2|=n$, and $N=2n$. Linear forests are spanning unless said otherwise, and bipartite with these parts, and their edges do not have to belong to $G_1\cup G_2$. An isolated vertex of a linear forest counts as a path of length $0$ and is an endpoint. We fix $0<\beta<1$ and set
\[
 D=2^{100},\qquad b=\beta n,\qquad r=1000b/D,\qquad K=\lceil2\log_DN\rceil+2,
\]
and we take $n$ large in terms of $\beta$, so that $K\log2\leq\log N$. We assume that $G_2$ is $b$-bipartite-joined. Boundedness is relative to a graph $H$ with $V(H)\supset V_1\cup V_2$. For a set $E$ of edges, we also write $G\cup E$ for $G+E$.

For a linear forest $\mathcal F$, spanning or not, and a set $U\subset V(\mathcal F)$, we call
\[
 int_{\mathcal{F}}(U)=\{ u\in U: N_{\mathcal{F}}(u)\subset U\}
\]
the \emph{interior} of $U$ in $\mathcal F$. As in Observation~2.4 of~\cite{C-expander}, if we delete a set $T$ from $U$, then at most $3|T|$ vertices are removed from $int_{\mathcal F}(U)$, and all of them are in $T$ or are neighbours in $\mathcal F$ of vertices of $T$.

\subsection{Definitions and endpoint growth}

\begin{definition}\label{def:rotation}
Let $\mathcal F$ be a linear forest. For an internal vertex $z$ of a path, we call the two paths obtained by deleting $z$ from it the \emph{branches} at $z$. Let $v\in End(\mathcal F)$ and let $vz$ be an edge of a given graph. A \emph{$1$-rotation} with \emph{old endpoint} $v$ and \emph{pivot} $z$ adds $vz$. If $z$ is isolated, it deletes no edge, and we let $y=z$. Otherwise it deletes the edge $zy$ of $\mathcal F$, where $y$ is the neighbour of $z$ given below.
\begin{enumerate}[label=\upshape(\alph{enumi})]
\item If $v$ and $z$ lie on the same path, $y$ is the neighbour of $z$ towards $v$.
\item If $z$ is an internal vertex of another path, $y$ is the neighbour of $z$ on the shorter branch at $z$, with ties broken arbitrarily (see Figure~\ref{fig:rotation}).
\item If $z$ is an endpoint of another path, $y$ is the neighbour of $z$ on that path.
\end{enumerate}
We call $y$ the \emph{new endpoint} and, when an edge is deleted, $zy$ the \emph{broken edge}.

As in Definition~3.2 of~\cite{C-expander}, a \emph{rotation} from $v$ is a sequence of consecutive $1$-rotations, called its \emph{steps}, such that the first starts at $v$ and the old endpoint of each step is the new endpoint of the previous one. The rotation must also satisfy the \emph{distance condition} that each pivot has distance at least three in the original forest $\mathcal F$ from $v$ and from every previous pivot, where distances between components are infinite. The new endpoint of the rotation is that of its last step, or $v$ if it has no steps.

We say that a rotation is \emph{in $U$} if all its pivots lie in $int_{\mathcal F}(U)$. Let $E_U^j(v,\mathcal F)$ be the set of new endpoints, including $v$, of rotations in $U$ from $v$ with at most $j$ steps and added edges from $G_1$.
\end{definition}

\begin{figure}[htbp]
  \centering
  \begin{tikzpicture}[x=.8cm,y=1cm,
      v/.style={circle,fill=black,inner sep=1pt},
      lab/.style={font=\scriptsize,inner sep=2pt}]
    \foreach \i in {0,...,4}{\node[v] (p\i) at (\i,1) {};}
    \draw (p0)--(p4);
    \node[lab,above] at (p4) {$v$};
    \foreach \i in {0,...,8}{\node[v] (q\i) at (\i+2,0) {};}
    \draw (q0)--(q2); \draw (q3)--(q8);
    \draw[black!35,densely dashed] (q2)--(q3);
    \node[lab,below] at (q2) {$y$};
    \node[lab,below] at (q3) {$z$};
    \draw[DarkDesaturatedBlue,line width=1pt] (p4)--(q3);
  \end{tikzpicture}
  \caption{A $1$-rotation in case~(b) of Definition~\ref{def:rotation}. The blue edge $vz$ is added and the dashed edge $zy$ is deleted, and $y$ is the new endpoint.}
  \label{fig:rotation}
\end{figure}

In case~(b) of Definition~\ref{def:rotation}, Definition~3.1 of~\cite{C-expander} allows either neighbour of $z$. The shorter branch is needed for Lemma~\ref{lem:log-orders}, which bounds how the orders of the paths change.

The next lemma collects basic facts about rotations. Part~(iii) bounds the number of pivots used to reach a set of new endpoints, so that a new pivot can be chosen at distance at least three from all these pivots.

\begin{lemma}\label{lem:rotation-basics}
Let $\mathcal F$ be a linear forest.
\begin{enumerate}[label=\upshape(\roman{enumi})]
\item In a rotation, every broken edge is an edge of the original forest $\mathcal F$.
\item Suppose $\mathcal F$ has no isolated vertices and a rotation starts at $v$ and has new endpoint $w$. Then the final forest has endpoint set $(End(\mathcal F)\setminus\{v\})\cup\{w\}$, its only possible isolated vertex is $w$, and $w$ can be isolated only if $w\in End(\mathcal F)$. Also, $v$ and $w$ lie in the same part.
\item Let $E_U^k(v,\mathcal F)\subset S\subset E_U^{k+1}(v,\mathcal F)$ for some $k\geq0$, and for each $x\in S$ fix a shortest rotation in $U$ from $v$ with new endpoint $x$ and added edges from $G_1$. These rotations use at most $2|S|$ pivots.
\end{enumerate}
\end{lemma}
\begin{proof}
Each rotation here is a rotation in the sense of Definitions~3.1 and~3.2 of~\cite{C-expander}, with a particular choice of the new endpoint. Part~(i) is noted in Definition~3.2 of~\cite{C-expander}. We recall why. By the distance condition, a pivot $z$ has the same neighbours as in $\mathcal F$ when its step is performed. Thus the broken edge $zy$ of its step is an edge of $\mathcal F$ if $z$ is not isolated in $\mathcal F$, and $y=z$ otherwise.

For (ii), we may assume that the rotation has at least one step. Then $w\neq v$ by the distance condition. By (B) and (C) of Lemma~3.3 of~\cite{C-expander}, applied with $U=V(\mathcal F)$ in the union of $\mathcal F$ and the given graph, the possible isolated vertex and the endpoint set are as stated. Since no pivot is isolated in $\mathcal F$, each step adds an edge $w'z$ at its old endpoint $w'$ and deletes an edge $zy$. It is easy to see that $w'$ and $y$ are in the same part, since all graphs here are bipartite with parts $V_1,V_2$. The old endpoint of each step is the new endpoint of the previous one. Thus $v$ and $w$ are also in the same part.

For (iii), each fixed rotation has at most $k+1$ steps. For $i\leq k$, its first $i$ steps form a rotation in $U$, and hence the new endpoint of its $i$th step is in $E_U^k(v,\mathcal F)\subset S$. The new endpoint of a $(k+1)$st step is the new endpoint $x\in S$ of the rotation. Hence the new endpoints of all steps are in $S$. By the proof of (i), each pivot is a neighbour in $\mathcal F$ of the new endpoint of its step, or equal to it if the pivot is isolated in $\mathcal F$. There are at most $2|S|$ pivots, since $\Delta(\mathcal F)\leq2$.
\end{proof}

\begin{lemma}\label{lem:endpoint-growth}
Let $\mathcal F$ be a linear forest with an endpoint $v$, and let $U$ be a set containing $v$ that is $(r,D)$-expanding into $int_{\mathcal F}(U)$ in $G_1$. Then $|E_U^K(v,\mathcal F)|\geq20b$.
\end{lemma}
\begin{proof}
We follow the proof of Lemma~3.4 of~\cite{C-expander}. Set $Y_j=E_U^j(v,\mathcal F)$. Since $v\in U$ and every new endpoint other than $v$ is a pivot in $int_{\mathcal F}(U)$ or a neighbour in $\mathcal F$ of one, we have $Y_j\subset U$. For each $x\in Y_j$ fix a shortest rotation in $U$ from $v$ with new endpoint $x$ and added edges from $G_1$. These rotations use at most $2|Y_j|$ pivots by Lemma~\ref{lem:rotation-basics}(iii) with $S=Y_j$. The set of vertices within distance two of these pivots or of $v$ has at most $5(2|Y_j|+1)\leq15|Y_j|$ elements. A pivot outside this set keeps the distance condition when it is added to any of these rotations.

Let $z\in N_{G_1}(Y_j)\cap int_{\mathcal F}(U)$ be outside this set, with a neighbour $w\in Y_j$. We add a step with pivot $z$ to the rotation for $w$, and we obtain a rotation in $U$ with at most $j+1$ steps. Since its new endpoint is a neighbour of $z$ in $\mathcal F$ or $z$ itself, at most two such pivots give the same endpoint. Hence
\begin{equation}\label{eq:endpoint-growth}
 |Y_{j+1}|\geq
 \tfrac12\bigl(|N_{G_1}(Y_j)\cap int_{\mathcal F}(U)|-15|Y_j|\bigr).
\end{equation}
Suppose $|Y_K|<20b$, which implies $|Y_j|<20b$ for $j\leq K$. First suppose that $|Y_j|\geq r$ for some $j<K$, and let $W\subset Y_j$ with $|W|=r$. By expansion, $|N_{G_1}(W)\cap int_{\mathcal F}(U)|\geq Dr$. As $N_{G_1}(W)\setminus Y_j\subset N_{G_1}(Y_j)$, this implies $|N_{G_1}(Y_j)\cap int_{\mathcal F}(U)|\geq Dr-|Y_j|$. Then by~\eqref{eq:endpoint-growth}, $Dr\geq1000b$ and $|Y_j|<20b$, we have $|Y_{j+1}|\geq(Dr-16|Y_j|)/2>(1000b-320b)/2>20b$, which is a contradiction.

Otherwise $|Y_j|<r$ for every $j<K$, so by expansion $|N_{G_1}(Y_j)\cap int_{\mathcal F}(U)|\geq D|Y_j|$, and~\eqref{eq:endpoint-growth} implies $|Y_{j+1}|\geq(D-15)|Y_j|/2\geq\sqrt D\,|Y_j|$. As $|Y_0|=1$, we get $|Y_K|\geq D^{K/2}>N$, again a contradiction.
\end{proof}

\subsection{Independent rotations and moving an endpoint}
The next two lemmas adapt Lemmas~3.5 and~3.6 of~\cite{C-expander}. In the proofs of these lemmas in~\cite{C-expander}, an edge is found by joinedness. We take this edge from $G_2$ instead, outside the set $Z$ given by Lemma~\ref{lem:max-degree}(iii) for $H$, so that it is $(2,b,H)$-bounded. In both lemmas, $\mathcal F$ is a linear forest with no isolated vertices.

\begin{lemma}\label{lem:independent-rotations}
Let $X$ and $Y$ be sets such that $X$ is $(r,D)$-expanding into $int_{\mathcal F}(X)$ and $Y$ into $int_{\mathcal F}(Y)$ in $G_1$, and no component of $\mathcal F$ meets both $X$ and $Y$. Suppose that $End(\mathcal F)\cap(X\cup Y)=\{x,y\}$, where $x\in X$ and $y\in Y$ lie in opposite parts.
Then $G_2$ has an edge $e$ such that $\{e\}$ is $(2,b,H)$-bounded, and $G_1\cup\mathcal F\cup\{e\}$ has a spanning linear forest $\mathcal F'$ with no isolated vertices, endpoint set $End(\mathcal F)\setminus\{x,y\}$, and
$|E(\mathcal F')\triangle E(\mathcal F)|\leq4K+1$.
\end{lemma}
\begin{proof}
By Lemma~\ref{lem:endpoint-growth}, $X'=E_X^K(x,\mathcal F)$ and $Y'=E_Y^K(y,\mathcal F)$ have at least $20b$ vertices each, and by Lemma~\ref{lem:rotation-basics}(ii) they are contained in the parts of $x$ and $y$. Let $Z$, with $|Z|\leq b$, be the set of Lemma~\ref{lem:max-degree}(iii). Since $G_2$ is $b$-bipartite-joined, $G_2$ has an edge $e=x'y'$ with $x'\in X'\setminus Z$ and $y'\in Y'\setminus Z$, and $\{e\}$ is $(2,b,H)$-bounded as $e$ avoids $Z$.

We perform a shortest rotation in $X$ with new endpoint $x'$ and added edges from $G_1$, and then one in $Y$ with new endpoint $y'$. The first changes only edges inside $X$, since its pivots belong to $int_{\mathcal F}(X)$. The second rotation can still be performed and the components of $x'$ and $y'$ are different, as no component meets both $X$ and $Y$. As $End(\mathcal F)\cap(X\cup Y)=\{x,y\}$, we have $x'=x$ or $x'\notin End(\mathcal F)$, and similarly for $y'$. Lemma~\ref{lem:rotation-basics}(ii) then shows that there is no isolated vertex and that the endpoint set is $(End(\mathcal F)\setminus\{x,y\})\cup\{x',y'\}$. Adding $e$ joins two paths, and we obtain the desired $\mathcal F'$. Each rotation changes at most $2K$ edges, and hence $|E(\mathcal F')\triangle E(\mathcal F)|\leq4K+1$.
\end{proof}

\begin{lemma}\label{lem:moving}
Let $U$ be a set that contains an endpoint $u$ of $\mathcal F$ and is $(r,D)$-expanding into $int_{\mathcal F}(U)$ in $G_1$. Let $V$ be a set with no endpoints of $\mathcal F$ such that $int_{\mathcal F}(V)$ has at least $50b$ vertices in each part. Then there are a $(2,b,H)$-bounded set $E\subset E(G_2)$ with $|E|\leq1$ and a rotation from $u$ of at most $K+1$ steps, with added edges from $G_1\cup E$, whose new endpoint $v$ lies in $V$. Moreover, the following hold, where $\mathcal F'$ is the resulting forest.
\begin{enumerate}[label=\upshape(\roman{enumi})]
\item $\mathcal F'$ has no isolated vertices and $End(\mathcal F')=(End(\mathcal F)\setminus\{u\})\cup\{v\}$.
\item If $E=\emptyset$, every pivot lies in $int_{\mathcal F}(U)$. Otherwise, the last step adds the edge of $E$ and has its pivot in $int_{\mathcal F}(V)$, and all other pivots lie in $int_{\mathcal F}(U)$.
\item For every $V'\subset V$, every vertex of $int_{\mathcal F}(V')$ except possibly the last pivot lies in $int_{\mathcal F'}(V')$.
\end{enumerate}
\end{lemma}
\begin{proof}
Write $Y_j=E_U^j(u,\mathcal F)$. If these sets reach $V$ within $K$ steps, we rotate into $V$ with edges of $G_1$ only. Otherwise we use one edge of $G_2$ for the last step.

\textbf{Case 1: $Y_K$ meets $V$.} Let $j$ be least with $Y_j\cap V\neq\emptyset$. We take a rotation in $U$ from $u$ with at most $j$ steps, added edges from $G_1$ and new endpoint $v\in Y_j\cap V$, and we let $E=\emptyset$.

\textbf{Case 2: $Y_K\cap V=\emptyset$.} By Lemma~\ref{lem:endpoint-growth}, $|Y_K|\geq20b$. Let $j$ be least with $|Y_j|\geq2b$, and choose $T$ with $|T|=2b$ and $Y_{j-1}\subset T\subset Y_j$. For each $x\in T$ fix a shortest rotation in $U$ from $u$ with new endpoint $x$ and added edges from $G_1$. Lemma~\ref{lem:rotation-basics}(ii) and~(iii), with $S=T$, imply that $T$ is a subset of the part of $u$ and that these rotations use at most $4b$ pivots.

Let $Z$ be the set of Lemma~\ref{lem:max-degree}(iii). Note that $T\setminus Z$ has at least $b$ vertices, as $|Z|\leq b$. We delete from $int_{\mathcal F}(V)$ the set $Z$ and the at most $5(4b+1)=20b+5$ vertices within distance two of these pivots or of $u$. Then each remaining vertex keeps the distance condition as a new pivot. At least $50b-b-(20b+5)\geq b$ vertices remain in $int_{\mathcal F}(V)$ in the other part. Since $G_2$ is $b$-bipartite-joined, we can choose an edge $xz\in E(G_2)$ with $x\in T\setminus Z$ and $z$ among these vertices. Then $E=\{xz\}$ is $(2,b,H)$-bounded, as $xz$ avoids $Z$. Adding a step with pivot $z$ to the rotation for $x$ yields a rotation whose new endpoint $v$ is in $N_{\mathcal F}(z)\subset V$, because its broken edge $zv$ is an edge of $\mathcal F$.

For (i), as $v\in V$, we have $v\notin End(\mathcal F)$, and hence Lemma~\ref{lem:rotation-basics}(ii) shows that the endpoint set is as stated and there is no isolated vertex. Part (ii) holds by construction.

For (iii), we first note that in both cases no old endpoint of a step is in $V$. Indeed, $u\notin V$ as $V\cap End(\mathcal F)=\emptyset$. The old endpoint of each later step is the new endpoint of a rotation in $U$ with added edges from $G_1$, and thus it belongs to $Y_{j-1}$ in Case~1 and to $Y_K$ in Case~2. The set $Y_{j-1}$ is disjoint from $V$ by the choice of $j$ in Case~1, and $Y_K\cap V=\emptyset$ in Case~2.

Now let $V'\subset V$. A vertex of $int_{\mathcal F}(V')$ that is not in $int_{\mathcal F'}(V')$ gains a neighbour outside $V'$. Only pivots and old endpoints gain neighbours. Since no old endpoint is in $V$, such a vertex is a pivot. Let $z$ be a pivot other than the last. The new endpoint of its step is the old endpoint of the next step and thus is not in $V$. Since it is a neighbour of $z$ in $\mathcal F$ by Lemma~\ref{lem:rotation-basics}(i), we have $z\notin int_{\mathcal F}(V')$. This proves (iii).
\end{proof}

\section{Forest reductions}\label{app:forests}
We now prove Lemmas~\ref{lemma:linearforest decomp1} and~\ref{lemma:linearforest decomp2}. For Lemma~\ref{lemma:linearforest decomp1}, we take a forest that minimizes a potential $Q$. We show that it has at most $n^{0.8}$ paths, as otherwise we can lower $Q$ further by adding a few edges of $G_2$. Some of these changes do not join paths and only make the orders of the paths more even, but each still uses an edge of $G_2$. Since $Q$ also counts the logarithms of the orders of the paths, these changes lower $Q$ too, and we can bound their number. In Lemma~\ref{lem:two-endpoints} we remove two endpoints, one in each part, with a few edges of $G_2$, following the proof of Lemma~3.7 of~\cite{C-expander}. Lemma~\ref{lemma:linearforest decomp2} follows by applying Lemma~\ref{lem:two-endpoints} repeatedly.

\subsection{The potential \texorpdfstring{$Q$}{Q}}
For a spanning linear forest $\mathcal F$ with $k(\mathcal F)$ components, let
\[
 \Lambda(\mathcal F)=\sum_{P\in\mathcal F}\log |V(P)|,
 \quad\text{and}\quad
 Q(\mathcal F)=40k(\mathcal F)-\frac{\Lambda(\mathcal F)}{\log N}.
\]
Note that as $0\leq \Lambda(\mathcal F)\leq k(\mathcal F)\log N$, we have $39k(\mathcal F)\leq Q(\mathcal F)\leq40k(\mathcal F)$. Thus $Q$ is about $40$ times the number of paths. Among forests with the same number of paths, $Q$ is smaller when the orders of the paths are more even.

The next lemma shows that $\Lambda$ drops by at most $\log2$ per step of a rotation, apart from the change in the order of the path that contains the current endpoint.

\begin{lemma}\label{lem:log-orders}
Suppose that in a rotation with $j$ steps, no pivot is an endpoint of another component. Let $\mathcal F_i$ be the forest after $i$ steps, and let $a_i$ be the order of the component of $\mathcal F_i$ that contains the new endpoint of the $i$th step, or the starting endpoint if $i=0$. Then the rotation preserves the number of components and satisfies
\begin{equation}\label{eq:log-orders}
 \Lambda(\mathcal F_j)\geq \Lambda(\mathcal F_0)+\log(a_j/a_0)-j\log2.
\end{equation}
Moreover, joining two different components by an edge decreases $\Lambda$ by at most $\log N$.
\end{lemma}
\begin{proof}
We show that $\Lambda(\mathcal F_i)-\log a_i$ decreases by at most $\log2$ in each step. Consider a step whose pivot and old endpoint are in different components. Let $a$ and $a'$ be the orders of the component of the old endpoint and of the component of the pivot. The pivot is internal by hypothesis. The step adds the edge from the old endpoint to the pivot and deletes an edge at the pivot. The pivot and one branch at it then join the component of the old endpoint, and the other branch becomes a separate component. By the definition of a $1$-rotation, this other branch contains the new endpoint and has order $c\leq(a'-1)/2$. The other new component has order $a+a'-c>a'/2$. As $a_i=a$ and $a_{i+1}=c$, we get
\[
 \bigl(\Lambda(\mathcal F_{i+1})-\log c\bigr)
 -\bigl(\Lambda(\mathcal F_i)-\log a\bigr)
 =\log\left(\frac{a+a'-c}{a'}\right)\geq-\log2.
\]
A step within one component replaces its path by a path on the same vertex set. The step does not change $\Lambda$, and $a_{i+1}=a_i$. In both cases the step preserves the number of components, and summing over the $j$ steps proves~\eqref{eq:log-orders}. Also, joining components of orders $a,a'\leq N$ changes $\Lambda$ by $\log(1/a+1/a')\geq-\log N$.
\end{proof}

\subsection{Good families}
The proof of Lemma~3.7 in~\cite{C-expander} uses rotations inside five families $\mathcal{H}_1,\dots,\mathcal{H}_4,\mathcal{H}_{\rm hop}$ of paths of length at most $\sqrt n$. We need such families too, with $\mathcal H_5$ in place of $\mathcal H_{\rm hop}$, and we find them in this subsection. They must have the following property, with $b$, $r$ and $D$ as in Appendix~\ref{app:rotations}.

\begin{definition}\label{def:good}
Let $\mathcal H_1,\ldots,\mathcal H_5$ be disjoint families of vertex-disjoint paths. We say that they are \emph{good} if, for all sets $W_1,\ldots,W_5$ with $W_i\subset V(\mathcal H_i)$ and $|W_i|\leq20b$, there are sets $U_i\subset V(\mathcal H_i)\setminus W_i$ such that, for all $i\in[5]$,
\begin{enumerate}[label=\upshape(\roman{enumi})]
\item $U_i$ is $(r,4D)$-expanding into $int_{\mathcal H_j}(U_j)$ in $G_1$ for all $j\in[5]$, and
\item each part contains at least $100b$ vertices of $int_{\mathcal H_i}(U_i)$.
\end{enumerate}
\end{definition}

In~\cite{C-expander}, the expanding subsets of the five families come from Lemma~2.5 of~\cite{C-expander}, which uses that a $C$-expander has an edge between any two disjoint sets of size $n/(2C)$. Our $G_1$ does not necessarily have this property, and rotating with edges of $G_2$ instead could use too many edges of $G_2$. We use the next lemma in its place. If paths of length at most $\sqrt n$ cover almost all vertices, it splits them into five families whose vertex sets contain expanding subsets in $G_1$. Lemma~\ref{lem:random-good} then shows that these families are good. In Appendix~\ref{app:bipartite-connectivity}, good families come from Lemma~\ref{lem:expanding-sets} instead.

\begin{lemma}\label{cl: partition into five}
Let $1/d_0\ll\eta,\beta\ll\alpha\ll1$, and let $n\in\mathbb N$ and $d\geq d_0$. Let $G_1$ be an $(\eta,\beta,d/n)$-sparse bipartite graph with parts $V_1,V_2$ of size $n$ and $\delta(G_1)\geq\alpha d$. Let $P_1,\ldots,P_t$ be vertex-disjoint paths of lengths between $1$ and $\sqrt n$ that respect this bipartition and together cover at least $(2-\alpha^{10})n$ vertices. Their edges do not have to belong to $G_1$. Then there is a partition $\mathcal{H}_{1},\dots, \mathcal{H}_{5}$ of $\{P_1,\dots, P_t\}$ with the following property. For all sets $W_1,\dots,W_5$ with $W_i\subset V(\mathcal{H}_i)$ and $|W_i|\leq\alpha^{10} n$, there exist sets $U_1,\dots,U_5$ with $U_i \subset V(\mathcal{H}_i)\setminus W_i$ such that for all $i\in [5]$,
 \begin{enumerate}[label= \upshape{(\roman{enumi})}]
     \item $U_i$ is $(\beta n,1/\alpha)$-expanding into $int_{\mathcal{H}_j}(U_j)$ in $G_1$ for all $j\in [5]$, and
     \item $|int_{\mathcal{H}_i}(U_i)\cap V_c|\geq\alpha^5n/4$ for $c=1,2$.
 \end{enumerate}
\end{lemma}
\begin{proof}
We choose the families at random, so that every vertex outside a set $B$ of at most $\alpha^{10}n$ vertices has at least $\alpha^5d$ neighbours in each family. Then we remove the vertices with less than $\alpha^5d/2$ neighbours in some interior, one at a time. After that, Lemma~\ref{lem:sparse->expander} gives the expansion.

\paragraph{The random partition.} Consider independent random variables $z_1,\dots, z_t$ taking values in $[5]$ uniformly at random, and let $\mathcal{H}_i=\{P_k:z_k=i\}$. Let $S=\bigcup_{k\in[t]}V(P_k)$, so that $|V(G_1)\setminus S|\leq\alpha^{10}n$. Let $B_0$ be the set of vertices $x$ with
\[
 |N_{G_1}(x)\setminus S|\geq\alpha d/2\quad\text{or}\quad|N_{G_1}(x)\cap V(P_k)|\geq\alpha^4d\text{ for some }k\in[t].
\]
A vertex outside $B_0$ has many neighbours in $S$, spread over many paths, so the random partition is likely to give it many neighbours in each family.

We first show that $|B_0|\leq\alpha^{10}n/2$. For this we apply Lemma~\ref{lem:few-heavy} to $G_1$ in two ways, with $2n$ playing the role of $N$. The hypothesis of Lemma~\ref{lem:few-heavy} on $h$ holds in both uses as $\eta\ll\alpha$ and $n\geq\delta(G_1)\geq\alpha d_0$ is large.
\begin{itemize}
\item With $h=\alpha d/2$ and $T=V(G_1)\setminus S$, it shows that at most $8\eta\alpha^9n$ vertices are of the first kind.
\item With $h=\alpha^4d$ and $T=V(P_k)$, it shows that for each $k\in[t]$ at most $4\eta|V(P_k)|/\alpha^4$ vertices have at least $\alpha^4d$ neighbours in $V(P_k)$.
\end{itemize}
Summing the second bound over $k\in[t]$ and using $\sum_{k}|V(P_k)|\leq2n$, we get $|B_0|\leq8\eta\alpha^9n+8\eta n/\alpha^4\leq\alpha^{10}n/2$.

Let $x\notin B_0$, and set $w_k=|N_{G_1}(x)\cap V(P_k)|$ for $k\in[t]$ and $w=\sum_{k\in[t]}w_k$. As $x\notin B_0$, we have $w\geq\alpha d/2$ and $w_k<\alpha^4d$ for every $k$. For $j\in[5]$, let
\[
 Y_j=|N_{G_1}(x)\cap V(\mathcal H_j)|=\sum_{k\in[t]}w_k\mathbf1_{\{z_k=j\}}.
\]
Then $\mathbb E[Y_j]=w/5\geq2\alpha^5d$ and $\sum_{k}w_k^2\leq\alpha^4dw$. McDiarmid's inequality (Lemma~\ref{lem:mcdiarmid}), with $c_k=w_k$ and with deviation $w/10$, gives
\[\mathbb P\left(Y_j<\alpha^5d\right)\leq2\exp\left(-\frac{w}{200\alpha^4d}\right)\leq2\exp\left(-\frac1{400\alpha^3}\right)\leq\frac{\alpha^{10}}{20}.\]
Summing over the $2n$ vertices and the five values of $j$, the expected number of vertices outside $B_0$ with less than $\alpha^5d$ neighbours in some $V(\mathcal H_j)$ is at most $2n\cdot5\cdot\alpha^{10}/20=\alpha^{10}n/2$. Fix an outcome in which there are at most $\alpha^{10}n/2$ such vertices, and let $B$ consist of $B_0$ and these vertices. Then $|B|\leq\alpha^{10}n$.

\paragraph{The sets $U_1,\dots,U_5$.} We show that $\mathcal{H}_1,\ldots, \mathcal{H}_5$ are as required. Let $W_1,\dots, W_5$ be fixed sets of size at most $\alpha^{10} n$, with $W_i\subset V(\mathcal H_i)$ for $i\in [5]$. Then $W^0=B\cup\bigcup_{i\in[5]}W_i$ has size at most $|B|+5\alpha^{10}n\leq10\alpha^{10}n$. For each $i\in [5]$, let $U^0_i = V(\mathcal{H}_i)\setminus W^0$ and note that, as $U^0_i\cap B=\emptyset$,
\begin{equation}\label{eq: degree}
 |N_{G_1}(x)\cap V(\mathcal H_j)|\geq\alpha^5d\quad\text{for all }x\in U^0_i\text{ and }j\in[5].
\end{equation}

We now add vertices one at a time, starting from $W^0$. Let $U^s_1,\ldots,U^s_5$ and $W^s$ be given. Suppose that for some $i,j\in[5]$ a vertex $x\in U^s_i$ has less than $\alpha^5d/2$ neighbours in $int_{\mathcal H_j}(U^s_j)$. We move $x$ from $U^s_i$ to $W^s$ and call the new sets $U^{s+1}_1,\ldots,U^{s+1}_5$ and $W^{s+1}$. The process stops after some $\ell$ steps, when there is no such vertex.

\begin{claim}\label{cl:W-small}
$|W^\ell|<20\alpha^{10}n$.
\end{claim}
\begin{claimproof}[Proof of Claim~\ref{cl:W-small}]
Suppose not, and let $s\le\ell$ be least with $|W^s|\geq20\alpha^{10}n$. Then $|W^s|<20\alpha^{10}n+1$, as the sets $W^s$ grow by one vertex at a time. Every vertex $x\in W^s\setminus W^0$ has less than $\alpha^5 d/2$ neighbours in $int_{\mathcal{H}_j}(U^s_j)$ for some $j\in[5]$, since interiors only shrink during the process. The set $W^s\cup N_{\mathcal H_j}(W^s)$ contains every vertex of $V(\mathcal H_j)$ outside $int_{\mathcal H_j}(U^s_j)$, as $U^s_j=V(\mathcal H_j)\setminus W^s$. Hence, by~\eqref{eq: degree}, $x$ has at least $\alpha^5d/2$ neighbours in $T=W^s\cup\bigcup_{j\in[5]}N_{\mathcal H_j}(W^s)$. We have $|T|\le11|W^s|\le300\alpha^{10}n$, as $|N_{\mathcal H_j}(W^s)|\le2|W^s|$ for each $j$ and $n$ is large. Lemma~\ref{lem:few-heavy}, applied as above with $h=\alpha^5d/2$ and this $T$ (its hypothesis on $h$ holds as $\eta\ll\alpha\ll1$), shows that $|W^s\setminus W^0|\le8\eta|T|/\alpha^5$. Hence
\[
 10\alpha^{10}n\le|W^s|-|W^0|\le|W^s\setminus W^0|\le8\eta|T|/\alpha^5\le2400\eta\alpha^5n,
\]
where the first inequality holds as $|W^s|\geq20\alpha^{10}n$ and $|W^0|\le10\alpha^{10}n$, and the last as $|T|\le300\alpha^{10}n$. This contradicts $\eta\ll\alpha$.
\end{claimproof}

Let $U_i:= V(\mathcal{H}_i)\setminus W^\ell$ for $i\in[5]$. When the process stops, every vertex $x\in U_i$ has at least $\alpha^5d/2$ neighbours in $int_{\mathcal{H}_j}(U_j)$, for all $i,j\in[5]$. Lemma~\ref{lem:sparse->expander}, with $2n$, $2d$, $\alpha^5/4$ and $1/\alpha$ playing the roles of $n$, $d$, $\alpha$ and $D$, now shows that every subset $X\subset U_i$ with $|X|\leq\beta n\leq\alpha^6n/8$ satisfies $|N_{G_1}(X)\cap int_{\mathcal{H}_j}(U_j)|\geq|X|/\alpha$. This proves (i).

Let $i\in[5]$ and $c\in\{1,2\}$, and let $X_c=(U_1\cup\dots\cup U_5)\cap V_c$. As $|V(G_1)\setminus S|\leq\alpha^{10}n$ and $|W^\ell|<20\alpha^{10}n$ by Claim~\ref{cl:W-small}, we have $|X_c|\geq(1-21\alpha^{10})n\geq2\eta n$. Every vertex of $X_c$ has at least $\alpha^5d/2$ neighbours in $int_{\mathcal H_i}(U_i)$, and they all belong to $V_{3-c}$. We apply the second assertion of Lemma~\ref{lem:sparse->expander} with the same parameters as for (i), and with $X_c$ and $int_{\mathcal H_i}(U_i)$ playing the roles of $X$ and $Y$. This gives $|int_{\mathcal H_i}(U_i)\cap V_{3-c}|\geq\alpha^5n/4$, and (ii) follows.
\end{proof}

\begin{lemma}\label{lem:random-good}
Let $G_1$ be as in Lemmas~\ref{lemma:linearforest decomp1} and~\ref{lemma:linearforest decomp2}. Then $G_1$ is $(r,4D)$-expanding. Moreover, let $n/10\leq m\leq n$ and let $G'$ be a subgraph of $G_1$ whose parts have size $m$ and lie in $V_1$ and $V_2$, such that $G'$ is $(\eta,\beta,d/m)$-sparse with $\delta(G')\geq\alpha d$. Then, for any paths as in Lemma~\ref{cl: partition into five} for $G'$ with $m$ in place of $n$, the families that this lemma gives are good.
\end{lemma}
\begin{proof}
Here $G_1$ is $(\eta,\beta,d/n)$-sparse with $\delta(G_1)\geq\alpha d$, and $n\geq\delta(G_1)\geq\alpha d_0$ is large. As $\eta,\beta\ll\alpha\ll1$, we may assume that $1/\alpha\geq4D$ and $1000b\leq\alpha^{10}n$. Lemma~\ref{lem:sparse->expander}, with $X=Y=V_1\cup V_2$ and with $2n$, $2d$, $\alpha/2$ and $4D$ playing the roles of $n$, $d$, $\alpha$ and $D$, shows that $G_1$ is $(r,4D)$-expanding, as $\eta,\beta\ll\alpha$.

Now let $m$ and $G'$ be as in the statement, and apply Lemma~\ref{cl: partition into five} to $G'$ with $m$ playing the role of $n$. The sets $W_i$ of Definition~\ref{def:good} are allowed in Lemma~\ref{cl: partition into five}, since $20b\leq\alpha^{10}n/50\leq\alpha^{10}m$. Also $4D\leq1/\alpha$, $100b\leq\alpha^{10}n/10\leq\alpha^5m/4$, and $r\leq\beta n/10\leq\beta m$ as $D=2^{100}$ and $n$ is large. Hence (i) and (ii) of Lemma~\ref{cl: partition into five} give (i) and (ii) of Definition~\ref{def:good}, as expansion in a subgraph of $G_1$ implies expansion in $G_1$.
\end{proof}

\subsection{Removing two endpoints}
In Lemma~\ref{lem:two-endpoints} we remove two endpoints $x,y$ in opposite parts. We move $x$ and $y$ by rotations into sets given by good families, and then we join them.

\begin{figure}[htbp]
  \centering
  \begin{tikzpicture}[x=1cm,y=1cm,
      v/.style={circle,fill=black,inner sep=1pt},
      lab/.style={font=\scriptsize,inner sep=2.5pt},
      box/.style={draw,rounded corners,minimum width=1.6cm,minimum height=.9cm},
      wavy/.style={->,decorate,decoration={snake,amplitude=.3mm,segment length=1.6mm,post length=1.5mm},shorten >=2pt},
      pa/.style={DarkDesaturatedBlue,line width=.8pt},
      pb/.style={VividOrange,line width=.8pt}]
    \node[box] (B1) at (2.5,1.6) {}; \node[lab,above] at (B1.north) {$U_1$};
    \node[box] (B2) at (2.5,0) {};   \node[lab,below] at (B2.south) {$U_2$};
    \node[box] (B5) at (5.5,1.6) {}; \node[lab,above] at (B5.north) {$U_5$};
    \node[box] (B3) at (8.5,1.6) {}; \node[lab,above] at (B3.north) {$U_3$};
    \node[box] (B4) at (8.5,0) {};   \node[lab,below] at (B4.south) {$U_4$};
    \node[v] (x) at (0,1.6) {}; \node[lab,left] at (x) {$x$};
    \node[v] (y) at (0,0) {};   \node[lab,left] at (y) {$y$};
    \node[v] (x1) at (2.5,1.6) {}; \node[lab,above] at (x1) {$x_1$};
    \node[v] (x2) at (2.5,0) {};   \node[lab,below] at (x2) {$x_2$};
    \node[v] (h) at (5.5,1.6) {};  \node[lab,above] at (h) {$h$};
    \node[v] (x3) at (8.5,1.6) {}; \node[lab,above left] at (x3) {$x_3$};
    \node[v] (x4) at (8.5,0) {};   \node[lab,below left] at (x4) {$x_4$};
    \draw[pa,wavy] (x)--(x1);
    \draw[pb,wavy] (y)--(x2);
    \draw[pa,wavy] (x1)--(h);
    \draw[pa,->,shorten >=2pt] (h)--(x3);
    \draw[pb,wavy] (x2)--(x4);
    \draw[densely dashed] (x3) to[bend left=40] node[lab,right] {join} (x4);
  \end{tikzpicture}
  \caption{The moves in the proof of Lemma~\ref{lem:two-endpoints}(i), those of $x$ in blue and those of $y$ in orange. Wavy arrows are moves by Lemma~\ref{lem:moving}, the straight arrow is a $1$-rotation, and the dashed line is the joining by Lemma~\ref{lem:independent-rotations}.}
  \label{fig:two-endpoints}
\end{figure}

\begin{lemma}\label{lem:two-endpoints}
Let $\mathcal F$ be a linear forest with no isolated vertices and at most $6b$ components, and suppose $G_1$ is $(r,4D)$-expanding. Let $\mathcal H_1,\ldots,\mathcal H_5$ be good families of components of $\mathcal F$ of length at most $\sqrt n$ each. Let $x,y\in End(\mathcal F)$ be in opposite parts. Then the following hold.
\begin{enumerate}[label=\upshape(\roman{enumi})]
\item There is a $(10,b,H)$-bounded set $E\subset E(G_2)$ such that $G_1\cup\mathcal F\cup E$ contains a spanning linear forest $\mathcal F'$ with no isolated vertices and
\[
 End(\mathcal F')=End(\mathcal F)\setminus\{x,y\},
 \qquad |E(\mathcal F')\triangle E(\mathcal F)|\leq10\log n.
\]
\item There is a $(10,b,H)$-bounded set $E\subset E(G_2)$ such that $G_1\cup\mathcal F\cup E$ contains a spanning linear forest $\mathcal F'$ with $Q(\mathcal F')\leq Q(\mathcal F)-20$.
\end{enumerate}
\end{lemma}
\begin{proof}
We first prove (i). We move $x$ and $y$ to new endpoints $x_3$ and $x_4$ in sets $U_3$ and $U_4$ of the third and fourth families such that no component meets both $U_3$ and $U_4$. Then we join them by Lemma~\ref{lem:independent-rotations}, which needs this condition (see Figure~\ref{fig:two-endpoints}). The first moves go into sets $U_1$ and $U_2$ of the first two families, and they may change paths anywhere. We therefore choose $U_3$, $U_4$ and a set $U_5$ of the fifth family only afterwards, on paths that are still untouched.

We have $|End(\mathcal F)|\leq12b\leq20b$, as $\mathcal F$ has at most $6b$ components. We apply the definition of good families with $W_i=End(\mathcal F)\cap V(\mathcal H_i)$, and let $U_1,U_2$ be the first two sets it gives, which contain no endpoint of $\mathcal F$.

\paragraph{Moving $x$ and $y$ into $U_1$ and $U_2$.}
We first apply Lemma~\ref{lem:moving} with $U=V_1\cup V_2$ and $V=U_1\cup U_2$ to move $x$ to $x_1$. This is possible, as $int_{\mathcal F}(V_1\cup V_2)=V_1\cup V_2$, $G_1$ is $(r,4D)$-expanding, and $int_{\mathcal F}(V)$ has at least $100b$ vertices in each part by property~(ii) of good families. We take $V=U_1\cup U_2$ rather than $U_1$, so that Lemma~\ref{lem:moving}(iii) applies with $V'=U_1$ and with $V'=U_2$. If necessary, we interchange the first two families so that $x_1\in U_1$.

We then move $y$ with $U=(V_1\cup V_2)\setminus\{x_1\}$ and $V=(U_1\cup U_2)\setminus\{x_1\}$, obtaining $w$. The new $U$ is still $(r,D)$-expanding into its interior, since deleting $x_1$ removes at most three interior vertices. The new $V$ contains no endpoint. Lemma~\ref{lem:moving}(iii) for the first move implies that its interior has at least $100b-4\geq50b$ vertices in each part.

If $w\in U_2$, let $x_2=w$. Otherwise $w\in U_1$, and we add a $1$-rotation from $w$ into $U_2$. By property~(i) of good families, $w$ has at least $4D$ neighbours in the original interior of $U_2$. At most two of them have left the current interior, at most one in each move by Lemma~\ref{lem:moving}(iii). At most two others lie within distance two of $w$ in the current forest, and they cannot be pivots by the distance condition. At least $4D-4>0$ neighbours remain, and we perform the $1$-rotation with one of them as pivot. The new endpoint $x_2$ is a neighbour of this pivot in the current forest and therefore belongs to $U_2$. Moreover, $x_2$ is not isolated, since $U_2$ contains no current endpoint. The step changes degrees only at $w$ and $x_2$.

By Lemma~\ref{lem:moving}(i) and this step, the resulting forest $\mathcal F_0$ has no isolated vertices and $End(\mathcal F_0)=(End(\mathcal F)\setminus\{x,y\})\cup\{x_1,x_2\}$. Each of the two moves has at most $K+1$ steps, each step changes at most two edges, and the $1$-rotation changes at most two more. Hence $\mathcal F_0$ differs from $\mathcal F$ in at most $2(2K+2)+2=4K+6$ edges. 

We need $U_1$ and $U_2$ to stay expanding, as the next moves start in them. Each of the two moves removes at most one original interior vertex from each of $U_1,U_2$ by Lemma~\ref{lem:moving}(iii), and the $1$-rotation removes at most one more. In the second move, Lemma~\ref{lem:moving}(iii) applies to $U_1$ only through $V'=U_1\setminus\{x_1\}$. Deleting $x_1$ removes at most three interior vertices, and hence $U_1$ may lose at most three more. Thus each of $U_1,U_2$ has lost at most six original interior vertices. As $4D|T|-6\geq D|T|$ for every nonempty $T$, both are $(r,D)$-expanding into their current interiors.

\paragraph{Moving the endpoints into $U_3$ and $U_4$.}
In the third, fourth and fifth families, we discard every original path incident with an edge of $E(\mathcal F_0)\triangle E(\mathcal F)$. Each such edge is incident with at most two paths, and thus at most $8K+12$ paths with at most $(8K+12)(\sqrt n+1)\leq b$ vertices are discarded. Since $b+12b\leq20b$, we may apply the definition of good families again, with $W_i$ the set of vertices of $V(\mathcal H_i)$ that lie on discarded paths or are endpoints of $\mathcal F$. Let $U_3,U_4,U_5$ be the third, fourth and fifth sets it gives. These sets are contained in original components that are untouched so far, and they contain no endpoint of $\mathcal F_0$. Their interiors and expansion properties in $\mathcal F_0$ are therefore as in $\mathcal F$.

We move $x_1$ to a vertex $h$ in $U_5$, then $x_2$ into $U_4$, and last $h$ into $U_3$ by a $1$-rotation. We need the move into $U_5$ when, for example, $x_1$ and $x_2$ lie on one path $R$. A direct move of $x_1$ into $U_3$ could attach a piece of a path of $\mathcal H_3$ to $R$. The move of $x_2$ into $U_4$ could then attach a piece of a path of $\mathcal H_4$ to the same component, which would then meet both $U_3$ and $U_4$. 

First we move $x_1$ with $U=U_1$ and $V=U_5$, obtaining $h$. This move changes edges only inside $U_1\cup U_5$, and the component of $h$ is a branch of one path of $\mathcal H_5$. Indeed, by Lemma~\ref{lem:moving}(ii), every pivot except possibly the last is in $int_{\mathcal F_0}(U_1)$. The rotation starts at $x_1\in U_1$, and the steps with their pivots in $int_{\mathcal F_0}(U_1)$ change only edges inside $U_1$. The new endpoints of these steps are in $U_1$, which is disjoint from $U_5$. Since $h\in U_5$, the last step does not have its pivot in $int_{\mathcal F_0}(U_1)$, which implies that it uses the edge of $G_2$ and has its pivot in the interior of $U_5$, on an untouched path $P$ of $\mathcal H_5$. The component of its old endpoint meets $U_1$, so it is not $P$. Thus the component containing $h$ is a branch of $P$ cut off at the pivot. Since all changed edges lie inside $U_1\cup U_5$, the forest is unchanged at $U_2$ and on every component meeting $U_3$ or $U_4$.

Next we move $x_2$ with $U=U_2$ and $V=U_4$, obtaining $x_4$. In the same way, the component of $x_4$ is a branch of an untouched path of $\mathcal H_4$, and all changed edges are inside $U_2\cup U_4$. Hence the component of $h$ and the paths of $\mathcal H_3$ meeting $U_3$ are unchanged. By property~(i) for the new sets, $h\in U_5$ has a neighbour in the interior of $U_3$. We perform a $1$-rotation with such a neighbour as pivot and get $x_3\in U_3$. Since the pivot is not in the component of $h$, it satisfies the distance condition. As before, $x_3$ is not isolated.

We now check the hypotheses of Lemma~\ref{lem:independent-rotations} for $X=U_3$ and $Y=U_4$.
\begin{itemize}
\item No component meets both $U_3$ and $U_4$. Indeed, before this $1$-rotation, the components meeting $U_3$ were untouched paths of $\mathcal H_3$. This rotation splits one of them and attaches one piece to the component of $h$, which is part of a path of $\mathcal H_5$. The other piece stays a separate component.
\item Each of $U_3,U_4$ has lost at most one original interior vertex, a pivot, and thus both are still $(r,D)$-expanding into their current interiors.
\item The only endpoints in $U_3,U_4$ are $x_3,x_4$. They are in the parts of $x,y$ by Lemma~\ref{lem:rotation-basics}(ii), that is, in opposite parts.
\end{itemize}

The current endpoint set is $(End(\mathcal F)\setminus\{x,y\})\cup\{x_3,x_4\}$, so Lemma~\ref{lem:independent-rotations} yields $\mathcal F'$ with no isolated vertices and $End(\mathcal F')=End(\mathcal F)\setminus\{x,y\}$. The four moves each select at most one edge of $G_2$, and the joining edge is one more. Thus $|E|\leq5$. In each of these applications of Lemmas~\ref{lem:moving} and~\ref{lem:independent-rotations}, we use $H$ together with the previously selected edges in place of $H$. Hence each selected edge increases $\Delta_b$ by at most two when added to $H$ and the previously selected edges, and Lemma~\ref{lem:max-degree}(iv) shows that $E$ is $(10,b,H)$-bounded. 

Reaching $\mathcal F_0$ changed at most $4K+6$ edges. The moves into $U_5$ and $U_4$ changed at most $2(K+1)$ edges each, the $1$-rotation changed two, and the joining by Lemma~\ref{lem:independent-rotations} changed at most $4K+1$. In total,
\[
 |E(\mathcal F')\triangle E(\mathcal F)|\leq(4K+6)+2\cdot2(K+1)+2+(4K+1)=12K+13\leq10\log n,
\]
as $K=\lceil2\log_DN\rceil+2$, $D=2^{100}$ and $n$ is large.

\paragraph{Proof of (ii).}
It is enough to add one edge that joins two components while $\Lambda$ drops by at most $20\log N$. We follow the construction above, with one change. If the pivot of a step is an endpoint of another component, we add its edge without deleting an edge, and we stop. Hence no pivot of a step performed before stopping is an endpoint of another component, as Lemma~\ref{lem:log-orders} needs.

There are at most eight rotations, namely the four moves, at most two $1$-rotations and the two rotations in the proof of Lemma~\ref{lem:independent-rotations}. Each has at most $K+1$ steps, and $\log(a_j/a_0)\geq-\log N$ as $1\leq a_i\leq N$. Hence, by Lemma~\ref{lem:log-orders}, each rotation, or the part of it performed before we stop, preserves the number of components and decreases $\Lambda$ by at most $\log N+(K+1)\log2$. Exactly one edge that joins two components is added, whether or not we stop early. This edge reduces the number of components by one and, by Lemma~\ref{lem:log-orders}, decreases $\Lambda$ by at most another $\log N$. In total, as $K\log2\leq\log N$, $\Lambda$ decreases by at most
\[
 9\log N+8(K+1)\log2\leq20\log N.
\]
Let $E$ be the set of the at most five edges of $G_2$ selected before stopping. As in the proof of (i), $E$ is $(10,b,H)$-bounded, and the resulting forest is a subgraph of $G_1\cup\mathcal F\cup E$. This forest has one component fewer and $\Lambda$ at least $\Lambda(\mathcal F)-20\log N$. Thus its value of $Q$ is at most $Q(\mathcal F)-40+20=Q(\mathcal F)-20$.
\end{proof}

\subsection{Reducing the number of paths}
Here we adapt Section~4 of~\cite{C-expander}. We replace the lexicographic order by $Q$, and Lemmas~4.5--4.7 there by Lemma~\ref{lem:long-paths}. In the proof of Lemma~\ref{lem:few-paths}, we take a forest that minimizes $Q$, and we get a contradiction by lowering $Q$ further in one of three ways.
\begin{itemize}
\item If there are more than $6b$ paths, an edge of $G_2$ joins two of them (Lemma~\ref{lem:endpoint-count}).
\item If long paths cover many vertices, a rotation from an end of a path of minimum order ends inside a long path. This makes the orders of the paths more even (Lemma~\ref{lem:long-paths}).
\item Otherwise short paths cover almost all vertices, and we remove two endpoints (Lemma~\ref{lem:two-endpoints}(ii)).
\end{itemize}
Recall that $Q(\mathcal F)=40k(\mathcal F)-\Lambda(\mathcal F)/\log N$. Joining two paths by an edge lowers $Q$ by at least $39$ (Lemma~\ref{lem:log-orders}), and a rotation that keeps the number of paths lowers $Q$ when it increases $\Lambda$. In the second and third ways we also use that a forest that minimizes $Q$ has no isolated vertices if it has at most $6b$ paths (Lemma~\ref{lem:no-isolated}). At the end of this subsection we prove Lemma~\ref{lemma:linearforest decomp1}.

\begin{lemma}\label{lem:endpoint-count}
Let $\mathcal F$ be a spanning linear forest with $k$ components.
\begin{enumerate}[label=\upshape(\roman{enumi})]
\item Each part contains at least $k/2$ endpoints of $\mathcal F$, and exactly $k$ if $\mathcal F$ has no isolated vertices.
\item Suppose $k>6b$. Then $G_2$ has an edge $e$ between endpoints of two different components of $\mathcal F$ such that $\{e\}$ is $(2,b,H)$-bounded. Moreover, $\mathcal F+e$ is a spanning linear forest and $Q(\mathcal F+e)\leq Q(\mathcal F)-39$.
\end{enumerate}
\end{lemma}
\begin{proof}
The sum of $2-d_{\mathcal F}(v)$ over either part is $2n-e(\mathcal F)=k$. An endpoint contributes one to it, or two if it is isolated, and all other vertices contribute zero. This proves (i). 

For (ii), let $Z$ be the set of Lemma~\ref{lem:max-degree}(iii). By (i), $V_1\setminus Z$ contains at least $k/2-b>2b$ endpoints, and we choose $b$ of them. A path with an endpoint in $V_1$ has at most one endpoint in $V_2$. The components of the chosen endpoints therefore contain at most $b$ endpoints in $V_2$. Again by (i), the part $V_2$ contains at least $k/2$ endpoints. The set $Z$ contains at most $b$ of them, and at most $b$ lie on the components of the chosen endpoints. As $k>6b$, at least $k/2-2b>b$ endpoints in $V_2\setminus Z$ belong to other components. Since $G_2$ is $b$-bipartite-joined, it has an edge $e$ between these endpoints and the chosen ones. Since $e$ avoids $Z$, the set $\{e\}$ is $(2,b,H)$-bounded. Adding $e$ joins two components, and hence $\Lambda$ decreases by at most $\log N$ by Lemma~\ref{lem:log-orders}, and $Q$ decreases by at least $39$.
\end{proof}

\begin{lemma}\label{lem:no-isolated}
Let $J$ be a bipartite graph on $V_1\cup V_2$ containing $G_1$, and let $\mathcal F$ minimize $Q$ among its spanning linear forests. Suppose $k(\mathcal F)\leq6b$ and $G_1$ is $(r,D)$-expanding. Then $\mathcal F$ has no isolated vertices.
\end{lemma}
\begin{proof}
We follow the proof of Lemma~4.4 of~\cite{C-expander}. We show that rotations from an isolated vertex keep $Q$ and the endpoint set, and that they then reach more endpoints than $\mathcal F$ has. Suppose $v$ is isolated. We first show that every neighbour of $v$ is the middle vertex of a path of order $3$. Since $G_1$ is $(r,D)$-expanding, $v$ has a neighbour in $J$. An edge of $J$ from $v$ to an endpoint of another component would decrease $Q$ by at least $39$ by Lemma~\ref{lem:log-orders}, contrary to minimality. Hence every neighbour $z$ of $v$ is internal, and its path has branches at $z$ of orders $a\leq c$. We add $vz$ and delete the edge between $z$ and the branch of order $c$. Then the orders $1$ and $a+c+1$ of $\{v\}$ and this path become $a+2$ and $c$. For $c\geq2$ we have $(a+2)c-(a+c+1)=(a+1)(c-1)>0$, which implies that the product of the two orders grows and $\Lambda$ increases, contradicting minimality. As $1\leq a\leq c$, we obtain $a=c=1$, and the $1$-rotation with pivot $z$ replaces orders $1,3$ by $3,1$. This preserves $Q$ and the endpoint set, and the new endpoint is isolated.

Iterating this argument, we find that every forest reached from $v$ by a rotation with added edges from $G_1$ minimizes $Q$, has the new endpoint of the rotation as an isolated vertex, and has endpoint set $End(\mathcal F)$. Lemma~\ref{lem:endpoint-growth} with $U=V_1\cup V_2$ shows that $|E^K_{V_1\cup V_2}(v,\mathcal F)|\geq20b$. All these new endpoints are endpoints of $\mathcal F$, but $|End(\mathcal F)|\leq2k(\mathcal F)\leq12b$, a contradiction.
\end{proof}

\begin{lemma}\label{lem:long-paths}
Let $J$ be a bipartite graph on $V_1\cup V_2$ containing $G_1$, and let $\mathcal F$ minimize $Q$ among its spanning linear forests. Suppose $n^{0.8}<k(\mathcal F)\leq6b$, $G_1$ is $(r,D)$-expanding, and paths of order at least $\sqrt n$ cover at least $1000b$ vertices of $\mathcal F$. Then $G_2$ has an edge $e$ such that $\{e\}$ is $(2,b,H)$-bounded, and $G_1\cup\mathcal F\cup\{e\}$ has a spanning linear forest $\mathcal F'$ with $Q(\mathcal F')\leq Q(\mathcal F)-1/5$.
\end{lemma}
\begin{proof}
We rotate from an endpoint of a path of $\mathcal F$ of minimum order to a vertex far from the ends of a long path. The component of the new endpoint then has order at least $\sqrt n/5$ instead of less than $2n^{0.2}$. This gain in $\Lambda$ is larger than the loss of at most $(K+1)\log2$ from the steps, so $Q$ decreases.

\paragraph{Rotations from $u$.} The forest $\mathcal F$ has no isolated vertices by Lemma~\ref{lem:no-isolated}, as Lemma~\ref{lem:moving} needs. Let $u$ be an endpoint of a path of $\mathcal F$ of minimum order $a_0\leq2n/k(\mathcal F)<2n^{0.2}$. Consider a rotation from $u$ with at most $K+1$ steps and added edges from $G_1$. Let $\mathcal F_i$ be the forest after $i$ steps, and let $a_i$ be the order of the component of $\mathcal F_i$ that contains the new endpoint of the $i$th step, or $u$ if $i=0$. 
\begin{claim}\label{cl:no-endpoint-pivot}
No pivot of the steps of such a rotation is an endpoint of another component, and $a_i\leq2^ia_0$ for all $i$. In particular $a_i\leq n^{1/3}$ for $i\leq K$.
\end{claim}
\begin{claimproof}[Proof of Claim~\ref{cl:no-endpoint-pivot}]
If some step has such a pivot, take the first one. The at most $K$ steps before it preserve the number of components. They decrease $\Lambda$ by at most $\log N+K\log2\leq2\log N$, by~\eqref{eq:log-orders} with $1\leq a_i\leq N$. Adding the edge of this step without deleting an edge joins two components and decreases $\Lambda$ by at most another $\log N$. The resulting forest is contained in $J$ and has one component fewer, so its value of $Q$ is at most $Q(\mathcal F)-40+3=Q(\mathcal F)-37$, contrary to minimality. Hence each $\mathcal F_i$ is a subgraph of $J$ and has $k(\mathcal F)$ components. By minimality $\Lambda(\mathcal F_i)\leq\Lambda(\mathcal F)$, and~\eqref{eq:log-orders} yields $a_i\leq2^ia_0$. Recall that $K=\lceil2\log_DN\rceil+2$ and $D=2^{100}$, which implies $2^{K+1}\leq16N^{0.02}$ and $a_i\leq2^{K+1}n^{0.2}\leq n^{1/3}$ for $i\leq K$.
\end{claimproof}

\begin{claim}\label{cl:unchanged-subpath}
For $i\leq K$, after $i$ steps of such a rotation, each original path $P$ of order at least $\sqrt n$ has an unchanged subpath in one current component, obtained from $P$ by deleting at most $i\,n^{1/3}$ vertices in total from its ends.
\end{claim}
\begin{claimproof}[Proof of Claim~\ref{cl:unchanged-subpath}]
We use induction on $i$. After $i<K$ steps, the subpath has at least $\sqrt n-i\,n^{1/3}>n^{1/3}$ vertices, as $(K+1)n^{1/3}=o(\sqrt n)$. Since the component of the new endpoint has order $a_i\leq n^{1/3}$ by Claim~\ref{cl:no-endpoint-pivot}, it does not contain the subpath. Hence, if the next broken edge is on the subpath, its pivot lies on a different component from the old endpoint. Then this edge splits the subpath into two pieces, and each piece contains one end of the subpath. The piece in the component of the next new endpoint has at most $n^{1/3}$ vertices, since this component has order at most $n^{1/3}$. We remove this piece from the subpath. The other piece is still one path. If the next broken edge is not on the subpath, the subpath stays in one component.
\end{claimproof}

\paragraph{The set $Y$.} For each path $P$ of $\mathcal F$ of order at least $\sqrt n$, take the vertices at distance at least $|V(P)|/4$ from both ends of $P$. Let $Y$ be the set of all these vertices. On each such $P$, these vertices form a subpath with at least $|V(P)|/2-2$ vertices, and each part contains at least $|V(P)|/4-2$ of them. The at most $2\sqrt n$ such paths cover at least $1000b$ vertices, so the number of vertices of $Y$ in each part is at least
\[
 \sum_P\bigl(|V(P)|/4-2\bigr)\geq1000b/4-4\sqrt n\geq100b,
\]
as $b\geq\beta n$ and $n$ is large.

The set $Y$ contains no endpoint of $\mathcal F$. Since only the two ends of each of the at most $2\sqrt n$ subpaths that form $Y$ are not in $int_{\mathcal F}(Y)$, the set $int_{\mathcal F}(Y)$ has at least $100b-4\sqrt n\geq50b$ vertices in each part. This allows us to apply Lemma~\ref{lem:moving} with $U=V_1\cup V_2$ and $V=Y$, as $G_1$ is $(r,D)$-expanding. We obtain a $(2,b,H)$-bounded set $E\subset E(G_2)$ with $|E|\leq1$ and a rotation from $u$ with at most $K+1$ steps, whose added edges are from $G_1\cup E$ and whose new endpoint $v$ is in $Y$.

\begin{figure}[htbp]
  \centering
  \begin{tikzpicture}[x=1cm,y=1cm,
      v/.style={circle,fill=black,inner sep=1pt},
      lab/.style={font=\scriptsize,inner sep=2.5pt},
      pa/.style={DarkDesaturatedBlue,line width=1.4pt}]
    \draw[gray,densely dotted] (0,0)--(1,0);
    \draw[gray,densely dotted] (9,0)--(10,0);
    \draw[line width=.8pt] (1,0)--(4.6,0);
    \draw[line width=.8pt] (5.2,0)--(9,0);
    \draw[densely dashed] (4.6,0)--(5.2,0);
    \node[v] (v) at (4.6,0) {}; \node[lab,below] at (v) {$v$};
    \node[v] (z) at (5.2,0) {}; \node[lab,below] at (z) {$z$};
    \node[v] (w) at (5.2,1.2) {}; \node[lab,above] at (w) {$w$};
    \draw[line width=.8pt] (w)--(6.4,1.2);
    \node[lab,right] at (6.4,1.2) {component of $w$ before the step, order at most $n^{1/3}$};
    \draw[pa] (w)--(z);
    \draw[|-|] (2.5,-.6)--(7.5,-.6) node[lab,midway,below] {$Y\cap V(P)$};
    \draw[|-|] (1,.35)--(4.6,.35) node[lab,midway,above] {in the new component of $v$};
    \node[lab,left] at (0,0) {$P$};
  \end{tikzpicture}
  \caption{The last step of the rotation in the proof of Lemma~\ref{lem:long-paths}. The blue edge $wz$ is added and the dashed edge $zv$ of $P$ is broken. The dotted ends of $P$ may have changed.}
  \label{fig:last-step}
\end{figure}

\paragraph{The component of $v$.} Let $\mathcal F'$ be the resulting forest. We first show that the component of $v$ in $\mathcal F'$ has order at least $\sqrt n/5$ (see Figure~\ref{fig:last-step}). Let $w$ and $z$ be the old endpoint and the pivot of the last step. Since $v\in Y$, it lies on a path $P$ of $\mathcal F$ of order at least $\sqrt n$, at distance at least $|V(P)|/4$ from both ends of $P$. Since the at most $K$ steps before the last step add edges of $G_1$ by Lemma~\ref{lem:moving}(ii), Claim~\ref{cl:unchanged-subpath} shows that the vertices at distance at least $|V(P)|/4-1$ from both ends of $P$ belong to the unchanged subpath of $P$. These include $v$ and both its neighbours on $P$. Thus the broken edge $zv$ is an edge of $P$, and $z$ is at distance at least $|V(P)|/4-1$ from both ends of $P$. As the unchanged subpath lacks at most $Kn^{1/3}$ vertices at the ends of $P$, both branches at $z$ in the forest before the last step have order at least
\[
 |V(P)|/4-Kn^{1/3}-3\geq\sqrt n/5.
\]
In particular the component of $z$ has more than $n^{1/3}$ vertices. By Claim~\ref{cl:no-endpoint-pivot}, the component of $w$ before the last step has order at most $n^{1/3}$ and cannot contain $z$. Hence $z$ is an internal vertex of another path, and $v$ lies on the shorter branch at $z$ by Definition~\ref{def:rotation}. This branch is the component of $v$ in $\mathcal F'$, and its order is $c\geq\sqrt n/5$. 

\paragraph{The gain in $\Lambda$.} The pivots of the steps before the last are not endpoints of other components by Claim~\ref{cl:no-endpoint-pivot}, and the last pivot $z$ is internal. Hence the rotation preserves the number of components by Lemma~\ref{lem:log-orders}, and~\eqref{eq:log-orders} gives
\[
 \Lambda(\mathcal F')-\Lambda(\mathcal F)\geq\log(c/a_0)-(K+1)\log2
 \geq\tfrac3{10}\log n-\log10-(K+1)\log2\geq\tfrac15\log N.
\]
The last inequality holds as $2^{K+1}\leq16N^{0.02}$ by the proof of Claim~\ref{cl:no-endpoint-pivot}. As $\mathcal F'$ has $k(\mathcal F)$ components, $Q(\mathcal F')\leq Q(\mathcal F)-1/5$. The minimality of $Q$ implies that the forest $\mathcal F'$ is not a subgraph of $J$. Thus $E\neq\emptyset$, and $E=\{e\}$ for an edge $e$ as required.
\end{proof}

\begin{lemma}\label{lem:few-paths}
Suppose that $G_1$ is $(r,4D)$-expanding, and that any vertex-disjoint paths in $G_1\cup G_2$ of lengths between $1$ and $\sqrt n$ that together cover more than $2n-1000b$ vertices have a partition into good families $\mathcal H_1,\ldots,\mathcal H_5$. Let $\mathcal F$ be a spanning linear forest in $G_1$ with $k$ paths. Then there is a $(400k,b,H)$-bounded set $E\subset E(G_2)$ such that $G_1\cup E$ has a spanning linear forest with at most $n^{0.8}$ paths and no isolated vertices.
\end{lemma}
\begin{proof}
We add edges of $G_2$ only when this also decreases $Q$. Consider the pairs $(E,\mathcal F')$, where $E\subset E(G_2)$ and $\mathcal F'$ is a spanning linear forest in $G_1\cup E$, possibly with isolated vertices, such that
\[
 \Delta_b(H+E)+10\,Q(\mathcal F')\leq\Delta_b(H)+400k.
\]
The forest $\mathcal F$ forms such a pair with $E=\emptyset$, as $Q(\mathcal F)\leq40k$. Choose a pair with $Q(\mathcal F')$ minimum, and write $J=G_1\cup E$ and $k'=k(\mathcal F')$. Note that $E$ is $(400k,b,H)$-bounded, as $Q(\mathcal F')\geq0$. A spanning linear forest in $J$ with a smaller value of $Q$ would form such a pair with $E$, which implies that $\mathcal F'$ minimizes $Q$ among the spanning linear forests in $J$.

If $k'\leq6b$, then $\mathcal F'$ has no isolated vertices by Lemma~\ref{lem:no-isolated}. As $n^{0.8}\leq6b$, we are done if $k'\leq n^{0.8}$. So suppose $k'>n^{0.8}$. In each of the three cases below, a lemma applied with $H+E$ in place of $H$ gives a set of edges of $G_2$. After we add them to $E$, the graph $G_1\cup E$ has a spanning linear forest with a smaller value of $Q$. The added edges increase $\Delta_b(H+E)$ by at most ten times the decrease of $Q$, so we get a pair with a smaller value of $Q$, contrary to the choice of $(E,\mathcal F')$.

\textbf{Case 1: $k'>6b$.} By Lemma~\ref{lem:endpoint-count}(ii), $G_2$ has an edge that joins endpoints of two components of $\mathcal F'$. Adding it increases $\Delta_b(H+E)$ by at most two and decreases $Q$ by at least $39$.

\textbf{Case 2: $k'\leq6b$, and the paths of $\mathcal F'$ of order at least $\sqrt n$ cover at least $1000b$ vertices.} Then we apply Lemma~\ref{lem:long-paths}. The selected edge increases $\Delta_b(H+E)$ by at most two, and $Q$ decreases by at least $1/5$.

\textbf{Case 3: $k'\leq6b$, and the paths of $\mathcal F'$ of order at least $\sqrt n$ cover less than $1000b$ vertices.} The paths of order less than $\sqrt n$ cover more than $2n-1000b$ vertices, and none of them is a single vertex, because $\mathcal F'$ has no isolated vertices. By the hypothesis of the lemma, they have a partition into five good families. As each part has $k'$ endpoints by Lemma~\ref{lem:endpoint-count}(i), we may apply Lemma~\ref{lem:two-endpoints}(ii) to two endpoints in opposite parts. The resulting forest has $Q$ at most $Q(\mathcal F')-20$, and the selected edges increase $\Delta_b(H+E)$ by at most ten.
\end{proof}

\begin{proof}[Proof of Lemma~\ref{lemma:linearforest decomp1}]
Consider vertex-disjoint paths in $G_1\cup G_2$ of lengths between $1$ and $\sqrt n$ that cover more than $2n-1000b$ vertices. As $\beta\ll\alpha$, they cover more than $(2-\alpha^{10})n$ vertices. Therefore they have a partition into five families as in Lemma~\ref{cl: partition into five}. These families are good and $G_1$ is $(r,4D)$-expanding, by Lemma~\ref{lem:random-good} with $m=n$. Hence the lemma follows from Lemma~\ref{lem:few-paths}, applied to the given forest.
\end{proof}

\subsection{Prescribing the endpoints}

\begin{proof}[Proof of Lemma~\ref{lemma:linearforest decomp2}]
Recall that $\mathcal F=\mathcal P\cup\mathcal Q$, where the $t$ paths of $\mathcal P$ have lengths at most $n^{0.19}$ and endpoint set $A\cup B$, and $\mathcal Q$ has $q\leq n^{0.8}$ paths. The forest $\mathcal F$ has $q$ endpoints outside $A\cup B$ in each part, and we remove them one pair at a time by $q$ applications of Lemma~\ref{lem:two-endpoints}(i). The main point is to find good families before each application, among the paths of $\mathcal P$ that are still components.

We first check the hypotheses of Lemma~\ref{cl: partition into five} for $G_1[R]$. Recall that $R=V(\mathcal P)$. As $\mathcal P$ has no isolated vertices and each part contains $t$ endpoints of $\mathcal P$, we have $e(\mathcal P)=2|R\cap V_c|-t$ for $c=1,2$, and thus $R$ is balanced. Let $m_R=|R|/2=|R\cap V_1|$. Lemma~\ref{lemma:linearforest decomp2}(ii) implies that $n/10\leq m_R\leq n$ and $\delta(G_1[R])\geq\alpha d$. Also $G_1[R]$ is $(\eta,\beta,d/m_R)$-sparse by Lemma~\ref{lem:restriction-normalization}.

Before each application, we apply Lemma~\ref{cl: partition into five} in $G_1[R]$, with $m_R$ playing the role of $n$, to the paths of $\mathcal P$ that are still components of the current forest. As $n/10\leq m_R\leq n$, the families it gives are good and $G_1$ is $(r,4D)$-expanding, by Lemma~\ref{lem:random-good}. Each application changes at most $10\log n$ edges by Lemma~\ref{lem:two-endpoints}(i), and each such edge is incident with at most two paths of $\mathcal P$. Hence, during the at most $q\leq n^{0.8}$ applications, all but at most $20q\log n$ paths of $\mathcal P$ remain components, and the other paths have at most $20q\log n\,(n^{0.19}+1)=o(n)$ vertices. The paths that remain components have positive lengths at most $n^{0.19}\leq\sqrt{m_R}$, and by this count they cover at least $2m_R-o(n)\geq(2-\alpha^{10})m_R$ vertices.

Next we remove the endpoints. Let $\mathcal F_i$ be the forest after $i$ applications, and let $k_i$ be its number of components, so that $k_0=t+q\leq n^{0.95}+n^{0.8}\leq6b$. Suppose that $\mathcal F_i$ has no isolated vertices and $k_i>t$. By Lemma~\ref{lem:endpoint-count}(i), it has exactly $k_i$ endpoints in each part. As $A\subset V_1$ and $B\subset V_2$ are sets of $t$ endpoints, both parts contain $k_i-t$ endpoints outside $A\cup B$. We choose one in each part. As $k_i\leq k_0\leq6b$, we may apply Lemma~\ref{lem:two-endpoints}(i) to them, with $H$ together with all previously selected edges in place of $H$. The new forest has no isolated vertices, and its endpoint set is the previous one without this pair. After $q$ applications the endpoint set is $A\cup B$. 

In each application, the edges of $G_2$ that we add increase $\Delta_b$ by at most ten. The union $E$ of these edges is therefore $(10q,b,H)$-bounded by Lemma~\ref{lem:max-degree}(iv). The final forest is contained in $G_1\cup\mathcal F\cup E$.
\end{proof}

\section{Hamilton connectivity of bipartite expanders}\label{app:bipartite-connectivity}
We prove Theorem~\ref{thm:bipartite-c-expander} by following the proof of Theorem~7.1 of~\cite{C-expander}. We use the rotation lemmas of Appendices~\ref{app:rotations} and~\ref{app:forests} with $G_1=G_2=G$ and with $H$ the empty graph on $V(G)$.

\subsection{Setup and expanding sets}
Let $C$ be a sufficiently large constant. In Lemma~\ref{lem:expanding-sets}, $G$ is a bipartite graph whose parts $V_1,V_2$ have size $n$, and $N=2n$ is large in terms of $C$. In Lemma~\ref{lem:expanding-sets} and in the lemmas of Appendices~\ref{app:rotations} and~\ref{app:forests}, we take $D=2^{100}$, $b=\beta n$ and $r=1000b/D$ for some $\beta\in[1/(2C),2/C]$. We assume that $G$ is $(r,4D)$-expanding and $b$-bipartite-joined. Note that, by joinedness, every set of at least $b$ vertices in one part has neighbours in all but at most $b-1$ vertices of the other part.

The next lemma plays for $G$ the role that Lemma~\ref{cl: partition into five} plays for sparse graphs. Its proof combines a bipartite version of Lemma~2.5 with the choice of five families in the proof of Lemma~3.7, both in~\cite{C-expander}.

\begin{lemma}\label{lem:expanding-sets}
Let $P_1,\ldots,P_s$ be vertex-disjoint paths in $G$ of lengths between $1$ and $\sqrt n$ that together cover at least $n/5$ vertices. Then $\{P_1,\ldots,P_s\}$ has a partition into good families $\mathcal H_1,\ldots,\mathcal H_5$, with $G$ playing the role of $G_1$.
\end{lemma}
\begin{proof}
Since each path has at most $\sqrt n+1$ vertices, we can choose four disjoint families $\mathcal H_1,\ldots,\mathcal H_4$ of the paths, each covering between $n/25$ and $n/25+\sqrt n+1$ vertices. The remaining paths form $\mathcal H_5$ and cover at least $n/5-4(n/25+\sqrt n+1)\geq n/30$ vertices.

We now check Definition~\ref{def:good}. Let sets $W_i\subset V(\mathcal H_i)$ with $|W_i|\leq20b$ be given, and let $\mathcal F$ be the linear forest formed by all the paths. Let $X_i=V(\mathcal H_i)\setminus W_i$. Since the set $V(\mathcal H_i)$ is a union of components of $\mathcal F$, we have $int_{\mathcal F}(V(\mathcal H_i))=V(\mathcal H_i)$. Each of the paths has at least a third of its vertices in each part, and deleting $W_i$ removes at most $3|W_i|$ vertices from this interior. Thus $int_{\mathcal F}(X_i)$ has at least $n/90-60b\geq n/500$ vertices in each part, as $\beta\leq2/C$ and $C$ is large.

We remove from each $X_i$ small sets with few neighbours in the interior of what remains of some $X_j$, and we take as many vertices as possible. Then no such set remains, as it could be added to the removed ones. For $i,j\in[5]$, choose sets $R_{ij}\subset X_i$ of size at most $4b$, with union $R$, that maximize the sum of their sizes subject to
\[
 |\Gamma_G(R_{ij})\cap int_{\mathcal F}(X_j\setminus R)|\leq5D|R_{ij}|
\]
for every pair $(i,j)$. The empty sets satisfy this, so such a choice exists.

Each set $R_{ij}$ has size less than $2b$. Indeed, $R$ contains at most $20b$ vertices of each $X_j$, and thus deleting $R$ removes at most $60b$ vertices from $int_{\mathcal F}(X_j)$. If $|R_{ij}|\geq2b$, then $R_{ij}$ has at least $b$ vertices in one part, so by joinedness it has at least $n/500-61b>20Db\geq5D|R_{ij}|$ neighbours in $int_{\mathcal F}(X_j\setminus R)$, a contradiction.

Set $U_i=X_i\setminus R$. We show that $|\Gamma_G(T)\cap int_{\mathcal F}(U_j)|\geq5D|T|$ for all $T\subset U_i$ with $|T|\leq b$ and all $j\in[5]$. Suppose that some such $T$ satisfies $|\Gamma_G(T)\cap int_{\mathcal F}(U_j)|<5D|T|$. We add $T$ to $R_{ij}$. As $T\subset X_i\setminus R$, this increases the sum of the sizes. The new set has size less than $2b+b\leq4b$. As $U_j=X_j\setminus R$, the neighbourhood of the new set in $int_{\mathcal F}(X_j\setminus R)$ has at most $5D|R_{ij}|+5D|T|$ vertices. Since $R$ grows, each interior $int_{\mathcal F}(X_{i'}\setminus R)$ can only shrink. Thus all the constraints still hold, contrary to maximality.

For every such $T$ and $j$, at most $|T|$ vertices of $\Gamma_G(T)\cap int_{\mathcal F}(U_j)$ belong to $T$, so at least $(5D-1)|T|\geq4D|T|$ lie outside $T$. Hence $U_i$ is $(b,4D)$-expanding into $int_{\mathcal F}(U_j)=int_{\mathcal H_j}(U_j)$, and so $(r,4D)$-expanding into it. This is property (i) of good families. Finally $|X_i\setminus U_i|<10b$, so $int_{\mathcal H_i}(U_i)$ has at least $n/500-30b\geq100b$ vertices in each part, which is (ii).
\end{proof}

\subsection{The linking structure and the Hamilton path}

\begin{proof}[Proof of Theorem~\ref{thm:bipartite-c-expander}]
Let $G$ be a bipartite $C$-expander whose parts $V_1,V_2$ have size $n$, and fix $x\in V_1$ and $y\in V_2$. We build a linking structure $\mathcal L$, cover the remaining vertices by a linear forest whose endpoints are terminals of $\mathcal L$, and join $\mathcal L$ and the forest by Lemma~\ref{lem:cycle-assembly}. In~\cite{C-expander}, the edge $xy$ is kept on a forest path and no rotation breaks it. Instead, we make $x$ and $y$ terminals of $\mathcal L$, so they are not on any forest path and no rotation uses them.

We first assume that $n$ is large in terms of $C$, and we treat smaller orders at the end. Let $D_L=8D$ and $\rho=1/D_L$, and fix $0<\delta<1/10$ with $\delta\log(e/\delta)\leq\rho/10^4$. Set
\[
 b=n/(2C),\quad r=1000b/D,\quad b_L=5n/C,\quad r_L=1000b_L/D_L,\quad D'=\lfloor C/1000\rfloor,\quad t=n^{0.9},
\]
and let $\ell_0$ be the largest odd integer with $\ell_0\leq n/(2t)$. We use $b$ and $r$ in Lemma~\ref{lem:expanding-sets} and in Appendices~\ref{app:rotations} and~\ref{app:forests}, $D_L,\rho,\delta,b_L,r_L$ in Lemma~\ref{lem: partition}, and $D'$ and $b_L$ when we grow the path $S_0$.

We grow a path $S_0$ from the edge $xy$ by the extendability method, cut $S_0-\{x,y\}$ into $t$ short paths and apply Lemma~\ref{lem: partition} to these paths. Extendability is with respect to $G$, and the edge $xy$ does not have to be an edge of $G$, as the definition of extendability allows. The graph consisting of the edge $xy$ is $(D',b_L)$-extendable. Indeed, a nonempty set $U$ in one part with $|U|<b$ has at least $C|U|-2\geq D'|U|$ neighbours outside $\{x,y\}$. If $b\leq|U|\leq2b_L$, then by joinedness $U$ has at least $n-b\geq2D'b_L$ such neighbours. In both cases this is at least the right side of~\eqref{eq:extendability}.

We extend the edge $xy$ beyond $x$ by adding leaves (Lemma~\ref{lem:add-leaf}) to a path $S_0$ with $t(\ell_0+1)+2$ vertices. This is possible, as each graph we build has at most $n/3$ vertices in each part. By joinedness, every set $U$ in one part with $b_L\leq|U|\leq2b_L$ then has more than $n-b-n/3>2D'b_L+1$ neighbours outside that graph.

We cut $S_0-\{x,y\}$ into consecutive subpaths of $\ell_0+1$ vertices and obtain a linear forest $\mathcal P^0$ of $t$ paths of length $\ell_0\leq n^{0.15}$. Let $U_1=V(\mathcal P^0)$ and $U_2=V(G)\setminus U_1$. As $\ell_0$ is odd, $U_1$ and $U_2$ are balanced. Also $x,y\in U_2$ and $n/2-t\leq|U_1|\leq n/2+t$.

By Lemma~\ref{lem:unused-expansion} applied to $S_0$, every set $U$ with $|U|\leq2b_L$ has at least $(D'-4)|U|$ neighbours outside $V(S_0)$, and so in $U_2$. We apply Lemma~\ref{lem: partition} with $G_1=G_2=G$, $H$ empty, the forest $\mathcal P^0$, $x_1=x$, $x_2=y$ and the constants $D_L,\rho,5/C,\delta$. Condition~\ref{eq: partition condition 1} of Lemma~\ref{lem: partition} holds, as $r_L\leq b_L$ and, since $C$ is large, $D'-4\geq2D_L$ and $(D'-4)r_L\geq Cb_L/(2D_L)>\rho n$. The other conditions hold, as $G$ is $b_L$-bipartite-joined, $t\leq\delta n$ and $5/C\leq\rho/10^8$.

Let $F,X,Y,Z,A,B$ be given by the lemma, and let $\mathcal L=F[X]$ and $\mathcal P=F[U_1\cup Y\cup A\cup B]$. By~\ref{eq: partition 2} and~\ref{eq: partition 3}, $\mathcal L$ is an $(A\cup\{x\},B\cup\{y\})$-linking structure with $x,y\notin A\cup B$ and $|A|=t$, and $\mathcal P$ consists of $t$ paths of length between $1$ and $3n^{0.15}$ with endpoint set $A\cup B$.

\paragraph{The forest outside $\mathcal L$.}
Let $W=U_1\cup Y\cup Z\cup A\cup B$, so that $W\cap X=A\cup B$. We cover $W$ by a linear forest with endpoint set $A\cup B$ in two steps. First we cover $Z$ by a small number of paths using Lemma~\ref{lem:few-paths}, and then we remove their endpoints one pair at a time using Lemma~\ref{lem:two-endpoints}(i).

We check that $G[Z]$ and $G[W]$ satisfy the assumptions made on $G$ before Lemma~\ref{lem:expanding-sets}. The sets $Z$ and $W$ are balanced. By~\ref{eq: partition 1} and $|U_1|\leq n/2+t$, each of $G[Z]$ and $G[W]$ has parts of some size $m\geq n-n/4-t-\rho n\geq n/2$, so $b/m\in[1/(2C),2/C]$ can play the role of $\beta$. Both graphs are $b$-bipartite-joined as $G$ is. They are also $(r,4D)$-expanding by~\ref{eq: partition 4}, as $b_L\geq5b$ gives $2r_L\geq r$, and $D_L/2=4D$.

We apply Lemma~\ref{lem:few-paths} with $G_1=G_2=G[Z]$ and any spanning linear forest of $G[Z]$. It gives a spanning linear forest $\mathcal Q$ of $G[Z]$ with $q\leq n^{0.8}$ paths and no isolated vertices. Its hypothesis on five families holds by Lemma~\ref{lem:expanding-sets}, as $|Z|-1000b\geq|Z|/10$.

The forests $\mathcal P$ and $\mathcal Q$ together form a spanning linear forest on $W$ with no isolated vertices and $t+q\leq6b$ components. As in the proof of Lemma~\ref{lemma:linearforest decomp2}, we make $q$ applications of Lemma~\ref{lem:two-endpoints}(i) in $G[W]$, each removing two endpoints outside $A\cup B$, one in each part.

Each application uses the five families that Lemma~\ref{lem:expanding-sets} gives for the paths of $\mathcal P$ that are still components. The lemma applies, as these paths have length at most $3n^{0.15}\leq\sqrt{|W|/2}$ and cover more than $|W|/10$ vertices. For the latter, the count in the proof of Lemma~\ref{lemma:linearforest decomp2}, with $3n^{0.15}$ in place of $n^{0.19}$, shows that the paths that are still components cover at least $|U_1|-o(n)\geq n/2-o(n)>|W|/10$ vertices. After these applications, we have a spanning linear forest on $W$ with no isolated vertices, and its endpoint set is exactly $A\cup B$.

Finally, as $X\cup W=V(G)$ and $X\cap W=A\cup B$, Lemma~\ref{lem:cycle-assembly}, applied to $\mathcal L$ and the final forest, gives a spanning $x,y$-path of $G$.

\paragraph{Small orders.}
Fix $C_0$ large enough for the preceding proof, and choose $N_0$ so that the proof applies to every bipartite $C_0$-expander of order at least $N_0$. Take $C\geq\max\{C_0,N_0\}$. A bipartite $C$-expander is also a bipartite $C_0$-expander. Indeed, by joinedness, a set $U$ in one part with $n/(2C)<|U|\leq n/(2C_0)$ has more than $n-n/(2C)\geq n/2\geq C_0|U|$ neighbours. If $2n\geq N_0$, the preceding proof applies. If $2n<N_0$, then $n/(2C)<1$, so $G$ is complete bipartite by joinedness. We take the vertices of the two parts alternately, starting at $x$ and ending at $y$, and obtain the desired Hamilton path.\end{proof}

\end{document}